\documentclass{article} 
\usepackage{graphicx}
\usepackage{amssymb}
\usepackage{amsthm} 
\usepackage{amsmath}
\usepackage{graphics} 
\usepackage{eepic}    
\usepackage[all,arc,poly,cmtip]{xy}

\newcommand{\nc}{\newcommand}
\nc{\look}{\marginpar{$\bullet$}}

\nc{\Section}{\section}
\nc{\SubSection}{\subsection}

\newtheorem{theo}{Theorem}[section]   
\newtheorem{ddef}[theo]{Definition}
\newtheorem{llem}[theo]{Lemma} 
\newtheorem{oobs}[theo]{Observation} 
\newtheorem{rrem}[theo]{Remark} 
\newtheorem{prop}[theo]{Proposition} 
\newtheorem{ccor}[theo]{Corollary}  
\newtheorem{claim}[theo]{Claim} 
\newtheorem{qquest}[theo]{Question} 
\newtheorem{fact}[theo]{Fact} 
\newtheorem{pprov}[theo]{Proviso}
\newtheorem{eexam}[theo]{Example} 

\nc{\bT}{\begin{theo}} 
\nc{\eT}{\end{theo}}
\nc{\bD}{\begin{ddef} \rm }
\nc{\eD}{\end{ddef}}
\nc{\bC}{\begin{ccor}}
\nc{\eC}{\end{ccor}}
\nc{\bCl}{\begin{claim}}
\nc{\eCl}{\end{claim}}
\nc{\bQ}{\begin{qquest}}
\nc{\eQ}{\end{qquest}}
\nc{\bL}{\begin{llem}}
\nc{\eL}{\end{llem}}
\nc{\bP}{\begin{prop}}
\nc{\eP}{\end{prop}}
\nc{\bR}{\begin{rrem}}
\nc{\eR}{\end{rrem}}
\nc{\bO}{\begin{oobs}}
\nc{\eO}{\end{oobs}}
\nc{\bF}{\begin{fact}}
\nc{\eF}{\end{fact}}
\nc{\bProv}{\begin{pprov}}
\nc{\eProv}{\end{pprov}}
\nc{\bE}{\begin{eexam} \rm }
\nc{\eE}{\end{eexam}}

\nc{\prf}{\begin{proof}}
\nc{\eprf}{\end{proof}}

\renewcommand{\phi}{\varphi}
\renewcommand{\geq}{\geqslant}
\renewcommand{\leq}{\leqslant}
\renewcommand{\preceq}{\preccurlyeq}
\renewcommand{\succeq}{\succcurlyeq}
\newcommand{\strictsubset}{\varsubsetneq}
\renewcommand{\subset}{\subseteq}
\newcommand{\subsetw}{\subset_{\ssc w}}

\renewcommand{\supset}{\supseteq}

\newenvironment{romanenumerate}%
{\begin{list}{(\roman{enumi})}{\usecounter{enumi}
\setlength{\labelwidth}{2cm}
\setlength{\itemindent}{0pt}
\setlength{\itemsep}{0.5\itemsep}
\setlength{\topsep}{\itemsep}
\setlength{\parsep}{0pt}
}}{\end{list}}
\nc{\bre}{\begin{romanenumerate}}
\nc{\ere}{\end{romanenumerate}}
\newenvironment{alphaenumerate}%
{\begin{list}{(\alph{enumii})}{\usecounter{enumii}
\setlength{\labelwidth}{2cm}
\setlength{\itemindent}{0pt}
\setlength{\itemsep}{0.5\itemsep}
\setlength{\topsep}{\itemsep}
\setlength{\parsep}{0pt}
}}{\end{list}}
\nc{\bae}{\begin{alphaenumerate}}
\nc{\eae}{\end{alphaenumerate}}
\newenvironment{numenumerate}%
{\begin{list}{(\arabic{enumiii})}{\usecounter{enumiii}
\setlength{\labelwidth}{2cm}
\setlength{\itemindent}{0pt}
\setlength{\itemsep}{0.5\itemsep}
\setlength{\topsep}{\itemsep}
\setlength{\parsep}{0pt}
}}{\end{list}}
\nc{\bne}{\begin{numenumerate}}
\nc{\ene}{\end{numenumerate}}

\nc{\ins}[1]{\bigskip\noindent
\framebox{\begin{minipage}{.95\textwidth} \sloppy \noindent \em #1 \end{minipage}}\bigskip}

\nc{\str}[1]{{\mathfrak{#1}}}
\nc{\brck}[1]{[\![ #1 ]\!]}
\nc{\restr}{\!\restriction\!}
\nc{\G}{\mathbb{G}}
\nc{\CG}{\mathbb{CG}}
\nc{\DG}{\mathbb{DG}}
\nc{\CGhat}{\mathbb{C}\hat{\mathbb{G}}}
\nc{\Ghat}{\hat{\mathbb{G}}}
\nc{\K}{\mathbb{K}}
\nc{\Khat}{\hat{\mathbb{K}}}
\nc{\Hhat}{\hat{\mathbb{H}}}
\nc{\Ehat}{\hat{E}}
\nc{\Shat}{\hat{S}}
\nc{\Vhat}{\hat{\mathbb{V}}}
\nc{\EEhat}{\hat{\E}}
\nc{\Ihat}{\hat{\mathbb{I}}}
\nc{\T}{\mathbb{T}}
\nc{\CE}{\mathsf{CE}}

\nc{\sym}{\mathrm{sym}}

\nc{\dom}{\mathrm{dom}}

\newcommand{\E}{\mathsf{E}}
\nc{\F}{\mathsf{F}}
\renewcommand{\H}{\mathbb{H}}
\newcommand{\V}{\mathbb{V}}
\newcommand{\I}{\mathbb{I}}

\nc{\barr}{\begin{array}}
\nc{\earr}{\end{array}}
\nc{\btab}{\begin{tabular}}
\nc{\etab}{\end{tabular}}

\nc{\nothing}{\rule{0em}{1ex}}
\nc{\highnothing}{\rule{0em}{3ex}}
\nc{\hnt}{\highnothing}
\nc{\nt}{\nothing}
\nc{\nnt}{\rule{.1pt}{0pt}}

\nc{\ssc}{\scriptscriptstyle}
\nc{\N}{{\mathbb N}}
\nc{\Z}{{\mathbb Z}}
\nc{\cym}{\mathbb{G}}

\renewcommand{\epsilon}{\varepsilon}

\begin{document}

\title{Acyclicity in Finite Groups and Groupoids}
\author{Martin Otto\thanks{Research partially supported 
by DFG grant OT~147/6-1: \emph{Constructions and Analysis in
  Hypergraphs of Controlled Acyclicity.}}
\\
Department of Mathematics\\
Technische Universit\"at Darmstadt}
\date{February~2024}

\maketitle

\begin{abstract}
\noindent
We expound a concise construction of finite groups and groupoids 
whose Cayley graphs satisfy graded acyclicity requirements.
Our acyclicity criteria concern cyclic patterns formed by coset-like
configurations w.r.t.\ subsets of the generator set 
rather than just by individual generators. The proposed constructions 
correspondingly yield finite groups and groupoids whose Cayley 
graphs satisfy much stronger acyclicity conditions than large girth. 
We thus obtain generic and canonical constructions of highly
homogeneous graph structures with strong acyclicity properties, 
which support known applications in finite graph and hypergraph 
coverings that locally unfold cyclic configurations. 
\end{abstract}

\pagebreak 

\tableofcontents 

\pagebreak

\section{Introduction}

The intimate connection between finite groups and graph-like structures 
is a long-standing theme that illustrates core concepts at the interface
of algebra and discrete mathematics. Groups arise as automorphism
groups of structures, and Frucht's theorem~\cite{Frucht} says that
every finite group arises as an automorphism group of a finite
graph; in particular, the given finite group -- an abstract group --
is realised as a permutation group, and thus as a subgroup of the full
symmetric group of some finite set, and in fact even as the full group of all
symmetries of a specifically designed discrete structure of a very
simple format. 

At a very basic level, permutation group actions 
can be defined through generators whose operation can be traced in
graph-like structures, which in turn determine the abstract group
structure~\cite{Cayley1,Cayley2}. The key notion in this correspondence is the
representation of the algebraic structure of the given group in its
Cayley graph: an edge-coloured directed graph that represents the 
internal group action of a chosen set of generators for the group. 

Interesting finite groups can be obtained from permutation
group actions induced by
graph-like extra structure on a finite set,
also in other ways than just as a group of symmetries. Specific graph
structures and carefully designed permutation group actions
can thus give rise to finite groups with desirable algebraic or combinatorial 
properties suggested by various applications. A very nice example of
this technique is a construction, due to Biggs~\cite{Biggs89} and
outlined in~\cite{Alon95}, 
of finite groups over a given set of generators that avoid short
cycles, i.e.\ in which non-trivial products of a small number of generators 
cannot evaluate to the neutral element. In terms of the Cayley graph of the
resulting group one obtains finite graphs of large girth that are not
only regular but (like any Cayley graph) highly symmetric in the
stronger sense of possessing a transitive automorphism group. 

Acyclicity criteria for groups matter in many natural applications. 
The free group over a given set of generators, which can be seen as
the unique fully acyclic group structure over the given generators, 
arises naturally in connection with universal coverings in the
classical topological context as well as in the context of discrete
structures, e.g.\ with tree unfoldings of transition systems. The
relevant coverings can be described as products with (the Cayley
graphs of) free groups. Of course free groups, and fully acyclic 
coverings in non-trivial settings, are necessarily infinite. 
Where finiteness matters and needs to be preserved, e.g.\ in finite
coverings, full acyclicity is typically unavailable. Here
graded degrees of acyclicity, like lower bounds on the girth
of the Cayley graph, are best possible and often can replace full
acyclicity, especially for local structural analysis -- just as a
graph of large girth is locally tree-like. Previous work, which arose
from applications in logic and the model theory of finite structures,
has led to the introduction of similar but much stronger measures
of graded acyclicity in Cayley graphs of finite groups. These notions
of acyclicity arise naturally in connection with covering
constructions for finite graphs and hypergraphs. Instead of
controlling just the length of shortest generator cycles, similar
control is achieved over the length of shortest cycles formed by
cosets w.r.t.\ generated subgroups. This generalisation involves a
passage from cycles at the level of individual generators to cycles 
formed by cosets, which a priori are not even bounded in size. 
In other words, this
is a shift in focus from first-order objects (generators) to
second-order objects (cosets) in the desired groups.
Corresponding constructions, which are inspired by Biggs' technique
but adapt the basic idea to the more complex technical setting, 
were first developed  for groups in~\cite{Otto04} and~\cite{Otto12JACM}
for specific applications of finite graph coverings. 
Generalisations of these techniques  
to the setting of groupoids offer a more direct route to hypergraph 
coverings (here necessarily branched, in a discrete analogue of
classical terminology from~\cite{Fox}). A main challenge and 
goal in these settings lies in the construction of corresponding
coverings that are generic and natural in the sense that they
do not break any symmetries of the underlying structure. 
This is essential for far-reaching applications, e.g.\ towards 
extension problems for local symmetries~\cite{lics2013,arXiv2018}. 

The goal here is a concise and generic combinatorial construction 
of groups and groupoids with strong acyclicity properties that control 
coset cycles rather than just generator cycles.
The present exposition not only serves to correct a serious mistake in 
the construction of the relevant groupoids that was sketched
in~\cite{lics2013}%
\footnote{Cf.\ acknowledgements at the end of this paper on this somewhat frayed
  history.}
but also to unify the treatment of finite groups and groupoids with the
desired acyclicity properties. 

One point in the generalisation from groups to groupoids stems from
the limitation of~\cite{Biggs89,Otto04}
to involutive generators, which does not directly fit the groupoidal setting.
In the current extension of the original idea we propose a
construction of highly acyclic finite groups, which is still based on involutive 
generators but yields stronger results for these groups;  
stronger notions of acyclicity, which are based on more general
patterns than plain coset cycles, allow us to give a self-contained
account in which the construction of highly acyclic finite groupoids 
can be reduced to the new, enriched construction
for groups with involutive generators. This yields a unified construction 
which offers a transparent view of the commonality between the two, 
seemingly so very different settings,
which may support further insights and applications.
A treatment that is essentially based on the constructions proposed
here, but geared more directly towards their algebraic uses in the
theory of finite groups and inverse semigroups has meanwhile been
given in~\cite{ABO}, 
which in contrast to the present treatment it does not rely on involutive generators. 
Concerning basic applications
we here discuss more general and more direct constructions of finite graph
and hypergraph coverings in Propositions~\ref{graphcoverprop} 
and~\ref{hypcoverprop}. 

\subsection*{Terminology and notation}

\paragraph*{Graphs and relational structures.}
In this paper we consider various kinds of graphs, some undirected, some directed,
often also allowing loops (reflexive edges), and in Section~\ref{groupoidsec}
also multi-graphs that may have more than one edge linking the same
two vertices. Notation should be standard, with small
adaptations to the specific formats that will be explicitly
stated where they occur. We mostly use a relational format for the
specification of a graph, with a binary edge relation, or with a separate edge
relation for each colour to encode edge-coloured graphs. In some
instances, and especially in Sections~\ref{groupoidsec} and~\ref{Nacycgroupoidsec}, 
it is natural to treat graphs and especially multi-graphs as 
two-sorted structures with a set of edges and a set of vertices 
linked by incidence maps that specify source and target vertices of each edge.
For \emph{subgraphs} we explicitly 
distinguish between \emph{induced subgraphs} (whose edge
relation is the restriction of the given edge relation to the
restricted set of vertices) and \emph{weak subgraphs} (whose edge
relation may be a proper subset of the given edge relation even 
in restriction to the smaller vertex set). Also more generally for relational structures
we use $\subsetw$ for the \emph{weak substructure} relationship, 
$\subset$ for the induced \emph{substructure} relationship. 
By a \emph{component} of a graph structure we mean 
an induced substructure that is closed w.r.t.\ the edge relation; 
a \emph{connected component} is a minimal  component.
The term \emph{reduct} refers to a restriction in the number of 
edge relations, or edge colours, which corresponds to the 
deletion of all edges of the colours to be eliminated.%
\footnote{Depending on context the corresponding edge relations 
can be thought of as erased (which produces a structure over a
smaller signature),  or just as emptied (which produces a weak substructure).}
 
We implicitly always assume that structures need only be identified up
to isomorphism. 
This will often allow us to avoid notational complications in the
interest of clarity. For instance, it often makes sense to suppress
explicit notation for isomorphic embeddings if instead we can w.l.o.g.\ treat the
pre-image as an actual (weak) substructure of the target structure;
relational structures of particular interest, especially Cayley graphs, will be
homogeneous in the sense that any two elements are linked by 
an automorphism, which implies that, up to isomorphism, 
explicit choice of a distinguished element can be suppressed and only
the relative position of two or more elements gives rise to meaningful distinctions.

\paragraph*{Algebraic structures.}
For structures like groups, semigroups, monoids or groupoids we 
adopt multiplicative notation 
and would typically write, for instance, $g \cdot h$ or just $gh$ 
for the result of the composition of group elements $g$ and $h$
w.r.t.\ the group operation, $1$ for the neutral
element and $g^{-1}$ for the inverse of $g$. When dealing with
subgroups of the symmetric group of some set $X$, we sometimes make 
the group operation explicit as in $h \circ g$ for the composition
of $g$ with $h$, which maps $x \in X$ to $h(g(x))$, and would in our
standard notation be rendered as $g \cdot h$ or $gh$  (!)
since we think of permutations
as operating from the right. 

Among standard terminology from other fields of 
mathematics we use some basic terms 
from formal language theory, especially to deal with \emph{words}
over a finite alphabet $E$ of letters; the set of all $E$-words is
the set of all finite (but possibly empty) strings or tuples of 
letters from $E$, denoted $E^\ast = \bigcup_{n \in \N} E^n$. 
As is common in formal language theory, we write a typical word 
of length $n \in \N$ as
$w = e_1 e_2 \cdots e_n \in E^n$ (rather than e.g., in tuple notation, as
$(e_1, e_2, \ldots, e_n)$), denoting its length as $n = |w|$.
We also write, e.g.\ just $w_1 w_2$ for the
\emph{concatenation} of the words $w_1,w_2 \in E^\ast$ (which is
often denoted as $w_1 \cdot w_2$ with explicit notation for the
concatenation operation as a monoidal semigroup operation). 
The \emph{empty word} $\lambda \in E^\ast$, which
is the unique $E$-word of length $0$, is the neutral element 
in the monoid $E^\ast$. 
Depending on the r\^ole of the letters $e \in E$, we may use
$E$-words to specify different objects of interest: 
thinking of $E$ as a set of generators of some group, an $E$-word  
is a generator word which can be read as a group product 
specifying a group element; thinking of $E$ as a set of colours in an
edge-coloured graph, an $E$-word is a colour sequence and can 
specify the class of walks that realise that colour sequence.
In some cases we also invoke a notion of \emph{reduced words}, which
are typically obtained by some cancellation operation. Especially if $E$ is a set of
generators of a group that is closed under inverses we may
(inductively) cancel factors $e e^{-1}$ in order to associate with
every $E$-word a unique reduced $E$-word that denotes the same
group element. In such contexts we often let $E^\ast$ stand for the
set of reduced words, endowed with the concatenation operation that 
implicitly post-processes plain concatenation by the necessary 
cancellation steps. More formally one could explicitly distinguish
between $E^\ast$ and its quotient $E^\ast/{\sim}$, but we suppress
this as an unnecessary distraction in our considerations.

\section{General patterns}
\label{genpatternsec}

\subsection{Cayley \& Biggs: the basic construction}
\label{Biggssec}

The fundamental idea to associate groups with 
permutation group actions and graphs can be
attributed to Arthur Cayley~\cite{Cayley1,Cayley2}. 
The Cayley graph of an abstract group, w.r.t.\ 
to a chosen set of generators, encodes the algebraic structural 
information about the algebraic group, and also represents the given 
group as a subgroup of the full symmetric group, and more
specifically as the automorphism group, of the Cayley graph. 
The natural passage between combinatorial properties of  
graph-like structures and group-like structures offers interesting
avenues for the construction of group-like and graph-like structures.
A classical example is the use of Cayley graphs in Frucht's
construction of (finite) graphs that realise a given abstract (finite)
group as their automorphism group~\cite{Frucht}.
In particular, Cayley graphs are, by construction, not just regular
but homogeneous in the sense of having a transitive 
automorphism group. So on one hand, Cayley graphs provide examples of graph 
structures with a particularly high degree of internal symmetry.
On the other hand, permutation group actions on suitably designed
graph structures generate groups that can display specific combinatorial properties
w.r.t.\ to a chosen set of generators -- and these groups in turn generate Cayley
graphs that reflect those group properties. It is one characteristic
feature of the inductive constructions to be expounded here that they are based
on a feedback loop built on this interplay. 

The idea to extract groups with certain acyclicity properties from
permutation group actions on suitably prepared graph structures is
best illustrated by the basic example of a construction of regular
graphs of high girth due to Biggs~\cite{Biggs89} and outlined
in~\cite{Alon95}. 

Let $E$ be a finite set of letters, $|E| = d \geq 2$, to be used to label involutive 
generators of a group to be constructed. With $E$ and a parameter $n
\geq 1$ in $\N$ associate a tree $\T(E,n)$ and a
group $\G(E,n)$ as follows. Let $\T(E,n)$ be a
$d$-branching, regularly $E$-coloured, finite undirected tree of depth $n$, as
represented by the set of all \emph{reduced words} $w \in E^{\leq n}
\subset E^\ast$, i.e.\ strings $w = e_1\cdots e_m$ of length
$|w|= m$, $0 \leq m \leq n$, 
with $e_i \in E$ for $1 \leq i \leq m$ and $e_{i+1} \not= e_{i}$
for $1 \leq i < m$. We regard the empty word $\lambda \in E^\ast$ 
as the root of $\T(E,n)$. More formally, we let 
\[
\T(E,n) = (V,(R_e)_{e \in E})
\]
be the tree structure with vertex set 
\[
V := \{ w \in E^\ast \colon |w| \leq n, w \mbox{ reduced } \}
\]
and undirected edge relation $R = \dot{\bigcup}_{e \in E} R_e$,
$E$-coloured by its partition into the 
\[
R_e :=  \{ (w,we), (we,w) \colon w,we \in V \}   
\] 
for $e \in E$. By construction, each vertex $w \in V$ with $|w| < n$
is an interior vertex of $\T(E,n)$ of degree $d = |E|$, with
precisely one $R_e$-neighbour for each $e \in E$; the remaining
vertices, viz.\ those $w \in V$ with $|w| = n$, are leaves of
$\T(E,n)$, each with an $R_e$-neighbour for a unique $e \in E$
(the last letter of $w$). Note that each $R_e$ is a partial
matching over $V$, and that $R_e$ and $R_{e'}$ are disjoint for $e
\not= e'$. With $e \in E$ we associate
the permutation $\pi_e \in \mathrm{Sym}(V)$ that swaps any pair of
vertices that are incident with a common $e$-coloured edge. This is
the involutive permutation of $V$ whose graph is the
matching $R_e$ augmented by loops in vertices not incident with an
$e$-coloured edge. 
The target of the construction is the group $\G(E,n)$,
which is the subgroup of
$\mathrm{Sym}(V)$ generated by these involutions:
\[
\G  = \G(E,n) := \langle \pi_e \colon e \in E \rangle  \subset \mathrm{Sym}(V).
\]

For the group operation we use the convention that the action by the 
generators is regarded as a right action via composition, i.e.\ for
$\rho \in G$:  
\[
\barr{rcl}
\rho \pi_e = \pi_e \circ \rho \colon V &\longrightarrow& V 
\\
w &\longmapsto& \pi_e(\rho(w)).
\earr
\]  

Its Cayley graph w.r.t.\ the generators $(\pi_e)_{e \in E}$ is an
edge-coloured graph $\CG$, with the set of group 
elements $\rho \in G$ as its vertex set, and with a family of 
edge relations 
\[
R_e^{G} := \{ (\rho, \rho \pi_e) \colon \rho \in G, e \in E \}
\subset G \times G,
\]
one for each $e \in E$.
Here these edge relations are symmetric due to the involutive nature of 
the $\pi_e$ in $\mathrm{Sym}(V)$, and they are irreflexive
and pairwise disjoint since $\mathrm{id}_V \not= \pi_e \not= \pi_{e'}$
for $e \not= e'$, as can be seen most easily by their action as 
permutations on $\lambda \in V$.  
So this Cayley graph is a $d$-regular finite graph, whose 
automorphism group acts transitively on the set of vertices. For the 
last claim consider the left action of the group on itself:
\[
\barr{rcl}
h \colon G &\longrightarrow& G 
\\
g &\longmapsto& hg,
\earr
\]  
which clearly induces an automorphism 
of the Cayley graph (albeit not of the group, which is rigid once we
label the generators). 
That the girth of the Cayley graph of $\G$ is at least
$4n+ 2$ can be seen as follows.
A reduced word $w \in E^k$ of length $k \geq 1$ can be written as 
$w = e_1u$. Let $v \in E^n$ be a leaf of $\T(E,n)$ whose reversal
$v^{-1}$ agrees with $u$ (up to $\mathrm{max}(n,|u|)$). 
Applying the corresponding permutation 
$\pi_w = \pi_u \circ \pi_{e_1}$ to $v$, we see that the action of the 
permutations prescribed by the first (up to) $n+1$ letters of $w$ takes
that leaf step by step towards the root $\lambda$, the next $n$
letters (if present) will take it step by step towards a different
leaf, where the very next letter (if present) can have no effect so
that it would take at least the action of another $2n$ letters after
that to bring this vertex back to where we started. In other words, no reduced word
of fewer than $n+1+n+1+ 2n = 4n +2$ letters can label a generator
sequence that represents the neutral element of the group, 
which is the identity in $\mathrm{Sym} (V)$. 

Considering what is essential for the passage from a 
graph like $\T(E,n)$ to a group like $\G(E,n)$, the only obvious
necessity is that each of the edge colours induces a partial matching
of the underlying vertex set in order to have well-defined
involutions $\pi_e$. Tree-likeness, by contrast, is of no special 
importance, not even for the bound on the girth of the resulting
group or Cayley graph. If $\T(E,n)$ were replaced,
for instance, by the disjoint union of all $E$-coloured line graphs 
corresponding to reduced words $w \in E^{2n}$, the above girth bound
of $4n+2$ persists with essentially the same argument.  
In the following paragraph we extract the basic format for the
generation of groups with involutive generators from edge-coloured
undirected graphs.

\subsection{E-graphs and E-groups}
\label{Egraphsandgroupssec}

In the following it is convenient to allow loops in the symmetric 
edge relation of an undirected graph $(V,R)$, 
and to let a loop at vertex $v$ contribute value~$1$ to the
degree of that vertex. A \emph{partial matching} is here cast as a
symmetric edge relation whose degree is bounded by $1$ at every
vertex, and may thus be thought of as the graph of a partial bijection
that is involutive (its own inverse); this involution has precisely
those vertices as fixed points at which the edge relation has loops,
and its domain $\mathrm{dom}(R)$ and range $\mathrm{rng}(R)$ 
consists of the set of the vertices of degree~$1$. A \emph{full
  matching} is a symmetric edge relation $R$ on $V$ 
such that every vertex $v \in V$ has a unique $R$-neighbour, which
in the case of a loop may be $v$ itself; it therefore corresponds to
the graph of an involutive permutation of the vertex set $V$.

\bD[$\E$-graph]
\label{Egraphdef}
\mbox{}\\
For a set $E$, an \emph{$\E$-graph} is an undirected 
edge-coloured graph $\H= (V,(R_e)_{e \in E})$ whose undirected 
edges are $E$-coloured in such a way that each $R_e$ is a partial
matching over the vertex set $V$.
The $\E$-graph $\H = (V,(R_e)_{e \in E})$ is \emph{strict} 
if there are \emph{no loops} (each $R_e$ is irreflexive) and 
\emph{no multiple edges} ($R_e \cap R_{e'} = \emptyset$ for $e \not= e'$). 
The $\E$-graph $\H = (V,(R_e)_{e \in E})$ is \emph{complete} if each $R_e$
is a full matching. 
The \emph{trivial completion} of an 
$\E$-graph $\H = (V,(R_e)_{e \in E})$ is the 
complete $\E$-graph $\bar{\H} = (V, (\bar{R}_e)_{e \in E})$ obtained by
putting $\bar{R}_e := R_e \cup \{ (v,v) \colon v \in V\setminus\mathrm{dom}(R_e)\}$.
\eD

We think of $R_e$-edges as edges of colour $e$ or as edges labelled
with $e$. In this sense an $\E$-graph is a special kind of 
$E$-coloured graph whose overall
edge relation would be $\bigcup_{e\in E} R_e$.

\medskip
For groups $\G = (G,\cdot\,,1)$ (in multiplicative notation),
an element $g \in G$ is an \emph{involution} if $g = g^{-1}$. 
A subset $E \subset G \setminus \{ 1 \}$ is a \emph{set of generators} 
for $\G$  if every group element $g \in G$ can be written as a 
product of elements from $E$ and their inverses.

\bD[$\E$-group]
\label{Egroupdef}
\mbox{}\\
For a set $E$, an \emph{$\E$-group} 
is any group $\G = (G,\cdot\,,1)$ that has $E \subset G$ as a set of
non-trivial involutive generators.%
\footnote{Clearly the elements $e \in E \subset G$ are pairwise 
distinct as elements of $\G$, and non-triviality means that $e \not= 1$.}
\eD

If $\G$ is an $\E$-group, we write $[w]_{\G} \in G$ for the group element
that is the group product of the generator sequence $w \in E^\ast$, so that 
\[
\barr{rcl}
[\,\nt\,]_{\G} \colon E^\ast &\longrightarrow& \G
\\
w = e_1\cdots e_n &\longmapsto& [w]_{\G} := \prod_{i=1}^n e_i = e_1 \cdots e_n
\earr
\]
is a surjective homomorphism from the free monoid structure of $E^\ast$,
with concatenation and neutral element $\lambda \in E^\ast$,
onto the group $\G$. 

\bO
\label{freeinvgroupobs}
The quotient of the free group generated by $E$ w.r.t.\ to the
equivalence relation induced by the identities $e = e^{-1}$ for 
$e \in E$ (as represented by reduced words in $E^\ast$) can be 
regarded as the \emph{free $\E$-group}. All other 
$\E$-groups are homomorphic images of this free $\E$-group.
\eO

\bD[$\sym(\H)$]
\label{symdef}
\mbox{}\\
For an $\E$-graph $\H= (V,(R_e)_{e \in E})$ we let $\sym(\H)$ be
the subgroup of $\mathrm{Sym}(V)$ that is generated by the involutive permutations
$\pi_e \colon V \rightarrow  V$ induced by the full matchings 
of its trivial completion $\bar{\H} = (V, (\bar{R}_e)_{e \in E})$.
\eD

Provided the $(\pi_e)_{e \in E}$ are pairwise distinct and distinct
from $\mathrm{id}_V$, we regard $\sym(\H)$ as an $\E$-group 
where we identify $e \in E$ with the generator $\pi_e$, for $e \in
E$. We shall always tacitly assume this whenever we use $\sym(\H)$.
A simple manner to force the necessary distinctions for the $\pi_e$ 
is to attach to $\H$, as a disjoint component, a copy of the hypercube
$2^E$ (and this modification will nowhere interfere with other
concerns of our constructions). 

The Biggs group $\G(E,n)$ as discussed above is $\sym(\T(E,n))$.

Recall that we let permutations act from the right. In terms of
the group product in $\sym(\H)$ this makes $\pi_e \pi_{e'} = \pi_{e'}
\circ \pi_e$. Extending this to arbitrary words $w = e_1 \cdots e_n \in
E^\ast$ over $E$ according to 
\[
\barr{rcl}
[\,\nt\,]_{\H} \colon E^\ast &\longrightarrow& \sym(\H)
\\
w &\longmapsto& [w]_{\H} := \pi_w := \prod_{i=1}^n \pi_{e_i} = \pi_{e_n} \circ \cdots \circ \pi_{e_1},
\earr
\]
yields a surjective homomorphism from the free monoid structure of $E^\ast$,
with concatenation and neutral element $\lambda \in E^\ast$, 
onto the group structure of $\sym(\H)$ with composition and
neutral element $\pi_\lambda = \mathrm{id}_V$. Factorisation w.r.t.\ 
the identities $e = e^{-1}$ turns this into a surjective group homomorphism
from the free $\E$-group onto $\sym(\H)$.

\bD[Cayley graph]
\label{Cayleydef}
\mbox{}\\
For an abstract group $\G = (G,\cdot\,,1)$ and any set $E \subset G$ of
generators, the \emph{Cayley graph} of $G$ w.r.t.\ $E$ is the directed 
edge-coloured graph $\CG := \mathrm{Cayley}(\G,E)= (G,(R_e)_{e \in E})$ 
with vertex set $G$ and edge sets 
\[
R_e := \{ (g,ge) \colon g \in G \}
\]
of colour $e$, for all $e \in E$.
\eD

The Cayley graph $\CG$ is undirected precisely if the generator set 
$E$ consists of involutions of $\G$. In general the $R_e$ will not be
symmetric, but each $R_e$ will always be the graph of a global
permutation $\pi_e$ of the vertex set $G$, viz.\ of right multiplication with
$e \in G$, $\pi_e \colon g \mapsto ge$. 
It is easy to check that, as an abstract group with generators $e \in E$, 
$\G$ is isomorphic to the subgroup of the full symmetric group 
$\mathrm{Sym}(G)$ over the vertex set $G$ generated by these
permutations $\pi_e$. In particular, in the case of a group 
$\G = (G,\cdot\,,1)$ 
that admits a set of involutive generators $E \subset G \setminus \{1\}$, the
associated Cayley graph $\CG = \mathrm{Cayley}(\G,E)$ 
is a complete and strict $\E$-graph in the sense of 
Definition~\ref{Egraphdef}, and 
\[
\G = (G, \cdot\,, 1) \;\simeq\; \sym(\CG).
\]

In the following it will be convenient, and without risk of confusion, to
identify the generators $e \in E \subset G$ of a group $\G$ with 
the maps $\pi_e \colon g \mapsto g e$ in $\G$ or in its Cayley
graph $\CG$. We similarly identify the family of generators 
$(\pi_e)_{e \in E}$ of $\sym(\H)$ with a subset $E \subset
\sym(\H)$ whenever $\sym(\H)$ is an $\E$-group,
by writing just $e$ instead of $\pi_e$ in this context.

\bD[generated subgroup]
\label{gensubgroupdef}
\mbox{}\\
For a subset $\alpha \subset E$ of the set of involutive generators
$E$ of an $\E$-group $\G$ we let 
$\G[\alpha]$ stand for the subgroup generated by $\alpha$, regarded as
an $\alpha$-group whose
universe is 
\[
G[\alpha] := \{ [w]_{\G} \colon w \in \alpha^\ast \} \subset G.
\] 

The Cayley graph $\CG[\alpha]$ of $\G[\alpha]$, correspondingly, is
regarded as an $\alpha$-graph, which is a weak subgraph 
$\CG[\alpha] \subsetw \CG$ of the Cayley graph of $\G$.%
\footnote{More specifically it is the $(R_e)_{e \in \alpha}$-reduct of the 
induced subgraph $\CG\restr G[\alpha]$ on $G[\alpha] \subset G$.}
\eD

\bD[$\alpha$-walk and $\alpha$-component]
\label{alphawalksconnectdef}
\mbox{}\\
For a subset $\alpha \subset E$ and an $\E$-graph $\H = (V,(R_e)_{e \in E})$,
an \emph{$\alpha$-walk} of length $n$ from $v$ to $v'$ 
is a sequence $v_0, e_1, v_1, e_2, \ldots, e_n, v_n$ of vertices and edge labels
where $v_i \in V$, $v= v_0$, $v_n = v'$, $e_i \in \alpha$ such that $(v_i,v_{i+1})
\in R_{e_{i+1}}$ for $i < n$. 
\\
The \emph{$\alpha$-connected component}, or just $\alpha$-component,
of $v \in
V$ consists of those vertices $v'$ that are linked to $v$ by
$\alpha$-walks. 
\eD

We typically write $\alpha(v) \subset V$ for this set of
vertices and $\H[\alpha;v]$ for the weak subgraph $\H[\alpha;v]
\subsetw \H$ obtained, as an $\alpha$-graph, as a reduct of 
the induced subgraph $\H \restr \alpha(v)$.
 
Note that the Cayley graph $\CG[\alpha] \subsetw \CG$ of 
$\G[\alpha]$ also arises as the $\alpha$-component 
of $1 \in G$ in the Cayley graph $\CG$. 
It is also useful to note that, if 
$v= v_0, e_1, v_1, e_2, \ldots, e_n, v_n=v'$ is an $\alpha$-walk from
$v$ to $v'$ in the $\E$-graph $\H$ such that $w = e_1\cdots e_n$ traces the edge
labels along this walk, then $v' = \pi_w(v) = [w]_{\H}(v)$ w.r.t.\ the 
permutation group action of $\sym(\H)$ on $\H$.

\subsection{Compatibility and homomorphisms}

The notion of a homomorphism
between $\E$-groups is the natural one. 
It requires compatibility with the group product \emph{and} with the
identification of the generators. 
We write $\Ghat \succeq \G$  or $\G \preceq \Ghat$ to 
indicate that there is a homomorphism from $\Ghat$ to $\G$.
If there is any homomorphism 
$h \colon \Ghat \rightarrow \G$ between 
$\E$-groups $\Ghat$ and $\G$ then it must be 
\[
\barr{rcl}
h \colon \;\; \Ghat\;\;  &\longrightarrow &\G
\\
\nt
[w]_{\Ghat}\!\! &\longmapsto & [w]_\G
\earr
\]
for all $w \in E^\ast$. So what matters is
well-definedness of this mapping, which is expressible as 
the condition that $[w]_\G = [u]_\G$ whenever 
$[w]_{\Ghat} = [u]_{\Ghat}$, or just that $[w]_\G =1$ in $\G$ whenever 
$[w]_{\Ghat}  =1$ in $\Ghat$. We may think of $\Ghat \succeq \G$ as an 
\emph{unfolding} of $\G$, which also captures the relationship between the
associated Cayley graphs $\CG$  and $\CG$.

\bD[compatibility]
\label{compatdef}
\mbox{}\\
For an $\E$-graph $\H$ and $\E$-group $\G$ we say that 
$\G$ is \emph{compatible with} $\H$ if there is a 
homomorphism of $\E$-groups from $\G$ to $\sym(\H)$, i.e.\
if $\sym(\H) \preceq \G$.  More generally, for $\alpha \subset \E$, 
$\G[\alpha]$ is \emph{compatible with} the $\E$-graph $\H$ if 
$\G[\alpha]$ is compatible with the $\alpha$-reduct $\H\restr \alpha$
of $\H$. 
\eD

In straightforward extension of this concept, a family of $\E$-groups 
is compatible with $\H$ if each member is; this will be of interest 
especially when certain families of (small) generated subgroups of $\G$, 
rather than $\G$ itself, are compatible with some $\E$-graph.

Recall that $[w]_{\H} = \pi_w = \prod_{i=1}^n \pi_{e_i} \in \sym(\H)$
for $w = e_1\cdots e_n \in E^\ast$.  
Compatibility of $\G$ with $\H$ precisely requires that the mapping
\[
\textstyle
[w]_\G 
\;\longmapsto\; 
[w]_\H \in \sym(\H) 
\]
is well-defined, i.e.~that 
$[w]_\G = 1$ in $\G$ implies $[w]_\H = 1$ in $\sym(\H)$.

Note that trivially $\sym(\H)$ is compatible with $\H$ and with
every connected component of $\H$. 
We collect some further simple but useful facts.

\bO
\label{easycompatobs}
\bre
\item
$\G[\alpha]$ is compatible with $\CG[\alpha]$ for $\alpha \subset E$. 
\item
$\G$ is compatible with the disjoint union $\bigoplus_i \H_i$ of
$\E$-graphs $\H_i$ if, and only if, it is compatible with each
component $\H_i$.
\ere
\eO

It also follows that $\G \preceq \Ghat$ if, and only if,  
$\Ghat$ is compatible with $\CG$. A version of this 
observation for generated subgroups will be crucial in the
construction of suitable $\E$-groups with specific acyclicity 
properties.

\bL
\label{GhatCGcompatlem}
Let $\Ghat \succeq \G$ be  $\E$-groups, 
$\Ghat = \sym(\H)$ for an $\E$-graph $\H$.
In this situation, the subgroups $\G[\alpha]$ and $\Ghat[\alpha]$ 
generated by $\alpha \subset E$ are isomorphic as $\alpha$-groups, 
$\Ghat[\alpha] \simeq \G[\alpha]$, if the homomorphism 
from $\Ghat$ to $\G$ is injective in restriction to $\Ghat[\alpha]$,
which is the case if, and only if, already $\G[\alpha]$ is compatible with 
every $\alpha$-component of $\H$ and hence with $\H$.
\eL

\prf
For the last claim, assuming that $\G \preceq \Ghat$, we need to show
that conversely $\Ghat[\alpha]\preceq \G[\alpha]$. Note that
$\Ghat[\alpha] = \sym(\H\restr \alpha)$, where 
$\H\restr \alpha$ stands for the $\alpha$-reduct of $\H$,  
which is an $\alpha$-graph. If $\G[\alpha]$ is compatible with every
connected component of $\H\restr \alpha$, then $[w]_\G = [u]_\G$
for $w,u \in \alpha^\ast$ implies that $\pi_w = \pi_u$ in
$\sym(\H\restr \alpha)$ and therefore also in
$\sym(\H)$, i.e.\ in $\Ghat$.
\eprf

Our constructions of finite $\E$-groups with special combinatorial
properties proceed by induction on the number of generators
for subgroups under consideration, i.e.\ by induction on the size $|\alpha|$ of
generator subsets $\alpha \subset E$. Subsets $\alpha$ of the same size are
always treated in parallel in order to guarantee that our constructions
are fully isomorphism-preserving and do not break any symmetries
(unlike a construction governed by, for instance, a chosen enumeration of the
generator set~$E$).
The desired properties will successively be achieved for subgroups generated 
by increasing numbers of generators. The induction is based on
suitable unfolding steps for the passage from $\G$ to $\Ghat \succeq \G$, where  
the desired property for $\Ghat[\alpha]$ relies on and preserves the behaviour of 
the subgroups $\G[\alpha']$ of $\G$ for $\alpha' \strictsubset
\alpha$; in particular, all generated subgroups of $\Ghat$ generated
by up to $k$ generators inherit the property at hand from the
subgroups of $\G$ that are generated by fewer than $k$ generators. So, for
$\alpha \subset E$, we let 
\[
\Gamma_\alpha := \{ \alpha' \subset E \colon
\alpha' \strictsubset \alpha \},
\]
and,  for $1 \leq k \leq |E|+1$,  
\[
\Gamma_k := \{ \alpha \subset E \colon
|\alpha| < k \},
\]

We denote corresponding families of generated subgroups and their
Cayley graphs in a given $\E$-group $\G$ as
\[
\barr{r@{\;:=\;}l}
\Gamma_k(\G) &
(\G[\alpha'] \colon \alpha' \in \Gamma_k)
\\
\hnt
\Gamma_k(\CG) &
(\CG[\alpha'] \colon \alpha' \in \Gamma_k)
\earr
\]

\bD
\label{preceqalphakdef}
For $\alpha \subset E$ and $\E$-groups $\G$ and $\Ghat$,  
$\G \preceq_\alpha \Ghat$ denotes the relationship that
$\G \preceq \Ghat$ 
and
\[
\Gamma_\alpha(\G) = \Gamma_\alpha(\Ghat),
\] 
i.e.\ that subgroups generated by proper subsets of $\alpha$ 
are preserved in the unfolding. 
For $1 \leq k \leq |E|+1$,  
$\G \preceq_k \Ghat$ analogously stands for an unfolding relationship 
$\G \preceq \Ghat$ in which 
\[
\Gamma_k(\G) = \Gamma_k(\Ghat),
\]
where $\G$ and $\Ghat$ agree on all 
subgroups generated by fewer than $k$ generators. 
\eD

If we think of the passage from $\G$ to $\Ghat$
in $\G \preceq \Ghat$ as an \emph{unfolding}, then 
$\G \preceq_k \Ghat$ says that the unfolding is trivial up to the 
level of $\Gamma_k$, or is \emph{conservative} w.r.t.\ $\Gamma_k$-generated subgroups. 

Recall that $\Ghat = \sym(\H \oplus \CG)$ guarantees
$\Ghat \succeq \G$. In this situation, morevover 
$\Ghat \succeq_k \G$ if already $\G$ itself is compatible
with all $\alpha$-components of $\H \oplus \CG$ for $\alpha \in \Gamma_k$.
It is clear that $\G[\alpha]$ is compatible with $\CG[\alpha]$ for all
$\alpha$; as we shall in Section~\ref{cosetcycsec},
compatibility of $\G$ (rather than $\G[\alpha]$) with all $\CG[\alpha] \in \Gamma_k(\CG)$
implies that $\CG$ does not admit $2$- or $3$-cycles formed by 
cosets generated by fewer than $k$ generators. In particular, so that 
compatibility of $\G$ with all $\CG[\alpha]$ implies $3$-acyclicity of
$\G$, as stated in Lemma~\ref{3acyclem}. 
The following summarises the content of Observation~\ref{easycompatobs} 
and Lemma~\ref{GhatCGcompatlem}.

\bL
\label{Ghatcompatlem}
Consider a functor $\F$ that maps $\E$-groups $\G$ to 
$\E$-graphs $\F(\G)$ in an isomorphism respecting manner
(i.e.\ $\F(\G) \simeq \F(\G')$ whenever $\G \simeq \G'$). 
Then for $\Ghat := \sym(\F(\G))$:
\bae
\item
If $\F(\G)$ has a component isomorphic to $\CG$, then $\G \preceq \Ghat$.
\item
If $\F$ is as in~(a) and such that all subgroups in 
$\Gamma_\alpha(\G)$
are compatible with $\F(\G)$, then $\Ghat[\alpha'] \simeq
\G[\alpha']$ for all $\alpha'\strictsubset \alpha$, 
i.e.\ $\G \preceq_\alpha \Ghat$.
\item
If $\F$ is as in~(a) and such that all generated subgroups 
$\G[\alpha] \in \Gamma_k(\G)$ are compatible with 
$\F(\G)$, then $\Ghat[\alpha] \simeq \G[\alpha]$ for  
$\alpha \in \Gamma_k$, i.e.\ $\G \preceq_k \Ghat$.
\eae
\eL

The pre-conditions in~(b) and~(c) should be seen as 
\emph{downward compatibility} conditions. They 
guarantee that the transition from $\G$ to $\Ghat$ 
is conservative w.r.t.\ to smaller generated subgroups, while possibly
unfolding $\G$ in non-trivial ways at the level of larger generated subgroups.

This opens up the potential for achieving successively more 
stringent structural conditions in an inductive fashion. Essentially,
the induction will be on the size $|\alpha|$ of the generator set of
subgroups $\G[\alpha]$ that may form certain 
obstructive patterns and progresses to exclude them by replacing $\G$
essentially by suitable $\Ghat := \sym(\CG \oplus \H)$. 
\begin{quotation}\noindent
The crux of the matter is to find $\H$ that eliminates the obstructive 
patterns at level $\Gamma_{k+1}$ while retaining the subgroups at
level $\Gamma_k$, which have already been righted in previous steps.
\end{quotation}

Slightly generalising the situation of the last lemma, 
consider some functor $\F$ that maps families of generated subgroups
$\G[\alpha]$ of $\E$-groups $\G$, for $\alpha \subset E$, to 
$\E$-graphs in an isomorphism respecting manner. 
Let $\F(\Gamma_k(\G))$ stand for the $\F$-images of families of
$\Gamma_k$-generated subgroups of $\G$, $\F(\G)$ for their union 
$\F(\G) = \bigcup_k \F(\Gamma_k(\G))$ and assume that $\F(\G)$ is
finite up to isomorphism, for every $\G$.

\bD
\label{conservedef}
We say that $\F$ is \emph{conservative} if, for
$\beta \in \Gamma_{k+1}$, all $\beta$-components of
$\E$-graphs in $\F(\Gamma_{k+1}(\G))$ 
are (isomorphic to) $\CG[\beta]$ or in $\F(\Gamma_k(\G))$. 
\eD

Examples of this include amalgamation
chains over $\alpha_i$-graphs $\CG[\alpha_i] \subsetw \CG$ (of some bounded
length) to be treated in Section~\ref{amalgsec.chains}, and free amalgamation 
clusters of $\alpha_i$-graphs 
$\CG[\alpha_i] \subsetw\CG$ to be treated in Section~\ref{amalgsec.clusters}.

\bP
\label{conserveprop}
Let $\F$ be a conservative functor in the above sense,
$\G$ an $\E$-group, $1 \leq k \leq |E|+1$. 
If $\G$ is compatible with the $\E$-graphs in $\F(\Gamma_k(\G))$,
then $\G$ admits an unfolding $\Ghat \succeq_k \G$
that is compatible with all $\E$-graphs in $\F(\G)$. 
\eP

\prf
Inductively, one
obtains a sequence of unfoldings
\[
\G = \G_k \preceq_k \G_{k+1} \preceq_{k+1} \cdots \preceq_{n-1}
\G_n \preceq_n \G_{n+1} =: \Ghat
\]
for $n = |E|$. Each step  
$\G_m \preceq_m \G_{m+1}$ is such that 
$\G_{m+1}$ is made compatible with $\F(\Gamma_{m+1}(\G_m))$:
$\G_k$ is compatible with $\F(\Gamma_k(\G))$ by assumption, and 
$\G_{m+1}$ is obtained as $\sym(\H_m)$ where $\H_m$ is the disjoint union
over the $\E$-graphs $\CG_m$ and $\F(\Gamma_{m+1}(\G_m))$.
On the one hand, the component $\CG_m$ in $\H_m$ guarantees $\G_{m+1}
\succeq \G_m$. On the other hand, 
the conservative nature of $\F$ makes sure that 
$\G_{m+1}[\alpha] \simeq \G_m[\alpha]$  is preserved for
$\alpha \in \Gamma_m$, since all non-trivial $\alpha$-components of 
$\E$-graphs in $\F(\Gamma_{m+1}(\G_m))$ are in 
$\F(\Gamma_{m}(\G_m))$ so that already $\G_m[\alpha] \subset \G_m$ is compatible
with those. 
\eprf

\paragraph*{Symmetries.} 
Some application contexts call for an analysis of symmetries of 
$\E$-groups $\G$ that are induced by permutations of the underlying 
set $E$ of generators. These are not covered by the notion of 
automorphisms of $\E$-groups since those, 
as special homomorphisms, need to fix the generators individually 
(model-theoretically they are treated as constants).
Similarly for $\E$-graphs $\H$, automorphisms of $\H$ viewed 
as a relational structure need to respect each $R_e$ 
individually, and do not account for symmetries induced by
permutations of the edge colours. In both cases, permutations $\rho
\in \mathrm{Sym}(E)$ induce what in model-theoretic terminology is a
\emph{renaming}, sending $\G$ to $\G^\rho$ and 
$\H$ to $\H^\rho$. For instance the $\rho$-renaming of the $\E$-graph  
$\H = (V, (R_e)_{e \in E})$ is $\H^\rho = (V, (R'_e)_{e \in
  E})$ with $R'_{\rho(e)} = R_e$. 
Such a renaming reflects a \emph{symmetry} if it
leaves the underlying structure invariant up to isomorphism.

\bD[symmetry over $E$]
\label{Esymmdef}
\mbox{}\\
A permutation $\rho \in \mathrm{Sym}(E)$ of the set 
$E$ is a \emph{symmetry} of an $\E$-group $\G$
if the renaming of generators according to $\rho$ 
yields an isomorphic $\E$-group, $\G^\rho \simeq \G$.  
Similarly, $\rho$ is a symmetry of the $\E$-graph $\H$
if the renaming of its edge relations according to $\rho$ 
yields an isomorphic $\E$-graph: $\H^\rho \simeq \H$.  
\eD

For instance, the trees $\T(E,n)$ in Biggs' construction 
are \emph{fully symmetric} in the sense that every 
$\rho \in \mathrm{Sym}(E)$ is a symmetry; the same 
is then true of the resulting $\E$-group $\G = \sym(\T(E,n))$ and its 
Cayley graph $\CG$.

\section{Coset cycles and acyclicity criteria}
\label{cosetcycsec}

A basic notion of $n$-acyclicity for $\E$-groups would just
forbid non-trivial \emph{generator cycles} (i.e.\ representations of $1 \in
\G$ by reduced generator words) of lengths up to $n$.
This account matches the graph-theoretic notion of \emph{girth} for the associated
Cayley graph, as the length of the shortest 
non-trivial generator cycle in $\G$ is the length of the shortest
graph cycle in $\CG$, i.e.\ its girth.
We are here interested in a more liberal notion of cycles, which leads to a
more restrictive notion of acyclicity that
forbids short \emph{coset cycles}, i.e.\ cyclic configurations of 
cosets $g_i \CG[\alpha_i]$.

\bD[coset cycle]
\label{cosetcycledef}
\mbox{}\\
Let $\G$ be an $\E$-group, $n \geq 2$.
A \emph{coset cycle} of length $n$ in $\G$ is a 
cyclically indexed sequence of pointed cosets 
$(g_i \G[\alpha_i],g_i)_{i \in \Z_n}$ w.r.t.\ subgroups 
$\G[\alpha_i]$ for $\alpha_i \subset E$ satisfying these conditions:
\bre
\item
(connectivity)
$g_{i+1} \in g_i \G[\alpha_i]$, i.e.\ $g_i \G[\alpha_i] =
g_{i+1} \G[\alpha_i]$; 
\item
(separation)
$g_i \G[\alpha_{i,i-1}] \cap g_{i+1} \G[\alpha_{i,i+1}] = \emptyset$,
\ere
where $\alpha_{i,j} := \alpha_i \cap \alpha_j$.
\eD

We sometimes put a focus on coset cycles whose constituent cosets 
stem from a restricted family of generated subgroups, and especially 
from $\Gamma_k(\G)$ for some $1 \leq k \leq |E|$. With terminology 
like \emph{coset cycle w.r.t.\ $\Gamma_k$} we then refer to coset
cycles $(g_i \G[\alpha_i])_{i \in \Z_n}$ with $\alpha_i \in \Gamma_k$,
i.e.\ with $|\alpha_i| < k$.

\bD[$N$-acyclicity]
\label{Nacycdef}
\mbox{}\\
For $N \ge 2$, an $\E$-group $\G$ or its Cayley graph $\CG$ are 
\emph{$N$-acyclic} if they admit no coset cycles of lengths up to $N$. 
\eD

Correspondingly, $N$-acyclicity w.r.t.\ $\Gamma_k$ 
forbids coset cycles $(g_i \G[\alpha_i])_{i \in \Z_n}$ 
of lengths $n \leq N$ with $\alpha_i \in \Gamma_k$.
Obviously, non-trivial generator cycles are very special coset cycles 
with singleton sets $\alpha_i = \{ e_i\}$.
So $N$-acyclicity w.r.t.\ $\Gamma_2(\G)$ 
precisely says that the girth of $\G$ or $\CG$ is larger than $N$. 
Also note that $N$-acyclicity w.r.t.\ $\Gamma_k(\G)$ in particular implies
outright $N$-acyclicity for $\G[\alpha]$ for all $|\alpha| \leq k$.
For $G$ itself outright  $N$-acyclicity is the same as 
$N$-acyclicity w.r.t.\ $\Gamma_{|E|}(\G)$. 

It is important to note that the graph-theoretic diameter of an 
$\alpha_i$-coset in the Cayley
graph or the cardinality of $\G[\alpha_i]$
cannot be uniformly bounded (e.g.\ in terms of $|\alpha_i|$). 
Therefore no level of generator acyclicity captures any fixed level of coset acyclicity.

The lowest level of coset acyclicity, viz.\ $N$-acyclicity for $N=2$,
is of special interest. It is easy to check that the condition for
$2$-acyclicity is equivalent to an intersection condition on pairs of
cosets, which is reminiscent of a notion
of simple connectivity.
This view also points to a natural 
generalisation of the concept from Cayley graphs to more general
$\E$-graphs, later in Definition~\ref{twoacycgendef}.

\bO
\label{2acycobs}
An $\E$-group$\G$ is $2$-acyclic if, and only if, for all
$\alpha_1,\alpha_2 \strictsubset E$,
\[
\G[\alpha_1] \cap \G[\alpha_2] = \G[\alpha_1 \cap \alpha_2].
\]
\eO

The following observation recasts $3$-acyclicity as a combinatorial
feature for clique-like configurations of cosets.

\bO
\label{3acycobs}
A $2$-acyclic $\E$-group $\G$ is $3$-acyclic if, and only if,
every finite collection of cosets
$(g_i\G[\alpha_i])_{i \in  I}$ with pairwise non-empty intersections
has non-empty intersection overall:
$\;g_i \G[\alpha_i] \cap g_j\G[\alpha_j] \not= \emptyset
  \mbox{ for all }  i,j \in I$ implies  
  $\;\bigcap_{i \in I} g_i \G[\alpha_i] \not= \emptyset$.
\eO 

\prf
For $3$ cosets in a $3$-acyclic $\E$-group the intersection claim follows
directly from the $3$-acyclicity criterion: the violation of the
separation condition between two cosets yields an element in the
intersection of all three. The apparently stronger consequence
for larger collections of pairwise intersecting
cosets follows by induction. For the induction step we may replace
two intersecting cosets that intersect all the others by their
intersection, which (by the base case) still intersects each one of the remaining
cosets. 
\eprf

As we shall see in Section~\ref{amalgsec.clusters},
any collection of $\alpha_i$-cosets with pairwise non-empty
intersections in a $3$-acyclic $\E$-group as above embeds into its
Cayley graph as a
free amalgamation cluster, with overall itersection isomorphic to   
$\CG\restr \bigcap_I g_i \G[\alpha_i] \simeq \CG\restr \bigcap_I
\G[\alpha_i] \simeq \CG[\bigcap_{i \in I} \alpha_i] \subseteq \CG$.

The following associates acyclicity criteria with closure 
properties and minimal supporting sets of generators -- another concept
that will be generalised from Cayley graphs of $\E$-groups to other
classes of $\E$-graphs in Sections~\ref{amalgsec} and~\ref{CEsec}.

\bR
\label{elementcosetsupportrem}
Consider an element $g \in \G$ and its $\beta$-component 
$B = g \G[\beta] \subset \G$ for some 
$\beta \subset E$ in an $\E$-group $\G$.
\bre
\item
If $\G$ is $2$-acyclic then there is a unique $\subset$-minimal generator set
$\alpha_0 \subset E$ such that $g \in \G[\alpha_0]$, viz.\
$\alpha_0 := \bigcap \{ \alpha \subset E \colon g \in \G[\alpha] \}$.
\item
If $\G$ is $3$-acyclic then there is a unique $\subset$-minimal generator set
$\alpha_0 \subset E$ such that $\G[\alpha_0] \cap B \not= \emptyset$, viz.\ 
$\alpha_0 := \bigcap \{ \alpha \subset E \colon \G[\alpha_0] \cap B \not= \emptyset \}$.
\ere
\eR

In connection with~(ii) it is not hard to see that, conversely, the
stated uniqueness property implies $3$-acyclicity. 

\prf
Observation~\ref{2acycobs} implies that the family of subsets $\alpha
\subset E$ for which $g \in \G[\alpha]$ is closed under intersections;
this immediately implies claim~(i). Towards claim~(ii) consider two
elements $g_i \in B$ and their minimal supporting sets of generators 
$\alpha_i$ according to~(i), for $i=1,2$. We need to show that also $\alpha_0 :=
\alpha_1 \cap \alpha_2$ supports $B$ in the sense that 
$\G[\alpha_0] \cap B \not= \emptyset$. For this consider the potential
$3$-cycle of cosets $1 \G[\alpha_1] = \G[\alpha_1]$, $g_1 \G[\beta]$ and
$g_2\G[\alpha_2] = \G[\alpha_2]$. As $\G$ does not admit $3$-cycles of
cosets, at least one of the three instances of the separation conditions 
in Definition~\ref{cosetcycledef} must fail. We argue that each such failure 
yields an element in $B \cap \G[\alpha_0]$. Failure in the link
$\G[\alpha_1]$ means that there is some
$g' \in g_1 \G[\alpha_1\cap \beta] \cap 1 \G[\alpha_1 \cap \alpha_2] $;
now $g' \in g_1 \G[\alpha_1\cap \beta]$ implies that $g' \in B$
since $g_1 \in B$, and $g' \in \G[\alpha_0]$ as 
$1 \G[\alpha_1 \cap \alpha_2] = \G[\alpha_0]$.
Failure in the link
$\G[\alpha_2]$ is entirely symmetric.
Failure in the link $g_1 \G[\beta]$ finally means that there is some
$g' \in g_1 \G[\alpha_1\cap \beta] \cap g_2 \G[\alpha_2 \cap \beta]$
which implies that $g' \in B$ as before. Moreover $g' \in \G[\alpha_i]$ 
as $g_i \in \G[\alpha_i]$ and $g' \in g_i \G[\alpha_i]$ for $i=1,2$; so 
$g' \in \G[\alpha_1] \cap \G[\alpha_2] = \G[\alpha_1 \cap \alpha_2] = \G[\alpha_0]$.
\eprf

\bL
\label{3acyclem}
If an $\E$-group $\G$ is compatible with $\CG[\alpha]$ for every
$\alpha \strictsubset E$ then $\G$ is $3$-acyclic.
\eL

\prf
Assume that $(g_i \G[\alpha_i],g_i)_{i \in \Z_3}$ 
formed a coset $3$-cycle in $\G$. 
The separation condition implies
that $\alpha_i \strictsubset E$ and that, for instance,
\[
(\ast) \qquad g_0 \G[\alpha_{0,2}]\cap g_1 \G[\alpha_{0,1}] = \emptyset.
\]

By the connectivity condition, there are $w_j \in \alpha_j^\ast$ such that 
\[
g_0^{-1} g_1 = [w_0]_\G, \;\; g_1^{-1} g_2 = [w_1]_\G, \;\;  g_2^{-1} g_0
= [w_2]_\G.
\]
Clearly $[w_2 w_0 w_1 ]_\G = 1$. 
Let  $w_{0,j} \in \alpha_{0,j}^\ast$ be the projection of the
word $w_j$ to $\alpha_0^\ast$, as obtained by deletion of all letters
$e \not\in \alpha_0$. If $\G$ is compatible with $\G[\alpha_0]$ 
then  $[w_2 w_0 w_1 ]_\G = 1$ implies that $[w_2w_0w_1]_\H = 1 \in
\sym(\H)$ for $\H = \CG[\alpha_0]$, which in turn implies that
$[w_{0,2} w_0 w_{0,1} ]_\G = 1$ since $[e]_\H$ is trivial for $e \not\in \alpha_0$. 
But the operation of the corresponding sequence of 
generators of $\sym(\CG[\alpha_0])$ maps the element 
$g := [w_{0,2}]^{-1}_\G$ (via $1$ and $g_0^{-1} g_1 = [w_0]_\G$) to  $g' := g_0^{-1} g_1 [w_{0,1}]_\G$. 
As $g \in \G[\alpha_{0,2}]$ and $g' \in g_0^{-1} g_1\G[\alpha_{0,1}]$, 
which by $(\ast)$ are disjoint subsets of $\G[\alpha_0]$, 
it follows that $[w_{0,2} w_0 w_{0,1} ]_\G \not= 1$, a contradiction. 
Analogously, $\G$ cannot admit $2$-cycles.
\eprf

\bL
\label{acychomlem} 
$N$-acyclicity of $\G[\alpha]$ is preserved under inverse 
homomorphisms that are injective on $\alpha'$-generated subgroups 
for all $\alpha'\strictsubset \alpha$.
\eL

\prf
If $h \colon \Ghat[\alpha] \rightarrow \G[\alpha]$ is a 
homomorphism of $\alpha$-groups that is a local isomorphism in restriction to each 
$\Ghat[\alpha']$ for $\alpha'\strictsubset \alpha$, then $h$ maps 
a coset cycle in $\Ghat[\alpha]$ to a coset cycle in $\G[\alpha]$. The connectivity
condition is obviously maintained under~$h$. The crux 
of the matter is the separation condition for links in a potential 
coset cycle. As each $\alpha_i$-coset in $\Ghat$, for $\alpha_i \strictsubset
\alpha$, is mapped bijectively onto its image coset in $\G$, so are
the disjoint critical $\alpha_{i,i\pm 1}$-cosets as its subsets.
\eprf

The method suggested by
Proposition~\ref{conserveprop}
is to be used to eliminate $N$-cycles of cosets with increasing numbers of
generators, in an inductive treatment. 
The relevant configurations in $\F(\G)$ will be suitable amalgams of $\E$-graphs
$\CG[\alpha]$ for $\alpha \in \Gamma_k$. 
In light of the proposition we need to focus on functors $\F$
that are conservative in the sence of Definition~\ref{conservedef}.
Towards coset acyclicity we can use as $\F$-images amalgamation 
chains that unfold potential cycles.

\section{Free amalgams of E-graphs}
\label{amalgsec}

We investigate free amalgams of $\E$-graphs, and especially of
Cayley graphs of generated subgroups $\G[\alpha]$. The goal is to provide elements of a 
structure theory for Cayley graphs that provides tools for partial
unfoldings that eliminate obstructions to desirable combinatorial
properties.

The core idea of free amalgams of  $\E$-graphs is to superpose copies 
of suitably matching $\E$-graphs in distinguished elements or regions 
with just minimal identifications. Those minimal identifications, or 
regions of overlap, are determined by shared generator edges.  
We first look at 
superpositions of two $\E$-graphs, then expand the idea to certain
chains or clusters of several $\E$-graphs and in particular of
$\E$-graphs of the form $\CG[\alpha_i]$ that arise as substructures of
a common $\E$-group $\G$. In each case, 
the isomorphism type of the resulting $\E$-graphs is fully determined
by the isomorphism types of the constituents --  and not by the
manner in which e.g.\ constituents $\CG[\alpha_i]$ are embedded in 
$\CG$ as weak subgraphs $\CG[\alpha_i] \subsetw \CG$. 
If $\G$ satisfies appropriate acyclicity conditions, however, the
free amalgam will be naturally isomorphic to the corresponding
embedded weak subgraph of $\CG$ on a union of $\alpha_i$-cosets of $\G$
(cf.~Observations~\ref{twoamalgobs},~\ref{amalgchainobs},~\ref{amalgclusterobs}).
For instance, if $\G$ is $2$-acyclic, then the
weak subgraphs of $\CG$ on overlapping cosets of the form
$g G[\alpha_1] \cup g G[\alpha_2] \subset G$
(with induced $\alpha_i$-edges on the $\alpha_i$-coset) will all be
isomorphic to the free amalgam $\CG[\alpha_1] \oplus \CG[\alpha_2]$
(Observation~\ref{twoamalgobs}).

\bD[free amalgam]
\label{twoamalggendef}
\mbox{}\\
A \emph{free amalgam} between two pointed $\E$-graphs $\H_1,v_1$ and
$\H_2,v_2$ is defined relative to a designated pointed 
$\E$-graph $\H_0,v_0$, which serves as the  
\emph{overlap}, if 
$\H_0$ is connected and admits isomorphic embeddings $\rho_i \colon 
\H_0,v_0 \rightarrow \H_i,v_i$ onto weak
substructures $\rho_i(\H_0) \subsetw \H_i$ 
that map $v_0$ to $v_i$, for $i=1,2$,
and if 
the following condition is satisfied for all $v \in \H_0$ and $e \in E$:
\[
(\ast) \quad
\Bigl( \rho_1(v)\in \dom(R_e^{\H_1})\mbox{ and } 
\rho_2(v) \in \dom(R_e^{\H_2}) \Bigr)
\;
\Rightarrow 
\;
v \in \dom(R_e^{\H_0}).
\]
In this situation we say that \emph{$\H_1,v_1$ and $\H_2,v_2$ admit a free
amalgam over $\H_0,v_0$}, and we obtain this free amalgam as
\[
\H_1,v_1 \oplus_{\H_0,v_0} \H_2,v_2
\]
from the disjoint union of $\H_1$ and
$\H_2$ by identification of vertices $\rho_1(v)$ and $\rho_2(v)$ for all
$v \in \H_0$. When referring to this free amalgam as a pointed
$\E$-graph we regard the copy of $v_0$, i.e.\ the identification of $v_1$
and $v_2$, as the new distinguished vertex. 
\eD

\begin{figure}
\[
\xymatrix{
& \H_1,v_1 \ar@{->}[dr]^{\sigma_1} &
\\
\H_0,v_0 \ar@{->}[ur]^{\rho_1} \ar@{->}[dr]_{\rho_2}& \qquad 
{\mbox{\large $\circlearrowright$}}
\qquad
& 
*++{\makebox(20,0)[l]{$\;\;\H,v = \H_1,v_1 \oplus_{\H_0,v_0} \H_2,v_2,v$}}
& \qquad & \qquad \nt
\\
& \H_2,v_2 \ar@{->}[ur]_{\sigma_2} & \nt
}
\]
\caption{Free amalgamation.}
\label{amalgfig}
\end{figure}

We also think of the constituents $\H_i,v_i$ for $i=1,2$ as well as of 
the overlap $\H_0,v_0$ as isomorphically embedded as weak substructures
of the free amalgam, via the maps $\sigma_i \colon \H_i \rightarrow
\H$ and $\sigma_i \circ \rho_i \colon \H_0 \rightarrow \H$, for both
$i=1$ and $i=2$,
as in the commuting diagram of Figure~\ref{amalgfig}.

Condition~$(\ast)$ serves to guarantee that the resulting 
structure is an $\E$-graph, i.e.\ that identifications of vertices
from $\H_1$ and $\H_2$ cannot introduce forking $e$-edges.
In all situations to be considered below, however, we shall be dealing
with more restricted settings where the embeddings $\rho_i$ are
onto induced substructures $\rho_i(\H_0) \subset \H_i$.
The most basic case occurs if the $\H_i$ arise as Cayley graphs of subgroups of a
common $\E$-group $\G$: $\H_i \simeq \CG[\alpha_i] \subsetw \CG$.
In this case the isomorphism type of the free amalgam of $\H_1$ and $\H_2$
does not depend on the choice of the 
distinguished vertices, and 
$\H_0$ is determined relative to $\H_1$ and $\H_2$. 
Together with the requirement of connectivity for
$\H_0$, condition~$(\ast)$ implies that the
overlap is $\H_0 \simeq \G[\alpha_0]$ for $\alpha_0 := \alpha_1 \cap \alpha_2$.
We correspondingly denote such free amalgams as just 
$\CG[\alpha_1] \oplus \CG[\alpha_2]$ whenever it is clear that 
we refer to an amalgam (rather than to a disjoint union).   
If  $\G$ is $2$-acyclic, then 
\[
\CG[\alpha_1] \oplus \CG[\alpha_2] 
\simeq \CG\restr (\alpha_1(1) \cup \alpha_2(1) )\subsetw
\CG, 
\]
and this (essentially unique) embedding even yields an induced substructure
relationship $\CG[\alpha_1] \oplus \CG[\alpha_2] \subset
\CG$ if $\G$ is $3$-acyclic.
For these  claims observe that no two elements of $G[\alpha_i]$
can be linked by any generator $e \not\in \alpha_i$ due to
$2$-acyclicity; similarly, $3$-acyclicity rules out
$e$-links for $e \not\in \alpha_1 \cup \alpha_2$ between
$G[\alpha_1] \setminus G[\alpha_2]$ and 
$G[\alpha_2] \setminus G[\alpha_1]$.

\bO
\label{twoamalgobs}
For $2$-acyclic $\G$ and $\alpha_1,\alpha_2 \subset E$, the free
amalgam $\CG[\alpha_1] \oplus \CG[\alpha_2]$ embeds into $\CG$ as a weak
substructure. For $3$-acyclic $\G$ it embeds as an induced
substructure. 
Moreover, $\G$ is $2$-acyclic if, and only if, 
all free amalgams of the form $\CG[\alpha_1] \oplus \CG[\alpha_2]$
embed into $\CG$ as weak substructures. 
\eO

Generalising $2$-acyclicity from Cayley graphs of 
$\E$-groups to arbitrary $\E$-graphs,
we obtain the following notion of simple connectivity in terms of 
$\alpha$-components.

\bD[$2$-acyclicity]
\label{twoacycgendef}
\mbox{}\\
An $\E$-graph $\H$ is \emph{$2$-acyclic} 
if, for all $\alpha_1,\alpha_2 \subset E$ with intersection 
$\alpha_0 := \alpha_1\cap \alpha_2$ and for all vertices $v$ of $\H$:
\[
\alpha_1(v) \cap \alpha_2(v) = \alpha_0(v).
\]
\eD

Note that $2$-acyclicity of $\H$ implies that any $\alpha$-component
$\H[\alpha;v]$ of $\H$, as an $\alpha$-graph, is not only a weak but an induced 
substructure $\H[\alpha;v]\subset \H$: distinct vertices of 
$\H[\alpha;v]$ are by definition
linked by an $\alpha$-walk; an additional $e$-edge between them
for $e \not\in \alpha$ would therefore constitute a $2$-cycle.

\bO
\label{2acyccompobs}
Let $\H_i,v_i$ for $i=0,1,2$ be such that the free amalgam
$\H,v := \H_1,v_1 \oplus_{\H_0,v_0} \H_2,v_2$ is well-defined, with 
embeddings $\rho_i \colon \H_0,v_0 \simeq \rho_i(\H_0),v_i
\subsetw \H_i,v_i$ for $i=1,2$.
If $\rho_i(\H_0) \subsetw \H_i$ is such that for $i=1,2$ every $\beta$-component 
of $\H_i$ intersects the overlap region $\rho(\H_0) \subsetw \H_i$ in 
at most a single $\beta$-component of $\rho_i(\H_0) \simeq \H_0$, then 
the $\beta$-components in $\H$ are (isomorphic to) 
$\beta$-components of just one of the $\H_i$ or isomorphic to
free amalgams of two $\beta$-components, one from each $\H_i$
over a $\beta$-component of $\H_0$; 
in particular 
\[
\H[\beta;v] \simeq \H_1[\beta;v_1] \oplus_{\H_0[\beta;v_0],v_0}
\H_2[\beta,v_2].
\]
\eO

The conditions in the observation are satisfied in particular 
for free amalgams of the form
$\CG[\alpha_1] \oplus \CG[\alpha_2]$, for $\CG[\alpha_i] \subsetw
\CG$ over $\CG[\alpha_0]$ for $\alpha_0 = \alpha_1 \cap
\alpha_2$, if  both $\G[\alpha_i]$ are $2$-acyclic. 
It is also implied by $2$-acyclicity of the constituent $\E$-graphs 
$\H_i$ for $i=1,2$ whenever these two $\E$-graphs allow for a free
amalgam over $\H_0,v_0 \simeq \H_i[\alpha_0;v_i],v_i$.

\prf
Let $B$ be the vertex set of the $\beta$-component in question.
If $B$ is fully contained 
in (the isomorphic copy of) one of the constituents $\H_i$ in $\H$, then it is 
isomorphic to a $\beta$-component of that $\H_i$.
Otherwise $B$ must contain a vertex in the
overlap of the constituent $\H_i$, and we consider w.l.o.g.\ the 
case of $v \in B$.  For $i=0,1,2$ let 
$B_i$ be the vertex set of the $\beta$-component
$\H_i[\beta;v_i]$ of $v_i$ in $\H_i$, viewed as a subset
of the vertex set $V$ of the amalgam $\H$ in terms of the natural embeddings.
Clearly $B_0 \subset B_i \subset  B$ for $i=1,2$. Under the
assumptions of the observation, moreover, $B_0 = B_1 \cap B_2$ and
this further implies that $B = B_1 \cup B_2$. It remains to argue that,
as a $\beta$-graph, $\H[\beta;v]$ is isomorphic to the free amalgam 
of $\H_1[\beta;v_1]$ and $\H_2[\beta;v_2]$ over $\H_0[\beta;v_0]$.
We first observe that this amalgam is well-defined, by means of the 
restrictions of the embeddings $\rho_i$ to $\H_0[\beta;v_0]\subsetw
\H_0$ (in particular condition~$(\ast)$ from Definition~\ref{twoamalggendef} 
naturally restricts to $e \in \beta$). Now $\H_i\restr B_i$ as a
$\beta$-graph is isomorphic to $\H_i[\beta;v_i]$, and by assumption, 
the embedded copies of these, for $i = 1,2$, in $\H$ 
overlap precisely in the embedded copy of $\H_0[\beta;v_0]$.
\eprf

We turn to two specific forms of amalgams of Cayley graphs of
generated subgroups: \emph{amalgamation chains} 
and \emph{amalgamation clusters} in \S~\ref{amalgsec.chains}
and \S~\ref{amalgsec.clusters}, respectively.  
Both patterns are motivated by desirable acyclicity properties w.r.t.\
local overlaps between $\G[\alpha_i]$-cosets in $\G$. In other words, 
they unfold overlaps among a family of $\CG[\alpha_i]$ in the
`tree-like' pattern encountered in sufficiently acyclic $\G$.

\subsection{Amalgamation chains} 
\label{amalgsec.chains}

\bD[amalgamation chain]
\label{amalgchaindef}
\mbox{}\\
Let $\G$ be an $\E$-group, $N \geq 1$ and 
$\alpha_i \strictsubset E$ for $1 \leq i \leq N$
with intersections $\alpha_{i,i+1} := \alpha_i \cap \alpha_{i+1}$
for $1 \leq i < N$. 
Consider the sequence of pointed 
$\E$-graphs $(\CG[\alpha_i],g_i)_{1 \leq i \leq N}$ and assume
that the $g_i \in \G[\alpha_i] \subsetw \CG$ are such that the cosets
$1 \G[\alpha_{i-1,i}] \subset \G[\alpha_i]$ and $g_i
\G[\alpha_{i,i+1}]\subset \G[\alpha_i]$ 
are disjoint in $\G[\alpha_i]$.
In this situation, the \emph{free amalgamation chain}
$\bigoplus_{i=1}^N (\CG[\alpha_i],g_i)$ is the $\E$-graph $\H$ obtained as the result of simultaneous free amalgamation of disjoint copies of relational 
structures $\H_i \simeq \CG[\alpha_i]$ via the shared 
embedded weak substructures $\H_{i,i+1} \simeq 
g_i \CG[\alpha_{i,i+1}] \subsetw \CG[\alpha_i]$ and 
$\H_{i,i+1} \simeq 1 \CG[\alpha_{i,i+1}] \subsetw \CG[\alpha_{i+1}]$.
\eD

Note that the precondition on disjoint overlaps w.r.t.\ the
next neighbours in the chain ensures that there is no interference
between the pairwise amalgamation processes between next neighbours;
it also means that (the images of) the cosets $g_i \G[\alpha_{i,i+1}]$ for
$1 \leq i < N$ are graph-theoretic separators along the chain.

Amalagamtion chains can be thought of as unfoldings of coset cycles.

\bO
\label{amalgchainobs}
Whenever the amalgamation chain $\bigoplus (\CG[\alpha_i],g_i)$
is defined, there is a unique homomorphism from $\bigoplus (\CG[\alpha_i],g_i)$ to 
$\CG$ that maps $1$ in (the isomorphic copy of) $\CG[\alpha_1]$ to
$1 \in \CG$.  If the $\E$-group $\G$ is $N$-acyclic then this 
homomorphism is injective and $\bigoplus (\CG[\alpha_i],g_i)$ 
is realised as a weak substructure of $\CG$ (via an essentially 
canonical isomorphism).
\\
The following are equivalent:
\bre
\item
  $\G$ is $N$-acyclic
\item
$\bigoplus_i (\CG[\alpha_i],g_i) \subsetw \CG$ for any 
sequence of up to $N$ many pointed generated 
subgroups $(\G[\alpha_i],g_i)$ such that
the amalgamation chain is defined.
\ere
\eO

For (i)~$\Rightarrow$~(ii)
consider an amalgamation chain 
$\bigoplus_{i=1}^n (\CG[\alpha_i],g_i)$ that is not 
injectively mapped into $\CG$ by the natural homomorphism. 
For an $\ell$ that is minimal with the property that 
$1 \leq k <  k + \ell \leq N$ and that the subchain $\bigoplus_{i=k}^\ell
(\CG[\alpha_i],g_i)$ is not injectively embedded, 
its homomorphic image in $\CG$ constitutes a coset cycle.

\medskip
We are interested in the structure of $\beta$-components of
amalgamation chains
towards compatibility guarantees in the inductive elimination of 
coset cycles over $\Gamma_k(\G)$. For the simple situation 
of just a binary free amalgam $\CG[\alpha_1] \oplus \CG[\alpha_2]$,
with overlap $\CG[\alpha_0]$ for $\alpha_0 := \alpha_1
\cap \alpha_2$, the following is just a special case of 
Observation~\ref{2acyccompobs}. 

\bL
\label{comptwoamalglem}
Let $\alpha_1,\alpha_2,\beta \strictsubset E$, $\G$ an $\E$-group.
If $\G[\alpha_1]$ and $\G[\alpha_2]$ are $2$-acyclic, then 
the $\beta$-connected components of vertices in
$\CG[\alpha_1] \oplus \CG[\alpha_2]$ 
are either isomorphic to one of the 
$\CG[\beta \cap \alpha_i]$ or to a free amalgam of the form 
$\CG[\beta \cap \alpha_1] \oplus \CG[\beta \cap \alpha_2]$. 
\eL

\bL
\label{compamalgchainlem}
Let $\alpha_i \strictsubset E$ for $1\leq i\leq N$, $\beta \subset E$, 
and let the $\E$-group $\G$ and elements $g_i \in \G[\alpha_i]$ 
be such that the free amalgamation chain 
$\bigoplus_{i=1}^N (\CG[\alpha_i],g_i)$ is defined. 
If the $\G[\alpha_i]$ are $2$-acyclic, then 
the $\beta$-connected components of vertices in the $\E$-graph
$\bigoplus_{i=1}^N (\CG[\alpha_i],g_i)$
are isomorphic to  free amalgamation chains of the form 
$\bigoplus_{i=s}^{t} (\CG[\beta_i],h_i)$ for
$\beta_i := \beta \cap \alpha_i$, 
some $1 \leq s \leq t \leq N$, and suitable choices 
of $h_i \in \G[\alpha_i]$.
\eL

\bC
\label{chainconservecor}
For any $\E$-group $\G$, the functor that maps 
families of $2$-acyclic $\CG[\alpha_i]$ for $\alpha_i \subset E$ 
to free amalgamation chains of length up to~$N$ is conservative in the sense of 
Definition~\ref{conservedef}. 
\eC

\prf[Proof of Lemma~\ref{compamalgchainlem}.]
The claim follows by induction, through repeated application of 
Observation~\ref{2acyccompobs}, based on the 
auxiliary claim that any $\beta$-component of 
$\bigoplus_{i=1}^n (\CG[\alpha_i],g_i)$ for $n < N$ can intersect 
the $(\alpha_n\cap \alpha_{n+1})$-component of the vertex 
representing $g_n \in \CG[\alpha_n]$ in $\bigoplus_{i=1}^n
(\CG[\alpha_i],g_i)$ (i.e.\ the overlap in the next amalgamation step) 
in at most a single $(\beta \cap \alpha_n\cap \alpha_{n+1})$-component 
(i.e.\ in  a single $\beta$-component of that overlap). The first
non-trivial instance, for the minimal $n$ for which the $\beta$-component
in question intersects the constituent $\CG[\alpha_n]$-copy, uses
$2$-acyclicity of that $\G[\alpha_n]$; the induction step from 
$n$ to $n+1$ similarly uses $2$-acyclicity of the next
$\CG[\alpha_i]$-copy in line, i.e.\ of $\G[\alpha_{n+1}]$.
\eprf

\subsection{Amalgamation clusters} 
\label{amalgsec.clusters}

While amalgamation chains reflect the free linear or path-like unfolding pattern 
of overlaps between cosets, amalgamation clusters reflect the overlap pattern
of cosets w.r.t.\ their branching in individual elements.  
The following definition 
of free amalgamation clusters
$\bigoplus_{i \in I} \CG[\alpha_i]$ for a collection of
generator sets $\alpha_i \subset E$ is similar in spirit
to Definition~\ref{twoamalggendef} in capturing the idea of a
minimal identification. 

\bD[free amalgamation cluster]
\label{amalgclusterdef}
\mbox{}\\
The \emph{free amalgamation cluster} 
$\bigoplus_{i \in I} \CG[\alpha_i]$ of Cayley graphs
stemming from subgroups 
$\G[\alpha_i] \subset \G$ of the same $\E$-group $\G$ is obtained 
from the disjoint union of these $\CG[\alpha_i]$ (as $\alpha_i$-graphs)
by identification of vertices from $\CG[\alpha_i]$ and $\CG[\alpha_{j}]$ 
precisely if they stem from $\G[\alpha_i \cap \alpha_{j}]$.
We speak of a \emph{free amalgamation cluster over
$\Gamma_k(\CG)$} if $\alpha_i \in \Gamma_k$ for
all $i \in I$. 
\eD

This definition does not immediately reduce to iterated
application of the binary free amalgamation operation 
in the most general case. It does reduce to that, though, and with a
result that is independent of the order and precedence of 
pairwise amalgamation steps, if the
constituents are $2$-acyclic, 
which will always be the case in further applications.
If $\G$ itself is $2$-acyclic, then it embeds all free amalgamation
clusters as weak substructures, and as induced substructures if it is
$3$-acyclic. This is another generalisation of Observation~\ref{twoamalgobs}.
 
\bO
\label{amalgclusterobs} 
If $\Ghat \succeq_k \G$ is at least $3$-acyclic w.r.t.\
$\Gamma_k$-generated cosets, i.e.\
admits no $2$- or $3$-cycles of $\Gamma_k$-generated cosets,
then any free amalgamation cluster over
$\Gamma_k(\CG)$ embeds as an induced substructure into $\CGhat$.%
\footnote{As always, the embedding is unique up to arbitrary choice of
  the image of any one vertex of the connected $\E$-graph to be embedded.}
\eO

\prf
Existence of an (essentially unique) embedding of each constituent
coset in the amalgamation cluster as a weak substructure
follows from the identity $\Gamma_k(\Ghat) = \Gamma_k(\G)$. 
That a free amalgamation cluster embeds as a weak substructure
is guaranteed by $2$-acyclicity of $\Ghat$ w.r.t.\
$\Gamma_k$-generated cosets: intersections of $\Gamma_k$-generated
cosets in $\Ghat$ are the cosets generated by the common generators,
as in free amalgams. That the embedding is onto an induced rather than
just a weak substructure follows from $3$-acyclicity 
of $\Ghat$ w.r.t.\ $\Gamma_k$-generated cosets:
if images of two vertices from the same constituent $\CG[\alpha]$ of the
amalgamation cluster where related by an $e$-edge for $e \not\in
\alpha$ then this configuration would establish a $2$-cycle of
an $\alpha$-generated and an $\{ e\}$-generated coset;
if images of two vertices from the symmetric difference of the images
of two distinct constituents $\CG[\alpha_i]$
where related by an $e$-edge for $e \not\in
\alpha_1 \cup \alpha_2$,  
this would establish a $3$-cycle formed by the 
$\alpha_i$-generated cosets and an $\{ e\}$-generated coset.
\eprf

\bC
\label{embedcor}
Let $\G$ admit a $3$-acyclic $\Ghat \succeq_k \G$.%
\footnote{We shall see in Section~\ref{acycsec} that this is guaranteed if  
$\G$ is compatible with the Cayley graphs of all its
$\Gamma_k$-generated subgroups, which by Lemma~\ref{3acyclem} 
implies that $\G$ itself does not admit $2$- or $3$-cycles of
$\Gamma_k$-generated cosets; in particular all $\Gamma_k$-generated
subgroups must be $3$-acyclic.}
Then any free amalgamation cluster over
$\Gamma_k(\CG)$ embeds into the Cayley graph $\CGhat$ of $\Ghat$ 
as an induced substructure.
\eC

Again, we are interested in the consequences of this observation for
connectivity phenomena within free amalgamation clusters over
$\Gamma_k(\CG)$. In essence,  if $\G$ is compatible with all
$\E$-graphs in $\Gamma_k(\CG)$ (i.e.\ of $\Gamma_k$-generated
subgroups), then internal connectivity patterns 
w.r.t.\ generator sets from $\Gamma_k$ are strongly constrained by
(in fact largely conform to) connectivity patterns in any $3$-acyclic
surrounding $\CGhat$, as provided by an embedding
according to Corollary~\ref{embedcor}.

\bL
\label{basicclusterlemone}
Let $\G$ admit a $3$-acyclic $\Ghat \succeq_k \G$.%
\addtocounter{footnote}{-1}\footnotemark\ 
This implies the following for any
free amalgamation cluster $\H = \bigoplus_{i \in I} \CG[\alpha_i]$ 
with $\alpha_i \in \Gamma_k$ and $\beta,\gamma \in \Gamma_k$:
\bre
\item
if $\H$ is embedded into the $3$-acyclic 
unfolding $\Ghat \succeq_k \G$ according to
Corollary~\ref{embedcor}, with vertex set $V \subset \Ghat$ for $\H$, 
then the vertex set $\beta(g) \subset V$ of the $\beta$-component 
of some $g$ in $\H$ is $V \cap g\Ghat[\beta]$.
\item 
any $\beta$-connected component of $\H$ 
is isomorphic to a free amalgamation cluster of the form 
$\H_0 = \bigoplus_{i \in J} \CG[\beta_i]$ for some $J \subset I$ and 
$\beta_i := \beta \cap \alpha_i$.
\item 
the intersection of  any $\beta$- and $\gamma$-components of $\H$ is
either empty or consists of a single $(\beta \cap \gamma)$-component of $\H$, i.e.\
$\H$ is $2$-acyclic.
\ere
\eL

\bC
\label{clusterconservecor}
The functor that maps families from $\Gamma_k(\CG)$ to free amalgamation clusters
is conservative in the sense of Definition~\ref{conservedef} in
restriction to all $\E$-groups $\G$ that admit
$3$-acyclic unfoldings $\Ghat \succeq_k \G$.%
\addtocounter{footnote}{-1}\footnotemark\ 
\eC

\prf[Proof of the lemma.]
We think of $\H$ as embedded as $\H \subset
\CGhat$ into the Cayley graph of a $3$-acyclic $\E$-group $\Ghat$
according to Corollary~\ref{embedcor}, rooted at $1 \in \Ghat$ say.
Let $B \subset V$ be the vertex set of the $\beta$-component
of $\H$ in question. All $\alpha_i$- and $\beta_i$-cosets in question can be 
equally regarded as cosets in $\Ghat$ as well as in $\G$, since 
$\G \preceq_k \Ghat$.
Assertions~(i) and~(ii) now follow from $3$-acyclicity of $\Ghat$ by use 
of Observation~\ref{3acycobs}. 
Obviously $B \subset g \Ghat[\beta]$ since $\H \subsetw
\CGhat$. Let $J \subset I$ be the subset of those $i \in I$ for which 
the coset $g \Ghat[\beta]$ intersects $\CG[\alpha_i]$. Then 
Observation~\ref{3acycobs} implies that $g \Ghat[\beta] \cap
\bigcap_{i \in J}\Ghat[\alpha_i] \not= \emptyset$, whence all these 
$\Ghat[\alpha_i]$ intersect within $g \Ghat[\beta]$. For $h \in
g \Ghat[\beta] \cap \bigcap_J\Ghat[\alpha_i]$ clearly $h \in B$
and $B = \bigcup_{i \in J} g_i \Ghat[\beta \cap \alpha_i]$ follows. 
The $\beta$-connected component $\H\restr B$ is therefore isomorphic
to $\bigoplus_{i\in J} \CGhat[\beta_i] \simeq \bigoplus_{i\in J}
\CG[\beta_i]$ for $\beta_i = \beta \cap \alpha_i$ as claimed in~(ii).
(iii) is a consequence of (i). 
The vertex set of a non-trivial intersection between a $\beta$-component  
and a $\gamma$-component is a connected component w.r.t.\ $\beta\cap\gamma$.
\eprf

The following slight generalisation will be useful for later constructions. 

\bL
\label{basicclusterlemtwo}
Let $\G$ admit a $3$-acyclic $\Ghat \succeq_k \G$.%
\addtocounter{footnote}{-1}\footnotemark\ 
Then all free amalgams of the form 
$\H,g \oplus \CG[\beta]$ between a free amalgamation cluster 
$\H = \bigoplus_{i \in I} \CG[\alpha_i]$ over $\Gamma_k(\CG)$ 
and $\CG[\beta]$ from $\Gamma_k(\CG)$ are well defined and embed into 
$\CGhat$ as
weak substructures:%
\footnote{The example of $\H,g := \CG[e_1] \oplus \CG[e_2], e_2$
  with commuting generators $e_1,e_2$ and $\beta := \{ e_1 \}$ shows
  that, without higher acyclicity requirements, 
  $\H,g \oplus \CG[\beta]$ may not embed as an induced substructure.}
$\H,g \oplus \CG[\beta]\subsetw
\CGhat$. 
Moreover, for $\gamma \in \Gamma_k$, any 
$\gamma$-component of free amalgams $\H,v \oplus
\CG[\beta]$ is isomorphic 
to $\CG[\gamma]$ or to a free amalgam of the form $\bigoplus_{i\in J}
\CG[\gamma_i], g' \oplus \CG[\gamma']$ for 
$J \subset I$, $\gamma_i = \gamma \cap \alpha_i$ and $\gamma' \subset 
\gamma \cap \beta$ (possibly empty).
\eL

\prf
The claims about existence and embeddability follow directly
from the previous lemma. For the claim about $\gamma$-components
of $\H,g \oplus \CG[\beta]$, we may use $2$-acyclicity of $\H$
(by~(iii) from the previous lemma) and of $\CG[\beta]$ 
to apply Observation~\ref{2acyccompobs}.
\eprf

\bC
\label{clusterplusextraconservecor}
The functor that maps  
families from $\Gamma_k(\CG)$ to free amalgams between
free amalgamation clusters over $\Gamma_k(\CG)$
and one individual $\CG[\alpha]$ from $\Gamma_k(\CG)$
is conservative in the sense of Definition~\ref{conservedef} in
restriction to all $\E$-groups $\G$ that admit
$3$-acyclic unfoldings $\Ghat \succeq_k \G$.%
\addtocounter{footnote}{-1}\footnotemark\ 
\eC

\section{Construction of N-acyclic E-groups}
\label{acycsec}

We may construct $N$-acyclic $\E$-groups based 
on the general considerations in Proposition~\ref{conserveprop}.
The functor $\F$, which collects the obstructions at a given level
$\Gamma_m$ into an $\E$-graph $\H$ that eliminates these obstructions 
at level $\Gamma_{m+1}$, just builds amalgamation chains of appropriate
length over $\Gamma_m(\CG)$. Its conservative nature was established
in Section~\ref{amalgsec.chains}, Lemma~\ref{compamalgchainlem}.

In essence we find, for fixed $N \geq 2$, a sequence of $\E$-groups
$\G_m$ such that 
\[
\G_m \mbox{ is $N$-acyclic w.r.t.\ } \Gamma_m(\G_m),
\]
i.e.\ $\G_m$ admits no coset cycles of length up to $N$ 
with constituents $\Gamma_m(\CG_m)$ ($\alpha$-cosets for $|\alpha| <
m$). In order to preserve the acyclicity property at lower levels in 
the induction step, we pass from $\G_m$ to $\G_m \preceq_m \G_{m+1}$ 
obtaining $\G_{m+1}$ as $\G_{m+1} := \sym(\H_m)$, with 
$\H_m$ the  disjoint union of $\CG_m$ and $\F(\Gamma_m(CG_m))$.

The following determines the appropriate length bound, which 
generalises Lemma~\ref{3acyclem}.

\bL
\label{fromcompattoacyclem}
For $n \geq 1$, let $\G[\alpha]$ be compatible with all free
amalgamation chains of length up to~$n$ 
with constituents $\CG[\alpha']$ for $\alpha' \strictsubset \alpha$.
Then $\G[\alpha]$ is $N$-acyclic for $N = n+2$.
\eL

\prf
The gist of the matter is that $\G \succeq \sym(\H)$ for every
free amalgamation chain $\H$ of length up to $n$. 
These free amalgamation chains unfold potential coset cycles. This rules out
corresponding cycles in $\sym(\H)$ and hence in $\G$ as follows. 
Suppose the pointed cosets $(g_i \G[\alpha_i],g_i)_{i \in \Z_m}$ for
$m \leq N$
formed a coset cycle in $\G$.  
Similar to the argument in Lemma~\ref{3acyclem}, we think of cutting
the cycle (this time in $g_0 = g_m$) and test the permutation group action
on a chain formed by the remaining $m-2$ links,  viz.\ on the free amalgamation
chain
\[
\H = \bigoplus_{i=1}^{m-2} (\CG[\alpha_i],g_i^{-1}g_{i+1})
\] 
of length $m-2 \leq n$. Let $\alpha_{i,j} := \alpha_i \cap \alpha_j$, 
and let $w_i \in \alpha_i^\ast$ be such that $[w_i]_\G = g_i^{-1} g_{i+1}$. 
For the links from and to $g_0$ in the cycle, 
$w_0 \in \alpha_0^\ast$ and $w_{m-1} \in \alpha_{m-1}^\ast$,
we also look at their projections on
the neighbouring constituents in $\H$. Let 
$w_{0,1} \in \alpha_{0,1}^\ast$ be the projection of $w_0$ to $\alpha_1$ 
and 
$w_{m-2,m-1} \in \alpha_{m-2,m-1}^\ast$ the projection of $w_{m-1}$ to $\alpha_{m-2}$.
Let $v$ be the element of the chain $\H$ that corresponds to
$[w_{0,1}]_\G^{-1}$ in its first constituent $\CG[\alpha_1]$, 
and consider the permutation group action of $[\prod_{i=0}^{m-1}
w_i]_\H \in \sym(\H)$,
on $v$ in $\H$. By the separation condition
for the $\alpha_0$-link of the cycle, the generator sequence 
$w_0$ has the same effect on $v$ as its projection $w_{0,1}$ and
maps $v$ to the element corresponding to $1$ in the first constituent 
$\CG[\alpha_1]$ of $\H$; the connectivity and separation conditions
for the $\alpha_i$-links up to $i =m-2$ imply that 
$[\prod_{i= 0}^{m-1} w_i]_\H$ maps $v$ to 
the element corresponding to $g_{m-2}^{-1}g_{m-1}$ in the
last constituent $\CG[\alpha_{m-2}]$ of $\H$; and the separation 
condition for the $\alpha_{m-1}$-link of the cycle shows that 
the final image of $v$ is an element in the $\alpha_{m-2,m-1}$-component 
of $g_{m-2}^{-1}g_{m-1}$ in that $\CG[\alpha_{m-2}]$ constituent of $\H$. But
this image is necessarily distinct from $v$, contradicting compatibility
of $\G$ with $\H$, as $[\prod_{i=0}^{m-1} w_i]_\G 
=\prod_{i= 0}^{m-1} (g_i^{-1} g_{i+1}) = 1$.
\eprf 

\bP
\label{inductiveNacycprop}
Let $n \geq 1$ and let $\F$ be the functor that maps an $\E$-group 
$\G$ to the collection of all free amalgamation chains of length up
to~$n$ over constituents $\CG[\alpha] \subsetw \CG$ for $\alpha
\subset E$. For fixed $1 \leq k \leq |E|+1$ consider an $\E$-group 
$\G$ that is compatible with the $\E$-graphs in $\F(\Gamma_k(\G))$.%
\footnote{Note that this is vacuously true at level $k=1$.}
Then $\G$ admits an unfolding $\Ghat \succeq_k \G$
that is $N$-acyclic for $N= n+2$.
\eP

\prf
This is a direct application of Proposition~\ref{conserveprop}
to the functor $\F$, which is conservative according 
to Lemma~\ref{compamalgchainlem}.
\eprf

Recall from Definition~\ref{Esymmdef} the notion of 
symmetries for $\E$-groups and $\E$-graphs that are induced by
permutations of the set $E$. 
It is clear from the construction steps in
Proposition~\ref{inductiveNacycprop}
above that they do not 
break any such symmetry in the passage from $\G_k$ to $\G_{k+1}$.
It follows that the passage from $\G$ to an $N$-acyclic $\Ghat \succeq
\G$ preserves all symmetries of $\G$. We state  this additional
feature in the theorem which otherwise just sums up the outcome 
of the proposition.

\bT
\label{symmacycEgroupthm}
For every finite $\E$-group $\G$ and $N \geq 2$ there is 
a finite $\E$-group $\Ghat \succeq \G$ that is $N$-acyclic and 
fully symmetric over $\G$ in the sense that every permutation 
of the generator set $E$ that is a symmetry of $\G$ extends to a
symmetry of $\Ghat$: $\G^\rho \simeq \G \Rightarrow \Ghat^\rho \simeq \Ghat$.
\eT

In particular, we obtain finite $N$-acyclic $\E$-groups $\Ghat$ that
are fully symmetric, admitting every permutation of $E$ as a
symmetry, if we start from the fully symmetric $\G := \sym(\H)$
for the hypercube $\H = 2^E$.

\section{Constraints on generator sequences}
\label{constraintsec}

The building blocks of plain coset cycles are 
generated subgroups of the form $\G[\alpha] \subset \G$, which may
also be seen as the images of $\alpha^\ast \subset E^\ast$ 
under the natural homomorphism
\[
\barr{rcl}
[\,\nt\,]_{\G} \colon E^\ast &\longrightarrow& \G 
\\
w = e_1\cdots e_n & \longmapsto & [w]_{\G} := \prod_{i=1}^n e_i = e_1 \cdots e_n
\earr
\]
that associates a group element with any (reduced) word over $E$.
This association naturally translates to cosets $g \G[\alpha]$. 
Alternatively, $\G[\alpha]$ and 
$g \G[\alpha]$ may be regarded as the $\alpha$-connected
components of $1$ or $g$ in the Cayley graph $\CG$ of $\G$.
A natural way of putting extra constraints on these weak subgraphs, 
with reasonable closure properties in terms of 
generator sets $\alpha \subset E$, is the following. Consider a 
fixed $\E$-graph $\I = (S,(R_e)_{e \in E})$ on vertex set $S$. 
We regard $\I$ as a template for systematic restrictions on patterns
of generator sequences, calling it a \emph{constraint graph}.

\bProv
\label{constraintpatternprov}
We fix an $\E$-graph $\I = (S,(R_e)_{e \in E})$ as a 
\emph{constraint graph}, and consider
only $\E$-groups $\G$ that are compatible with $\I$, 
i.e.\ with $\G \succeq \sym(\I)$. 
\eProv

\bR
\label{compatconstraintrem}
The restriction to $\E$-groups that are compatible with $\I$, 
$\G \succeq \sym(\I)$, does imply that there is a well-defined group
action of $\G$ on $\I$. But w.r.t.\ this group action, any $s \in S$
that is not incident with an $e$-edge is a fixed point of 
$\pi_e \in \sym(\I)$. As $\I$ will typically not be a complete
$\E$-graph, the Cayley graph of $\CG$ does not map 
homomorphically to $\I$. But $\alpha$-walks from $s \in \I$ do have 
unique lifts to $\alpha$-walks from any $g \in \CG$. Compatibility says 
that lifts at the same $g \in \CG$ of different walks from 
$s \in S$ can only meet in $\CG$ above positions in which the 
given walks meet in $\I$.
\eR

The idea is to regard $\I$ as a template for edge patterns of walks 
and correspondingly restricted notions of reachability and connected 
components. Recall from Definition~\ref{alphawalksconnectdef} 
the notion of $\alpha$-walks in $\E$-graphs, which we shall
now refine in connection with the constraint graph $\I$.
According to Definition~\ref{alphawalksconnectdef}, a 
\emph{walk} in an $\E$-graph $\H$ 
is a sequence of vertices and edge labels
of the form   
\[
s_0,e_1,s_1,e_2,\ldots, s_{n-1},e_n,s_n
\]
with $s_i \in S$ and $e_i \in E$ where 
$(s_i,s_{i+1}) \in R_{e_{i+1}}$ for $0 \leq i < n$ for some 
$n \in \N$. The edges $(s_i,s_{i+1}) \in R_{e_{i+1}}$ are the edges 
traversed by this walk. The above walk is a walk of length~$n$ from 
the source $s= s_0$ to the target $t= s_n$
and its \emph{edge label sequence} is the word 
$w = e_1\cdots e_n \in E^\ast$. We also say that this word $w$
\emph{labels a walk} (of length $n$, from $s_0$ to $s_n$) in $\H$,
and describe the walk in question as a $w$-walk, or as an $X$-walk 
for some language $X \subset E^\ast$ if $w \in X$. In the case of 
$X = \alpha^\ast$ we also speak of $\alpha$-walks instead of 
$\alpha^\ast$-walks.

Note that a word $w \in E^\ast$ can label at most one walk 
from a given vertex $v$ in any $\E$-graph $\H$. That all $w \in E^\ast$ label walks 
from all $v \in \H$ is equivalent to $\H$ being complete
(i.e.\ to each $R_e$ being a full matching).

\bD[$\I$-words and $\I$-walks]
\label{Ialphawalksdef}
\mbox{}\\
For $\alpha \subset E$ and $s \in S$ let $\I[\alpha,s] \subsetw \I$ 
be the weak substructure of $\I = (S,(R_e)_{e \in E})$ whose universe
is the $\alpha$-connected component of $s$ (i.e.\ the connected
component of $s$ in the $\alpha$-reduct $\I\restr \alpha$), and 
with induced $R_e$ for all $e \in \alpha$. 
Natural sets of $E$- or $\alpha$-words that occur as edge label sequences
along walks on the constraint graph $\I$ are defined as follows:
\bre
\item[--]
$\alpha^\ast[\I] \subset E^\ast$ consists of those 
$w \in \alpha^\ast$ that label a walk in $\I$;
\item[--]
$\alpha^\ast[\I,s] \subset E^\ast$
consists of those $w \in \alpha^\ast$ that label a walk 
from $s$ in $\I$;
\item[--]
$\alpha^\ast[\I,s,t] \subset E^\ast$
consists of those $w \in \alpha^\ast$ that label a walk 
from $s$ to $t$ in $\I$.
\ere
\eD

In the following, $\alpha^\ast[\I]$-walks will also be addressed as
$\I[\alpha]$-walks, especially if the restriction of the
available generators to $\alpha \subset E$ seems more important than the
actual labelling.

\bD[$\I[\alpha,s]$-component]
\label{Ialphacompdef}
\mbox{}\\
For an $\E$-graph $\H$, $\alpha \subset E$ and $s \in S$, 
the \emph{$\I[\alpha,s]$-component} of a vertex $v \in \H$
is the following weak substructure $\H[\I,\alpha,s;v] \subsetw \H$:
its vertex set consists of those vertices that are reachable 
on an $\alpha^\ast[\I,s]$-walk from $v$ (i.e.\ on a walk in 
$\H$ whose edge label sequence is in $\alpha^\ast[\I,s]$); 
its edge relations comprise those $R_e$-edges for $e \in \alpha$
that are traversed by $\alpha^\ast[\I,s]$-walks from $v$ in $\H$. 
\eD

When we speak of an $\I[\alpha]$-component we mean an 
$\I[\alpha,s]$-component for some $s \in S$, which is left
unspecified, but keep in mind that, e.g.\ $\I[\alpha]$-walks inside
such $\I[\alpha]$-components implicitly refer to a fixed choice of 
an anchor point that determines which edges are available. 
For later use we define a direct product of $\I$ with the Cayley graph
of an $\E$-group, which reflects $\I$-reachability, as follows.

\bD[direct product]
\label{directproddef}
\mbox{}\\
Let $\G$ be an $\E$-group 
that is compatible with 
the constraint graph $\I= (S,E)$, $\CG$ the Cayley graph of $\G$.
Then the \emph{direct product} $\I \otimes \CG$ is the
$\E$-graph 
\[
\I \otimes \CG = (V,(R_e)_{e \in E})
\]
with vertex set $V = S \times G$ and edge relations
\[
R_e = \{ ((s,g),(s',g')) \colon (s,s') \in R_e \mbox{ in } \I \mbox{
  and } g' = ge \mbox{ in } \G\} \quad \mbox {for $e \in E$.}
\]
\eD

Note that all $\alpha$-walks in $\I \otimes \CG$, by definition of the
edge relations, trace the lifts of $\alpha$-walks in $\I$. 
For the following also compare Remark~\ref{compatconstraintrem}
above on lifts of walks from $\I$ to $\CG$.

\bO
\label{directprodobs}
Compatibility of $\G$ with $\I$ implies that for any $(s,g)$ and $(s',g)$ in the same
connected component of $\I \otimes \CG$ we must have $s = s'$.
It follows that connected components of $\I \otimes \CG$ are isomorphic
to weak subgraphs of $\CG$, and that the connected components of $\I \otimes \CG$ 
reflect $\I$-reachability in $\CG$ in the sense that 
$(s,g)$ and $(t,g')$ are in the same $\alpha$-connected 
component of $\I \otimes \CG$ if, and only if, 
$g' \in \CG[\I,\alpha,s;g]$. 
\eO

We turn to cosets, coset cycles and acyclicity criteria relative to
the given constraint graph $\I$. 

\bD[$\I$-coset]
\label{Icosetdef} 
\mbox{}\\
In an $\E$-group $\G$ 
that is compatible with $\I$ 
the \emph{$\I[\alpha,s]$-coset} of $g \in \G$
is the $\I[\alpha,s]$-component $\CG[\I,\alpha,s;g] \subsetw \CG$ 
of $g$. 
\eD

We drop mention of $g$ for cosets at $g = 1$, writing e.g.\ just 
$\CG[\I,\alpha,s] \subsetw \CG$ for $\CG[\I,\alpha,s;1] \subsetw \CG$. 
For the following compare Definition~\ref{cosetcycledef}
for plain coset cycles.

\bD[$\I$-coset cycle]
\label{Icosetcycledef} 
\mbox{}\\
Let $\G$ be an $\E$-group, $n \geq 2$.
An \emph{$\I$-coset cycle} of length $n \geq 2$ in $\G$ or $\CG$ is a cyclically
indexed sequence $(\CG[\I,\alpha_i,s_i;g_i],g_i)_{i \in \Z_n}$ of pointed
$\I[\alpha_i,s_i]$-cosets for $\alpha_i \strictsubset E$, $s_i \in S$ and 
$g_i \in G$ satisfying these conditions:
\bre
\item
(connectivity)
there is an $\alpha_i^\ast[\I,s_i, s_{i+1}]$-walk from $g_i$ to
$g_{i+1}$,
\footnote{Equivalently, $g_{i+1} \in \CG[\I,\alpha_i,s_i;g_i]$, or
$\CG[\I,\alpha_i,s_i;g_i] = \CG[\I,\alpha_i,s_{i+1};g_{i+1}]$}
\item
(separation)
$\CG[\I,\alpha_{i,i-1},s_i;g_i] \cap \CG[\I,\alpha_{i,i+1}, s_{i+1}; g_{i+1}] = \emptyset$,
\ere
where $\alpha_{i,j} := \alpha_i \cap \alpha_j$.
\eD 

\bD[$N$-acyclicity over $\I$]
\label{INacycdef}
\mbox{}\\
An $\E$-group 
that is compatible with $\I$ 
is called \emph{$N$-acyclic
over $\I$} if it does not admit any $\I$-coset cycles of length up to $N$. 
\eD

The above definitions generalise corresponding definitions 
in the unconstrained setting. Those definitions, 
Definitions~\ref{cosetcycledef} and~\ref{Nacycdef}, 
are subsumed in the above as special cases for
the trivial constraint graph having a single vertex with loops for all $e \in E$.

As in the unconstrained case, $2$-acyclicity over $\I$ is akin 
to a notion of simple connectivity, being equivalent to the
requirement that, for all $\alpha_i$ and $s$, 
\[
\CG[\I,\alpha_1,s] \cap \CG[\I,\alpha_2, s]
 = \CG[\I,\alpha_1\cap \alpha_2, s].
\]

In order to boost degrees of plain coset acyclicity to acyclicity 
w.r.t.\ $\I$-cosets we aim to employ compatibility with unfoldings
of potential $\I$-coset cycles in a manner similar to the treatment
of Section~\ref{acycsec}. The idea is to unfold potential $\I$-coset cycles
with constituents $\CG[\I,\alpha_i,s_i,g_i] \subsetw \CG[\alpha_i]$   
through the unfolding of the surrounding plain cosets $\CG[\alpha_i]$
into free amalgamation chains. This is the same unfolding into
amalgamation chains as discussed in Section~\ref{amalgsec}
and used towards plain coset acyclicity in Section~\ref{acycsec}.
The new challenge is to safeguard this unfolding against unwanted 
shortcuts by $\alpha$-walks that do not correspond to 
$\I[\alpha]$-walks.
At the same time, the consecutive unfoldings need to preserve
compatibility with smaller $\alpha$ in order to render the process conservative
w.r.t.\ to what has been achieved in earlier stages. 

It is to this end that we consider the relationship between 
$\I[\alpha]$-components and the halo around them that is
formed by overlapping (unconstrained) $\alpha'$-connected components for 
$\alpha' \strictsubset \alpha$. The \emph{$\I$-skeletons} to be
defined below play the r\^ole of $\I[\alpha]$-connected 
components that unfold all those $\alpha$-links that are available 
according to $\I$, and none else. These arise naturally as 
$\I[\alpha]$-components of Cayley graphs $\CG$ provided 
$\G$ is compatible with $\I$, 
as we assumed with Proviso~\ref{constraintpatternprov}. But we also
consider individual such $\I$-skeletons as scaffoldings for the
construction of further $\E$-graphs, 
which are to be fed into the unfolding process that successively 
unclutters the embedded $\I[\alpha]$-components and unfolds short
$\I$-coset cycles for increasing numbers of generators.

\bD[$\I$-skeleton]
\label{Ialphaskeletondef}
\mbox{}\\
For $\alpha \subset E$ an \emph{$\I[\alpha,s]$-skeleton} is a
connected $\alpha$-graph $\H$ that admits a surjective homomorphism 
$h \colon \H \rightarrow \I[\alpha,s]$ onto
an $\alpha$-connected component $\I[\alpha,s]$ of $\I$
with the following lifting property: 
whenever $h(v) \in S$ is incident with an $e$-edge in $\I[\alpha,s]$
then so is $v$ in $\H$.
Across all $s \in S$ we speak of just $\I[\alpha]$-skeletons. 
\eD

When speaking of an $\I[\alpha]$-skeleton, it is important to keep in
mind that it is an $\I[\alpha,s]$-skeleton for some $s \in S$.

\bE
\label{skeletobex}
Every $\CG[\I,\alpha,s;g] \subsetw \CG[\alpha] \subsetw \CG$ 
is an $\I[\alpha,s]$-skeleton, via  
the unique homomorphism $h$ that maps $g \in \G$ to $s \in S$.
This map is well-defined since $\G$ is assumed compatible with $\I$ 
(cf.\ Remark~\ref{compatconstraintrem} 
and comments after Definition~\ref{directproddef}). 
We refer to such instances as \emph{embedded $\I$-skeletons}.
\eE

The lifting property in the definition guarantees that every
$\I[\alpha]$-walk from $s'$ in $\I[\alpha,s]$ has a unique lift to
an $\alpha^\ast[\I]$-walk from $v$ for every $v$ in the pre-image of $s'$.
In alternative terminology, the homomorphism establishes 
$\H$ as an unbranched bisimilar cover of $\I[\alpha,s]$ (or of $\I\restr
\alpha$).

\bProv
\label{lazyShprov}
In order not to overburden notation, we often leave implicit the
$S$-marking of the elements of an $\I[\alpha,s]$-skeleton $\H$ 
that comes with a given homomorphism $h \colon \H \rightarrow \I[\alpha,s]$.
Speaking of an $\I[\alpha]$-walk in $\H$, for instance, we only
admit walks labelled by $w \in \alpha^\ast[\I,s,t]$ from 
$v$ to $v'$ in $\H$ if these are marked accordingly, i.e.\ for
$h(v) = s$, which then implies $h(v') = t$.
\eProv

\bD[freeness]
\label{Gfreeoverskeletdef}
\mbox{}\\
The embedded $\I[\alpha]$-skeleton $\H := \CG[\I,\alpha,s;g] \subsetw
\CG[\alpha] \subsetw \CG$
is \emph{free} in $\CG[\alpha]$ or in $\CG$ if for any two 
elements $v_1,v_2 \in \H$, any cosets 
$v_1 \G[\alpha_1]$ and $v_2 \G[\alpha_2]$,
for $\alpha_i \strictsubset \alpha$, overlap in the surrounding 
$\CG$, if, and \emph{only if}, the $\alpha_i$-connected components 
of these elements $v_i$ overlap within $\H$:
 \[
\alpha_1(v_1)\, \cap\, \alpha_2(v_2) = \emptyset
   \mbox{ in $\H$ } \; \Longrightarrow \; 
v_1 \G[\alpha_1]\cap v_2 \G[\alpha_2]
= \emptyset
   \mbox{  in $\CG[\alpha]$.} 
 \]
The $\E$-group $\G$ is \emph{free over $\I$} if every 
embedded $\I[\alpha]$-skeleton for $\alpha \subset E$ is free in $\CG$.
\eD

Already the special case of $\alpha' := \alpha_1 = \alpha_2 \strictsubset \alpha$
gives an indication of the strength of this freeness condition:
the freeness requirement here says that $\alpha'$-reachability inside the embedded 
$\I[\alpha]$-skeleton agrees with $\alpha'$-reachability in the surrounding $\CG$. 
The following is an analogue of Lemma~\ref{acychomlem} for freeness. 

\bL
\label{freehomlem} 
If $\G[\alpha]$ is $2$-acyclic, then freeness of $\G[\alpha]$ 
over $\I$ is preserved under inverse homomorphisms that are injective on
$\alpha'$-generated subgroups for $\alpha'\strictsubset \alpha$.
\eL

\prf
Let $h \colon \Ghat \rightarrow \G$ be a homomorphism of $\E$-groups
(both assumed to be compatible with $\I$:
Proviso~\ref{constraintpatternprov}), its restrictions to members of
the family $\Gamma(\Ghat[\alpha])$ injective. 
The image of an embedded $\I[\alpha,s]$-skeleton in 
$\CGhat$ under $h$ is an embedded $\I[\alpha,s]$-skeleton in $\CG$.
Assume that an $\I[\alpha,s]$-skeleton $\Hhat \subsetw \CGhat$
violates the freeness condition through some 
$\hat{g} \in v_1 \Ghat[\alpha_1] \cap v_2 \Ghat[\alpha_2]\not= \emptyset$ where
$\CGhat[\I,\alpha_1,s_1; v_1] \cap \CGhat[\I,\alpha_2,s_2;v_2] = \emptyset$.
By assumption $h$ is injective in restriction to each one of the two cosets 
$v_i\Ghat[\alpha_i]$ and $g_2 \Ghat[\alpha_2]$. By
Lemma~\ref{acychomlem} $\Ghat[\alpha]$ is
$2$-acyclic so that these two cosets intersect precisely in the coset 
$\hat{g} \Ghat[\alpha_0]$ for $\alpha_0 = \alpha_1 \cap \alpha_2$. By
$2$-acyclicity of $\G[\alpha]$, $h$ must be injective in restriction to the union 
of these overlapping cosets. Therefore the violation of freeness in
$\CGhat[\alpha]$ would be isomorphically mapped onto a violation in $\G[\alpha]$. 
\eprf

An immediate observation relates freeness to reachability 
notions that are crucial in connection with the separation condition 
for coset cycles. To make the connection with
Lemma~\ref{cosettocosetoverIlem} below, consider  
$\alpha,\alpha_1,\alpha_2$ to play the r\^oles
of $\alpha_i, \alpha_{i,i-1}, \alpha_{i,i+1}$ in a potential
coset cycle or $\I$-coset cycle.

\bO
\label{freenessobs2}
For $\alpha_1,\alpha_2 \strictsubset \alpha$ consider vertices
$v_1,v_2$ of an embedded $\I[\alpha]$-skeleton 
$\H := \CG[\I,\alpha,s_1;v_1] = \CG[\I,\alpha,s_2;v_2] 
\subsetw \CG[\alpha]$ that is free in $\CG$. 
If $\CG[\I,\alpha_1,s;v_1] \cap \CG[\I,\alpha_2,s';v_2] =\emptyset$
(separation as in $\I$-coset cycles)
then also $v_1 \G[\alpha_1] \cap v_2 \G[\alpha_2] =\emptyset$
(separation as in plain  coset cycles).
\eO

Finally, the condition on overlaps in Definition~\ref{Gfreeoverskeletdef}
takes an especially neat and concrete form if $\G$ is $2$-acyclic.

\bO
\label{freenessobs3}
Let the embedded $\I[\alpha]$-skeleton $\H := \CG[\I,\alpha,s;g] \subsetw \CG$
be free in $\CG[\alpha]$ and let $\G[\alpha]$ be $2$-acyclic.
Then two cosets $v_1 \G[\alpha_1]$ and $v_2 \G[\alpha_2]$
for $v_i \in \H$ and $\alpha_i \strictsubset \alpha$ 
are either disjoint or their overlap
is a single coset of the form $v_0 \G[\alpha_0]$ for
$\alpha_0 := \alpha_1\cap \alpha_2$ and an element 
$v_0 \in \H$ that is 
$\alpha_i$-reachable from $v_i$ inside the skeleton, for $i=1,2$.
\eO

\prf
Due to $2$-acyclicity a non-empty intersection of 
$\alpha_i$-cosets in $\CG[\alpha]$ is a single
$\alpha_0$-coset. By freeness it must intersect 
the embedded $\I[\alpha]$-skeleton in an element 
$v_0$ that is $\alpha_i$-reachable from $v_i$ in the
skeleton. 
\eprf

\bL
\label{cosettocosetoverIlem}
If $\G$ is $N$-acyclic (does not admit plain coset cycles of length up
to $N$) and free over $\I$, 
then $\G$ is $N$acyclic over $\I$
(does not admit $\I$-coset cycles of length up
to $N$).
\eL

\prf
It suffices to show that every $\I$-coset cycle in $\G$ induces
a plain coset cycle if the links $\CG[\I,\alpha_i,s_i;g_i]$ are
extended to the surrounding $g_i\CG[\alpha_i]$. 
The crux is the separation condition: 
separation w.r.t.\ $\I$-reachability in 
links $\CG[\I,\alpha_i,s_i;g_i]\subsetw g_i\CG[\alpha_i]$ 
needs to translate into the stronger separation condition for 
the encompassing $g_i\G[\alpha_i]$ themselves. This
is precisely the content of Observation~\ref{freenessobs2}.  
\eprf

\section{Small coset amalgams over I-skeletons}
\label{CEsec}

We think of a fixed constraint graph $\I$ and its generated subgraphs
$\I[\alpha]$, and only consider $\E$-groups $\G$ that are compatible 
with $\I$ according to Proviso~\ref{constraintpatternprov}.

For \emph{small coset amalgams} we want to amalgamate
constituent $\CG[\alpha']$-copies for $\alpha' \strictsubset \alpha$
along an $\I[\alpha]$-skeleton $\H[\alpha]$ to form an $\E$-graph,
in as free a manner as possible.
In every vertex $v$ of the skeleton, we mean to amalgamate
$v$-tagged copies $v \CG[\alpha']$ of $\CG[\alpha']$ with minimal
identifications as are enforced by the requirements for $\E$-graphs.
Before looking at embedded realisations of such amalgams in a surrounding
$\CG[\alpha]$ (with appropriate acyclicity conditions) we define these 
amalgams in isolation (with only intrinsic acyclicity conditions for
the skeleton $\H[\alpha]$ in relation to the constituents
$\CG[\alpha']$ for $\alpha' \strictsubset \alpha$). 
The pre-requisite~(iii) in the following definition 
corresponds to the freeness condition for embedded $\I$-skeletons of
Definition~\ref{Gfreeoverskeletdef}, here expressed for 
$\alpha'$-components $\H[\alpha';v] \subsetw \H[\alpha]$ for $\alpha' \strictsubset \alpha$. 

\bD[small coset amalgam]
\label{smallcosetamalgCEdef}
\mbox{}\\
An $\I[\alpha]$-skeleton $\H := \H[\alpha]$ and the family of subgroups
$\G[\alpha']$ for $\alpha' \strictsubset \alpha$ \emph{admit 
a small coset amalgam} if
\bre
\item
  $\H[\alpha]$ is $2$-acyclic;
\item 
the $\G[\alpha']$ are $2$-acyclic; 
\item
 all $\alpha'$-components of $\H[\alpha]$ embed 
 isomorphically (as $\alpha'$-graphs) into $\CG[\alpha']$ 
 in a free manner,
 i.e.\ there is an embedding $\rho \colon \H[\alpha';v] \rightarrow \CG[\alpha']$
 such that for $\alpha_1',\alpha_2'
 \strictsubset \alpha'$:
 \[
   \alpha_1'(v_1) \cap \alpha_2'(v_2) = \emptyset
   \mbox{ in $\H[\alpha';v]$ } \; \Longrightarrow \; 
   \alpha_1'(\rho(v_1)) \cap \alpha_2'(\rho(v_2)) = \emptyset
   \mbox{  in $\CG[\alpha']$.}
 \]
 \ere
Under these conditions 
the \emph{small coset amalgam}
\[
 \CE(\H[\alpha],\G,\alpha) := \bigoplus_{v \in V, \alpha'\strictsubset
   \alpha} \!\!\!\!\! v\CG[\alpha'] \;\big/\, \approx,
\]
is defined as the quotient of the disjoint union of $v$-tagged 
constituents $v\CG[\alpha']$ w.r.t.\ 
$\approx$, the equivalence relation induced by the following identifications, 
cf.\ Figure~\ref{cosetextnpatternfig} and recall
Proviso~\ref{lazyShprov} for walks in $\H[\alpha]$: 
\[
  (\dagger)
  \qquad
  \mbox{\begin{minipage}{.75\textwidth}
$v_1 g_1 \in v_1\CG[\alpha_1']$ is identified with 
$v_2 g_2 \in v_2 \CG[\alpha_2']$ if for some $v_0$ in $\H$ and
$g_0 \in \G[\alpha_0']$ for $\alpha_0' := \alpha_1'\cap\,\alpha_2'$, 
$v_0$ is reachable from $v_i$ inside 
$\H$ by an 
$\I[\alpha_i]$-walk labelled $w_{i0}$ 
such that $v_0 = v_i [w_{i0}]_{\H}$ and 
$g_{i0} g_0 := [w_{i0}]_\G g_0 = g_i$ in $\G[\alpha_i']$ for $i=1,2$.
\end{minipage}}
\]
\eD

Note that, in particular, constituents $v\CG[\alpha']$ and $v'\CG[\alpha']$ coincide whenever $v$ and $v'$ are
linked by an $\I[\alpha']$-walk in $\H[\alpha]$ ($\alpha'$-connected
in $\H[\alpha]$).

\begin{figure}
\[
\xymatrix{\raisebox{-20pt}{$v_1\CG[\alpha_1']$} &&\circ && \raisebox{-20pt}{$v_2\CG[\alpha_2']$}
\\
\\
v_1 \ar@{->}[uurr]|*+{g_1}\ar@{->}[rr]|*+{g_{10}}
&& v_{0} 
\ar@{->}[uu]|*+{g_{0}} 
&& v_2 \ar@{->}[uull]|*+{g_2}\ar@{->}[ll]|*+{g_{20}} \ar@{}[r]|*{\H[\alpha]} &
}
\]
\caption{Pattern for identifications in $\CE(\H[\alpha],\G,\alpha)$,
  see
   $(\dagger)$ in Definition~\ref{smallcosetamalgCEdef}:
  for $i=1,2$, 
$g_{i0} = [w_{i0}]_\G$ links $v_0$ to $v_i$ in their
$\alpha_i'$-component of $\H$.}
\label{cosetextnpatternfig}
\end{figure}

Pre-condition~(iii) in the definition is essential for the analysis of 
$\approx$ for identification of elements according to 
$(\dagger)$. Let $\sim$ denote the relation on 
the disjoint union of the $v$-tagged copies of constituent 
$v\CG[\alpha']$: $v_1g_1 \sim v_2g_2$ if $(\dagger)$, where we 
write $v_ig_i$ instead of $(v_i,g_i)$ for the element 
$g_i$ in the $v_i$-tagged copy of $\G[\alpha_i']$ that forms the 
constituent $v_i\CG[\alpha_i']$. Then $\approx$ is the transitive 
closure of $\sim$. But condition~(iii) entails that $\sim$ itself is
transitive so that $(\dagger)$ directly describes all identifications
in $\CE(\H[\alpha],\G,\alpha)$. Transitivity of $\sim$ is shown as follows,
cf.\ Figure~\ref{transpatternfig}.
\begin{figure}
\[
\xymatrix{ & 
&*+{\circ} & 
&
\\
&
\raisebox{-20pt}{$\qquad\qquad\qquad$} 
&
\raisebox{-20pt}{$\qquad\qquad\qquad$} 
&
\raisebox{-20pt}{$\qquad\qquad\qquad$} 
&
\\
v_1 \ar@{-}[r]|*{\alpha_1'}
\ar@{->}[uurr]|*+{g_{1}}
& v_{1,2} \ar@{->}[uur]|*+{g_{1,2}}
\ar@{-}[r]|*{\alpha_{2}'}
&  v_2 \ar@{-}[r]|*{\alpha_{2}'}
\ar@{->}[uu]|*+{g_{2}}
& v_{2,3} \ar@{->}[uul]|*+{g_{2,3}}
\ar@{-}[r]|*{\alpha_{3}'}
& v_3
\ar@{->}[uull]|*+{g_{3}} 
}
\]
\caption{Pattern for two consecutive identifications.}
\label{transpatternfig}
\end{figure}
Let $v_1g_1 \sim v_2g_2 \sim v_3g_3$
where $g_i \in \G[\alpha_i']$ is an element of the 
constituent $v_i \CG[\alpha_i']$, for 
$\alpha_i' \in \Gamma_\alpha$, i.e.\ 
$\alpha_i' \strictsubset \alpha$. 
By the definition of $\sim$ in $(\dagger)$ this implies 
the existence of intermediaries $v_{i,j}g_{i,j}$ for $(i,j) = (1,2), (2,3)$
according to $(\dagger)$. I.e.\ 
$v_{i,j}$ is in the intersection of the $\alpha_i'$-component of
$v_i$ and the $\alpha_j'$-component of $v_j$ in $\H$ and 
$g_{i,j} \in \G[\alpha_{i,j}]$ with $\alpha_{i,j} := \alpha_i' \cap
\alpha_j'$ are such that 
(cf.\ Figure~\ref{transpatternfig})
\bre
\item[--]
$v_{1,2}g_{1,2}$ is identified with 
$v_1g_1$ in the $\CG[\alpha_1']$-copy at $v_1$
and with $v_2g_2$ in the $\CG[\alpha_2']$-copy at $v_2$;
\item[--]
$v_{2,3}g_{2,3}$ is identified with 
$v_2g_2$ in the $\CG[\alpha_2']$-copy at $v_2$
and with $v_3g_3$ in the $\CG[\alpha_3]$-copy at $v_3$.
\ere

We look at the identifications that occur in the central
$\CG[\alpha_2']$-copy. That the $\alpha_2'$-component
of $\H$ containing $v_{1,2}, v_2$ and $v_{2,3}$ satisfies
condition~(iii) means that at least one of the following must hold
\bne
\item
not $\alpha_{1,2} \strictsubset \alpha_2'$, which implies 
$\alpha_2' \subset \alpha_1'$, or 
\item
not $\alpha_{2,3} \strictsubset \alpha_2'$, which implies $\alpha_2' \subset \alpha_3'$, 
or
\item
the $\alpha_{1,2}$-component of $v_{1,2}$ and 
the $\alpha_{2,3}$-component of $v_{2,3}$ intersect
\\
in the $\alpha_2'$-component of $v_2$ in $\H$.
\ene

In either one of the first two cases the identification of $v_1g_1$
with $v_3g_3$ is obviously mediated by a direct
identification through $\sim$ as defined in $(\dagger)$.
The same is true in the third case, since it
implies that also the $\alpha_1'$-component of $v_1$ and the 
$\alpha_3'$-component of $v_3$ in $\H$ intersect in some vertex $v_2'$
for which there is a $g_2' \in \G[\alpha_{1,2} \cap \alpha_{2,3}]$ that
mediates a direct identification (note that $\H[\alpha_1';v_1]
\supseteq \H[\alpha_{1,2};v_{1,2}]$ and
$\H[\alpha_3';v_3] \supseteq \H[\alpha_{2,3};v_{2,3}]$). 
So the crucial assumption~(iii) implies that
$\sim$ is itself transitive, hence equal to $\approx$. 

Moreover, the isomorphic embeddings of $\alpha'$-components of
$\H[\alpha]$ into $\CG[\alpha']$ for $\alpha' \strictsubset \alpha$,
which are required in~(iii) of Definition~\ref{smallcosetamalgCEdef},
rule out that $\sim$ induces any identifications within the skeleton
$\H[\alpha]$: $v_11 \sim v_21$ can only be mediated by some 
$v_0g_0$ such that $g_0 = g_{10}^{-1} = g_{20}^{-1} \in \G[\alpha_0']$
for $\alpha_0' = \alpha_1' \cap \alpha_2'$. 
By condition~(iii), once more,  the $v_i$ for $i=0,1,2$ are in the same 
$\alpha_0'$-component of $\H[\alpha]$,
which embeds isomorphically into $\CG[\alpha_0']$; it follows that $v_1=v_2$.
So $\H[\alpha] \subsetw \CE(\H,\G,\alpha)$. We summarise
these findings in the following remark.

\bR
\label{smallfootrem}
$\K := \CE(\H[\alpha],\G,\alpha)$ as in 
Definition~\ref{smallcosetamalgCEdef} embeds 
the $\I[\alpha]$-skeleton $\H[\alpha]$ as a weak substructure,
$\H[\alpha] \subsetw \K$, 
and the equivalence relation $\approx$ identifies elements
precisely if they are related as in $(\dagger)$.
\eR

We now look at  strong \emph{freeness} criteria for the manner in which
$\I[\alpha]$-skeletons and small coset amalgams
embed into $\CG[\alpha]$. The basic freeness condition, which 
we just repeat in context as~(i) here, is as in Definition~\ref{Gfreeoverskeletdef}. 
It also occurred as an intrinsic condition for the constituents $\CG[\alpha']$ for
$\alpha'\strictsubset \alpha$ in~(iii) of Definition~\ref{smallcosetamalgCEdef} above.

\begin{figure}
  \[
\xymatrix{ 
\circ \ar@{-}@(ur,ul)@{--}[rr]|*+{\beta}& 
& \circ & 
& \circ && \circ & 
  \\
& \quad && \raisebox{5pt}{$\CE(\H[\alpha],\G,\alpha)$}  &  & \bullet
\ar@{-}[ul]|*{\beta_1} \ar@{-}[ur]|*{\beta_2}& &
\\
v_1 \ar@{-}[rr] \ar@{->}[uu]|*+{g_1}
&&
v_2 \ar@{->}[uu]|*+{g_2}
&
\H[\alpha]
&
v_1 \ar@{-}[r] \ar@{->}[uu]|*+{g_1}
\ar@{}[ur]^*{\alpha_1}
& 
v_0 \ar@{-}[r] \ar@{->}[u]|*+{g_0}
&
v_2 \ar@{->}[uu]|*+{g_2}
\ar@{}[ul]_*{\alpha_2} &
}
\]
\caption{Bridge-freeness.}
\label{bridgefreefig}
\end{figure}

\bD[freeness and bridge-freeness]
\label{bridgefreedef}
\mbox{}\\
An embedded
$\I[\alpha]$-skeleton $\H[\alpha] \subsetw \CG[\alpha]\subset \CG$ is
\bre
\item
  \emph{free} if, for all $\alpha_1,\alpha_2 \strictsubset \alpha$ and
  $v_1,v_2$ from $\H[\alpha]$:
 \[
   \alpha_1(v_1)\, \cap\, \alpha_2(v_2) = \emptyset
   \mbox{ in $\H[\alpha]$ } \; \Longrightarrow \; 
   \alpha_1(v_1) \cap \alpha_2(v_2) = \emptyset
   \mbox{  in $\CG[\alpha]$.} 
 \]
\item 
\emph{bridge-free} if it is free and 
$\K := \CE(\H[\alpha],\G,\alpha)\subsetw \CG[\alpha]$
is well-defined as a weak substructure of $\CG[\alpha]$
over the embedded skeleton $\H[\alpha] \subsetw \CG[\alpha]$,
and if 
for all $\beta \strictsubset \alpha$ and
$v_1g_1, v_2g_2 \in \K$:
if $v_1g_1$ and $v_2g_2$ are $\beta$-connected in $\G[\alpha]$
(i.e.\ if $v_1g_1\G[\beta]  = v_2g_2\G[\beta]$), 
then they are $\beta$-connected inside $\K$
in the following simple form  (cf.\ Figure~\ref{bridgefreefig}): 
\[
  (\ast)
  \qquad
  \mbox{\begin{minipage}{.8\textwidth}
for some $\alpha_1,\alpha_2 \strictsubset \alpha$ such that $g_i \in
\G[\alpha_i]$ for $i=1,2$ there is some $v_0 g_0 \in \K$, 
where $g_0 \in \G[\alpha_0]$ for $\alpha_0 := \alpha_1 \cap \alpha_2$
and such that $v_i g_i \in v_0 g_0
\G[\beta_i]$ for $i =1,2$ where
$\beta_i := \beta \cap \alpha_i$.
\end{minipage}}
\] 
\ere
\eD

In essence, the connectivity condition in~(ii) implies that 
$\beta$-connectivity inside the embedded 
$\K = \CE(\H[\alpha],\G,\alpha) \subsetw \CG[\alpha]$ coincides with 
$\beta$-connectivity in the surrounding $\G[\alpha] \subset \G$.

Provided that $\G[\alpha]$ and 
$\H[\alpha] \subsetw\CG[\alpha]$ are $2$-acyclic,
condition~(i) implies that $\K$ 
embeds isomorphically as a weak substructure 
\[
\CE(\H[\alpha],\G,\alpha)\subsetw \CG[\alpha]
\]
over the embedded $\H[\alpha]$. 
In this case, condition~(ii) then additionally implies, by application to singleton 
 $\beta = \{ e\} \subset \alpha$,
 an induced substructure relationship 
\[
\CE(\H[\alpha],\G,\alpha)\subset \CG[\alpha].
\]

We also note that, for $2$-acyclic $\G[\alpha]$,
the freeness of $\H[\alpha'] \subsetw \CG[\alpha']$ for
$\alpha'\strictsubset \alpha$ implies that $\H[\alpha]$ is
$2$-acyclic (condition~(i) of Definition~\ref{smallcosetamalgCEdef}).
For this consider vertices $v_1,v_2$ of $\H[\alpha]$ that are
$\alpha_i$-connected in $\H[\alpha]$ for $\alpha_1,\alpha_2
\strictsubset \alpha$. If $v_1$ and $v_2$ were not part of the same
$\alpha_0$-component of $\H[\alpha]$, for $\alpha_0 := \alpha_1 \cap
\alpha_2$, then by freeness of
$\H[\alpha_1] \subsetw \CG[\alpha_1]$, $v_1\CG[\alpha_0] \cap v_2
\CG[\alpha_0] = \emptyset$; and this would contradict $2$-acyclicity
of $\CG[\alpha]$ which implies that $v_1$ and $v_2$ are
$\alpha_0$-connected in $\G[\alpha]$.  We state these observations
for further reference.

\bO
\label{CEexistobs}
If $\H[\alpha'] \subsetw \CG[\alpha']$ is 
free for all $\alpha'\strictsubset \alpha$ and $\G[\alpha]$ is $2$-acyclic,
then $\K := \CE(\H[\alpha],\G,\alpha)$ is
well-defined.

If $\H[\alpha] \subsetw \CG[\alpha]$ itself
is free, then $\K$ embeds as a weak substructure
$\CE(\H[\alpha],\G,\alpha) \subsetw \CG[\alpha]$; if
$\H[\alpha] \subsetw \CG[\alpha]$ is even bridge-free, 
then $\K$ 
is embedded as an induced substructure.
\eO

We say that a set $\alpha$ of generators \emph{supports} 
some vertex over some embedded skeleton
if that vertex is linked to the skeleton
by an $\alpha$-walk. We extend this notion to sets of vertices, 
and consider minimal supports as in the following definition.
Compare Remark~\ref{elementcosetsupportrem} for essentially the same
idea in the context of plain Cayley
graphs.

\bD[support and unique minimal support]
\label{supportdef}
\mbox{}\\
Let the small coset amalgam $\K := \CE(\H[\alpha],\G,\alpha)$ be defined according to Definition~\ref{smallcosetamalgCEdef}.
For a vertex $w$ of $\K$ 
or for  a subset $B$ of $\K$,  
or a subset $B \subset \G[\alpha]$
in case that $\K \subsetw\CG[\alpha]$,
we define the notion of \emph{support} as follows:
\bre
\item 
a $\CG[\alpha']$-constituent $v \CG[\alpha']$ of $\K$ 
for $\alpha' \strictsubset \alpha$ \emph{supports} a vertex $w$ 
in $\K$ if $w \in v\CG[\alpha']$; in this case 
we also say that the generator set $\alpha'$ provides support for $w$; 
similarly for a set $B$ of vertices, the requirement is that some
element of $B$ be supported in this manner.
\item
the set of generators $\alpha' \strictsubset \alpha$ provides 
\emph{minimal support} for a vertex $w$ 
in $\K$ if some $\CG[\alpha']$-constituent of
$\K$ supports~$w$, and if $\alpha'$ is a subset of every such
generator set; 
similarly, in case of a set $B$ of vertices, w.r.t.\
containment in all generator sets that support some
element of $B$. 
\item
$w$ has \emph{unique minimal support} if the
supporting constituent $v\CG[\alpha']$ in
$\K$ with such minimal $\alpha'$ is unique (up to the obvious
ambiguity in the choice of $v$ within
its $\alpha'$-component in $\H[\alpha]$); in the case of a set $B$,
the unique minimal supporting constituent $v\CG[\alpha']$ is 
required to be contained in the minimal supporting
constituent of any element of
$B \cap \K$.
\ere
\eD

Vertices of $\CE(\H[\alpha],\G,\alpha)$ have unique
minimal supports by construction, viz.\ the intersection over all 
supporting constituents $v\CG[\alpha']$, which is a supporting
constituent (over the intersection of the
corresponding generator sets). 
Of particular importance to us 
is the case of $\beta$-components of $\CE(\H[\alpha],\G,\alpha)$
in the r\^ole of $B$ in~(ii). Unique minimal support for that case
is part of the \emph{cluster property}  in
Definition~\ref{clusterpropdef} below.

For trivial $\I$ without edges, $\I$-skeletons consist of isolated vertices, and
a small coset amalgam becomes a disjoint union of free amalgamation clusters.
In that case, $\beta$-components $B$ of free
amalgamation clusters (when embedded into a $3$-acyclic $\E$-group)
have unique minimal support: 
$\alpha_0 := \bigcap \{ \alpha \subset E \colon \G[\alpha] \cap B
\not= \emptyset \}$  is a minimal set of supporting generators, and
$\CG[\alpha_0] \subset \CG$ a unique minimal support (cf.\
Remark~\ref{elementcosetsupportrem}). 
We here try to lift this and related insights from free amalgamation clusters
to small coset amalgams over non-trivial skeletons.

\bD[cluster property]
\label{clusterpropdef}
\mbox{}\\
Let the small coset amalgam $\K := \CE(\H[\alpha],\G,\alpha)$ be defined according to Definition~\ref{smallcosetamalgCEdef}.
This small coset amalgam 
has the
\emph{cluster property} if, 
for  $\beta \strictsubset \alpha$, any $\beta$-component
of $\K$ that does not directly intersect the 
embedded skeleton $\H[\alpha]$ 
has unique minimal support, say $v \CG[\alpha_0]$,
and is isomorphic (as a $\beta$-graph) to a free
amalgamation cluster
of $\beta'$-cosets for $\beta' \subset \beta$ all containing
a common $(\beta\cap \alpha_0)$-coset of $v \CG[\alpha_0]$.
\eD

The exception of the simpler case that the $\beta$-component in
question intersects the skeleton, and hence coincides with a
$\beta$-constituent of $\CE(\H[\alpha],\G,\alpha)$, is
necessary. While this $\beta$-constituent is uniquely determined,  
it cannot be \emph{minimal} in the sense of Definition~\ref{supportdef}
since typical elements have smaller supports, some even empty support. 
In the non-degenerate case of a $\beta$-component that is disjoint
from the skeleton, the cluster property states that the unique minimal
supporting consituent (with a necessarily non-empty generator set
$\alpha_0$) serves as a unique bottleneck between the skeleton and
that $\beta$-component. Moreover that $\beta$-component branches out,
as a cluster, into $\beta_i$-cosets of constituents
$v\CG[\alpha_i] \supset v\CG[\alpha_0]$, for all
$\alpha_0 \subset \alpha_i \strictsubset \alpha$, from a single
$\beta_0$-coset of the minimal supporting constituent
$v\CG[\alpha_0]$ for $\beta_i := \beta \cap \alpha_i$.

There is an apparent similarity between the conditions expressed as
bridge-freeness in Definition~\ref{bridgefreedef}~(ii) and the cluster
property as just defined. One may think of bridge-freeness as a
condition that guarantees the analogue of the cluster property 
not just for $\beta$-components inside small coset amalgams
$\CE(\H[\alpha],\G,\alpha)$, but w.r.t.\ the intersections 
between $\CE(\H[\alpha],\G,\alpha)$ and full $\beta$-cosets 
in a surrounding $\CG[\alpha]$.

\bR
\label{clustertwoacycrem}
Let the small coset amalgam 
$\K := \CE(\H[\alpha],\G,\alpha)$ be defined according to
Definition~\ref{smallcosetamalgCEdef}, and  
let $\G[\alpha]$ be $3$-acyclic.
If $\K$ has the cluster property, 
then for 
$\beta,\gamma \strictsubset \alpha$, the $\beta$- and $\gamma$-components of
elements $w$ in $\K$
intersect 
as $\beta(w) \cap \gamma(w) = \delta(w)$
 for $\delta := \beta\cap \gamma$; i.e.\ 
$\CE(\H[\alpha],\G,\alpha)$ is $2$-acyclic.
\eR

\prf
Let $B$ and $C$ be the vertex sets of the $\beta$- and
$\gamma$-components of $w$ in $\K := \CE(\H[\alpha],\G,\alpha)$. 
The claim is trivial if both $\beta(w)$ and $\gamma(w)$ intersect
$\H[\alpha]$. If just one of them, say $B$ intersects
$\H[\alpha]$, the other one has non-trivial minimal supporting 
constituent $v_0\CG[\alpha_0]$ in $\K$, 
and is given as a free amalgamation cluster of $\G[\gamma']$-cosets for 
$\gamma' \subset \gamma$ centered on some $v_0g_0 \in  v_0\CG[\alpha_0]$.
The cluster property implies that the minimal supporting constituent of
the vertex $w \in C$ is of the form $v_0\CG[\alpha_0']$ for some $\alpha_0'$ with 
$\alpha_0 \subset \alpha_0' \subset \beta$ (cf.\ Definition~\ref{supportdef}~(iii)).
As $w \in v_0\CG[\alpha_0'] \subset v_0\CG[\beta]$, 
$B = v_0\CG[\beta]$ follows.
By $2$-acyclicity of any one of the constituents of $\K$ that contain 
$v_0 g_0 \G[\gamma']$ for any one of the $\gamma'$-cosets in the cluster $C$,
the intersection of $B$ with that $\gamma'$-coset is just 
$v_0g_0 \G[\beta \cap \gamma']$. 

In the remaining case, both  $B$ and $C$ are disjoint form
$\H[\alpha]$, hence correspond to free amalgamation 
clusters of $\G[\beta']$- and $\G[\gamma']$-cosets, respectively, 
centered on some $v_1g_1 \in  v_1\CG[\alpha_1]$ and $v_2g_2 \in
v_2\CG[\alpha_2]$ for minimal supporting constituents 
$v_1\CG[\alpha_1]$ for $B$ and $v_2\CG[\alpha_2]$ for $C$.
As $w$ is contained in both, the cluster property now implies that its
minimal supporting constituent is of the form $v_0\CG[\alpha_0]$ such
that $v_i\CG[\alpha_i] \subset v_0\CG[\alpha_0]$ for $i=1,2$. 
In particular $\alpha_1,\alpha_2 \subset \alpha_0$. Suppose now that these
clusters meet not just in $w$ but also in some $w'$ so that $B = \beta(w)=
\beta(w')$ and $C = \gamma(w) = \gamma(w')$; for the claim of the
remark, it suffices to show that $w$ and $w'$ are $\delta$-connected
in $\K$. Just as for $w$, we find that the minimal supporting
constituent for $w'$ must be of the form  $v_0'\CG[\alpha_0']$
such that $\alpha_1,\alpha_2 \subset \alpha_0'$.
It follows that $w$ and $w'$ are elements of the constituents
$v_0\CG[\alpha_0]$ and $v_0'\CG[\alpha_0']$.
These constituents
$v_0\CG[\alpha_0]$ and $v_0'\CG[\alpha_0']$ both contain the elements
$v_1g_1 \in  v_1\CG[\alpha_1]$ and $v_2g_2 \in v_2\CG[\alpha_2]$.
W.r.t.\ to $\alpha_0\cup \alpha_0'$, 
these two constituents form a free amalgam of copies
of $\CG[\alpha_0]$ and $\CG[\alpha_0']$, and $w$ and $w'$ are
both $\beta$- and $\gamma$-connected inside this amalgam:
$\beta$-connected via $v_1g_1 \in  v_1\CG[\alpha_1]$, and  
$\gamma$-connected via $v_2g_2 \in  v_2\CG[\alpha_2]$.
As $\G[\alpha]$ is $3$-acyclic, the free amalgam in question embeds 
isomorphically as an induced substructure into $\CG[\alpha]$. 
By $2$-acyclicity of $\G[\alpha]$, $w$ and $w'$ are 
$\delta$-connected inside this free amalgam, hence also in $\K$.
\eprf

The following can serve to lift the cluster property for 
small coset amalgams at level $\Gamma_\alpha$, i.e.\ for 
$\alpha' \strictsubset \alpha$, to level $\alpha$.
It does, however, crucially require bridge-freeness at level
$\Gamma_\alpha$. That will necessitate a second lifting process,
geared towards the freeness criteria that concern the embedding of $\H[\alpha]$ in
$\G[\alpha]$, rather than intrinsic properties of
$\CE(\H[\alpha],\G,\alpha)$.
Recall form Observation~\ref{CEexistobs} that under the assumptions of
the lemma, which concern level $\Gamma_\alpha$, the small coset
amalgam $\CE(\H[\alpha],\G,\alpha)$ at level $\alpha$ is well-defined.

\bL
\label{firstclusterlem}
Let $\G[\alpha]$ be $3$-acyclic and such that,
for all $\alpha' \strictsubset \alpha$,
the embedded $\I[\alpha']$-skeletons 
$\H[\alpha']  \subsetw \CG[\alpha']$ 
are bridge-free and such that the embedded small coset amalgams
$\CE(\H[\alpha'],\G,\alpha')$ have the cluster property
for all $\alpha' \strictsubset \alpha$.
Then the small coset amalgam $\CE(\H[\alpha],\G,\alpha)$
has the cluster property. 
\eL

\prf
Let $V$ be the vertex set of 
$\K := \CE(\H[\alpha],\G,\alpha)$, $B \subset V$ the vertex set 
of the $\beta$-component of some $w \in V$ in $\K$. 
As $w \in \K$, some $\alpha' \strictsubset \alpha$ supports $w$,
say $w = vg \in v \CG[\alpha']$. For the case of interest we assume
that $B$ does not intersect the skeleton $\H[\alpha]$ of $\K$. 

Suppose first that, for some constituent $v\CG[\alpha']$
with $|\alpha'| = |\alpha|-1$ (maximal in $\CE(\H[\alpha],\G,\alpha)$)
that supports $B$, the set 
$B' := B \cap v\CG[\alpha'] \not= \emptyset$ is not supported by any
strict subset of the generator set $\alpha'$.
We show that in this case $w=vg \in v\CG[\alpha']$ implies 
$B = B' = vg \CG[\beta']$ for $\beta' := \beta \cap \alpha'$. 
Otherwise some $\beta$-edge would have to link some vertex
in $B'$ to some vertex outside $v\CG[\alpha']$; w.l.o.g.\ let these
vertices be $w = vg$ and $v' \in B \setminus B'$, linked by a $b$-edge,
$b \in \beta \setminus \alpha'$, stemming from some
$\CG[\alpha'']$-constituent; but then $b \in \alpha'' \setminus \alpha'$ 
implies that $w \in B'$ would be supported by $\alpha' \cap \alpha''
\strictsubset \alpha'$, contradicting the assumption. 
So $B = B' = vg\G[\beta']\subset v \CG[\alpha']$. Uniqueness of
$v \CG[\alpha']$ as the minimal supporting constituent follows.

For the remaining cases, where $B$ contains elements
supported by smaller generator sets, the above argument shows that
all its non-empty intersections with constituents $v_i\CG[\alpha_i]$ 
of $\K$ for $|\alpha_i| = |\alpha|-1$ will be 
$\beta_i$-cosets for $\beta_i := \beta \cap \alpha_i$ 
supported by generator sets $\alpha_{0i} \strictsubset
\alpha_i$. Bridge-freeness for the embedded small coset amalgam
$\CE(\H[\alpha_i],\G,\alpha_i)$ within this
$\CG[\alpha_i]$-constituent
implies that the $\beta_i$-coset in question has a unique minimal
supporting sub-constituent $v_i\CG[\alpha_{0i}]$ for 
$\alpha_{0i} \strictsubset \alpha_i$ 
(after shifting $v_i$ within its $\alpha_i$-component of $\H[\alpha]$
if necessary, to avoid unnecessary distinctions). 
The $\beta$-component $B$ therefore is the union of 
$\beta_i$-cosets of the form 
\[
v_i g_i \G[\beta_i] \subset v_i \CG[\alpha_i]
\; \mbox{ for } \; g_i \in \G[\alpha_{0i}], \alpha_{0i} \strictsubset
\alpha_i \strictsubset \alpha, |\alpha_i| = |\alpha|-1.
\] 

We need to show that these parts overlap in the form of a 
free amalgamation cluster centered on some $\beta'$-coset 
of a unique minimal supporting constituent for $B$ in $\K$. 
Consider two such $\beta_i$-cosets $v_i g_i \G[\beta_i]\subset v_i
\CG[\alpha_i]$ with $g_i \in \G[\alpha_{0i}]$,
say for $i=1,2$, that overlap in $\K$. 
A common vertex $w \in  v_1 g_1 \G[\beta_1]
\cap v_2 g_2 \G[\beta_2]$ is of the form $v_0g_0$ where 
$g_0 \in \G[\alpha_0]$ for a minimal supporting set $\alpha_0$ of
generators for $w$, so that in particular $\alpha_0 \subset \alpha_1
\cap \alpha_2$. (The situation is as 
represented in Figure~\ref{bridgefreefig} on the right.)

Looking at $v_ig_i$ and $v_0g_0$, which are $\beta_i$-connected in 
$v_i\CG[\alpha_i]$, we can invoke bridge-freeness 
in this $\CG[\alpha_i]$-constituent to find 
$v_{0i} g_{0i} \in v_i g_i \G[\beta_i]$ where $g_{0i} \in
\G[\alpha_{0i}']$ for $\alpha_{0i}' \subset \alpha_{0i} \cap \alpha_0$,
by minimality of those two generator sets.%
\footnote{We invoke bridge-freeness rather than the cluster property
here, since it is not clear that $v_ig_i$ and $v_0g_0$
are $\beta_i$-connected in $\CE(\H[\alpha_i],\G,\alpha_i)$, even though
they are $\beta_i$-connected in $v_i \CG[\alpha_i]$.}
So the cluster property 
for this $\CG[\alpha_i]$-constituent
implies that $v_i\CG[\alpha_{0i}] \subset \alpha_0$. 
The cluster property for the
$\CG[\alpha_0]$-constituent implies that $v_1 \CG[\alpha_{01}] = v_2
\CG[\alpha_{02}] \subset v_0\CG[\alpha_0]$ and that 
$v_1 g_1\G[\beta_1]$ and $v_2 g_2\G[\beta_2]$ overlap in $v_0 \CG[\alpha_0]$.
So any two directly overlapping $\beta_i$-parts of $B$ share the same unique minimal
supporting constituent in their $\CG[\alpha_i]$-constituents (for maximal $\alpha_i
\strictsubset \alpha$) and overlap inside that minimal constituent. 
This shared minimal constituent $v_0\CG[\alpha_0]$ is therefore minimal and unique 
also in $\K$. The nature of $B$ as a union of
$\beta_i$-cosets that pairwise overlap in this minimal constituent
implies that it forms a free amalgamation cluster. 
\eprf

It follows that, under suitable conditions, the classes of free
amalgamation clusters, small coset amalgams, and of free amalgams
between them, satisfy useful closure conditions
w.r.t.\ passage to $\beta$-components for smaller $\beta$.
For the second part of the following corollary we use 
$2$-acyclicity of $\CE(\H[\alpha],\G,\alpha)$ from 
Remark~\ref{clustertwoacycrem} and appeal to Observation~\ref{2acyccompobs}. 
We also note that any $\beta$-component that intersects the skeleton
of any small coset amalgam with constituents from $\Gamma_k(\G)$
and for $\beta\in\Gamma_k$ is necessarily isomorphic to a constituent
$\CG[\beta]$-copy, and cannot possibly extend beyond that in any 
$\E$-graph extensions.

\bC
\label{closurecor}
If $\G[\alpha]$ is $3$-acyclic and for all $\alpha' \strictsubset
\alpha$ the $\I[\alpha']$-skeletons 
$\H[\alpha'] \subsetw \CG[\alpha']$ are bridge-free and 
their embedded small coset amalgams
have the cluster property,
then the small coset amalgam
$\CE(\H[\alpha],\G,\alpha)$ is well-defined, has the cluster property
and is $2$-acyclic. Moreover, all the free amalgams 
\[
\CE(\H[\alpha],\G,\alpha),v \oplus \CG[\alpha']
\]
are well-defined for all vertices $v$ of 
$\CE(\H[\alpha],\G,\alpha)$ and all $\alpha' \strictsubset \alpha$. 
Their $\beta$-components, for
$\beta \strictsubset \alpha$,
are isomorphic to $\CG[\beta]$ itself or to
 a free amalgam between a
  free amalgamation cluster over $\Gamma(\beta)$ and some $\CG[\beta']$ for
$\beta' \strictsubset \beta$ (including, for 
$\beta'=\emptyset$, plain free amalgamation clusters over $\Gamma(\beta)$).
\eC

We note that this corollary establishes downward compatibility  
in the sense that the passage to small coset amalgams and their free
amalgams with individual small $\CG[\alpha']$ has the desired downward
closure property:
the functor that associates free amalgamation 
clusters and their free amalgams with individual $\CG[\alpha']$
is conservative in the sense of Definition~\ref{conservedef}
(as
established in Corollaries~\ref{clusterconservecor}
and~\ref{clusterplusextraconservecor}). By the above corollary
the isomorphism types of the $\beta$-components of 
extended small coset amalgams 
are determined by $\Gamma_\alpha(\CG)$, across all skeletons 
$\H[\alpha]$ that support the construction of these small coset
amalgams according to Definition~\ref{smallcosetamalgCEdef}.
The underlying conservative functor $\F$ that produces all these
isomorphism types in the context of Corollaries~\ref{clusterconservecor}
and~\ref{clusterplusextraconservecor} therefore also safeguards the
inclusion of extended small coset amalgams as components 
in $\E$-graphs $\K$ for inductive unfolding steps $\G \preceq_k \Ghat
= \sym(\K)$ when dealing with $|\alpha| = k$, i.e.\ $\alpha \in \Gamma_{k+1}$.

\medskip
It requires one more argument to lift the (intrinsic) cluster
property for $\CE(\H(\alpha],\G,\alpha)$ to bridge-freeness 
(in relation to some surrounding $\E$-group $\Ghat$ that unfolds 
$\G$ and hence w.r.t.\ embedded $\I[\alpha]$-skeletons) while preserving 
all structure at the level of $\Gamma_\alpha$. This is the purpose
of the following lemma.

\bL
\label{bridgelem}
Suppose that $\G[\alpha]$ is $3$-acyclic and such that 
all embedded $\I[\alpha']$-skeletons $\H[\alpha'] \subsetw
\CG[\alpha']$ for $\alpha' \strictsubset \alpha$ are bridge-free
and that the embedded small coset amalgams
$\CE(\H[\alpha'],\G,\alpha')$ have the cluster property. 
Let $\Ghat \succeq_\alpha \G$
be $3$-acyclic and compatible with all $\alpha$-graphs of the form
\[
 \CE(\H[\alpha],\G,\alpha),v \oplus \CG[\alpha']
\; \mbox{ for }\alpha' \strictsubset \alpha.
\]
Then any embedded $\I[\alpha]$-skeleton
$\Hhat[\alpha] \subsetw \CGhat[\alpha]$ 
is bridge-free.
\eL

Note that $\CE(\Hhat[\alpha],\Ghat,\alpha) \simeq
\CE(\Hhat[\alpha],\G,\alpha)$ as $\Ghat \succeq_\alpha \G$, but
$\Hhat[\alpha]$ may be  a non-trivial unfolding of $\H[\alpha]$.
By the previous lemma, $\CE(\H[\alpha],\G,\alpha)$ has the cluster
property, and the amalgams $\CE(\H[\alpha],\G,\alpha),v \oplus
\CG[\alpha']$ above are well-defined. The conclusion of the 
lemma holds for any $\Ghat' \succeq_k \Ghat$ if $\alpha \in \Gamma_k$.

By Observation~\ref{CEexistobs} the conclusion of the 
lemma implies that the small coset amalgam 
$\CE(\Hhat,\Ghat,\alpha)$ embeds as an induced substructure into
$\CGhat$ and, by Lemma~\ref{firstclusterlem}, it also 
has the cluster property.

\prf[Proof of Lemma~\ref{bridgelem}]
We first observe that
$\H[\alpha'] \simeq \Hhat[\alpha']$ for $\alpha' \strictsubset \alpha$,
due to $\Ghat \succeq_\alpha \G$. Together with
$2$-acyclicity of $\Ghat$, this further implies
that the embedded $\I[\alpha]$-skeleton
$\Hhat[\alpha]$ is $2$-acyclic: any $\alpha_i$-components
$\alpha_i(v)$ for $\alpha_1,\alpha_2 \strictsubset \alpha$ in $\Hhat[\alpha]$ 
that also meet in a vertex $v'$ of $\Hhat[\alpha]$ outside
$\alpha_0(v)$, 
for $\alpha_0 := \alpha_1 \cap \alpha_2$,
would establish a violation of freeness for $\H[\alpha_1] \simeq \Hhat[\alpha_1]$ 
w.r.t.\ $\alpha_0 \strictsubset \alpha_1$ in $\CG[\alpha_1]$.

\begin{figure}
  \[
\nt\hspace{-1.4cm}
\xymatrix{&&v'&& & && g &&
  \\
  &&&& \ar@{<~>}[r] &&&&&
\\
v_1' \ar@{->}[uurr]|*{g_1}
&& v_0'
\ar@{->}[ll]|*{u_{01}} \ar@{->}[rr]|*{u_{02}} \ar@{->}[uu]|*{g_{0}}
&& v_2' \ar@{->}[uull]|*{g_2}
&
v_1 \ar@{-}[rrrr]
\ar@{->}[uurr]|*{g_1}
&& 
&& v_2
\ar@{->}[uull]|*{g_2}
\\
&&\CE(\H[\alpha],\G,\alpha) && &&&
\CE(\Hhat[\alpha],\Ghat,\alpha) &&
}
\]
\caption{Walks in relevant $\alpha_i$-components of 
$\CE(\H[\alpha],\G,\alpha)$ and $\CE(\Hhat[\alpha],\Ghat,\alpha)$.}
\label{patternfig}
\end{figure}

\medskip\noindent
\emph{Claim~1.}
$\Hhat[\alpha] \subsetw \CGhat[\alpha]$ is free in
$\CGhat$, so that $\CE(\Hhat[\alpha],\Ghat,\alpha)$ embeds
isomorphically as a weak substructure of $\CGhat[\alpha] \subset \CGhat$.

\medskip
This relies on compatibility of $\Ghat$ with
the $\E$-graphs of the form $\CE(\H[\alpha],\G,\alpha)$ (which are included in the
stated collection for $\alpha' =\emptyset$). For the following compare
Figure~\ref{patternfig}.
Suppose that for $\alpha_i \strictsubset \alpha$ and $v_1,v_2$ from
$\Hhat[\alpha]$
there is some $g \in v_1 \G[\alpha_1] \cap v_2 \G[\alpha_2]$ in $\Ghat$.
Let $g_i \in \G[\alpha_i]$ and $u \in \alpha^\ast[\I]$ be such that 
$v_2 =  v_1 [u]_{\Ghat}$ and $g = v_i g_i$ for $i=1,2$ (cf.\
Figure~\ref{patternfig}, right). It follows that 
\[
g_1 = [u]_{\Ghat} g_2 \; \mbox{ in } \Ghat,
\]
which by compatibility transfers to the operation of $\Ghat$ on 
$\K := \CE(\H[\alpha],\G,\alpha)$. So there are vertices
$v_1',v_2'$  in $\H[\alpha] \subseteq \CG[\alpha]$ such that 
$v_2' =  v_1' [u]_{\K}$ with a non-trivial intersection of
the constituents $v_i' \CG[\alpha_i]$ in $v' = v_1' g_1 = v_2' g_2$.%
\footnote{Any veretx that matches
  the sort~$s$ of $u \in \alpha^\ast[\I,s]$ in
  $\CE(\H[\alpha],\G,\alpha)$ can serve as $v_1'$.}
By definition of $\CE(\H[\alpha],\G,\alpha)$, the $\alpha_i$-components of the $v_i'$
overlap in the skeleton  $\H[\alpha]$, say in $v_0'$. Let
$ v_1' = v_0'  [u_{01}]_{\K}$ and $v_2' = v_0' [u_{02}]_{\K}$ for
$u_{0i} \in \alpha_i^\ast[\I]$, and
$v' = v_0' g_0$ for some $g_0 \in \G[\alpha_0]$ where
$\alpha_0 = \alpha_1 \cap \alpha_2$.
Then
\[
  [u_{01}]_{\G} g_1 = g_0 =   [u_{02}]_{\G} g_2
\] 
in $\G$. This transfers to $\Ghat$, since both $\G$ and $\Ghat$ are
$2$-acyclic, and $\CG[\alpha_i] \simeq \CGhat[\alpha_i]$, which
implies that free amalgams between such constituents produce
isomorphic results in $\CG$ and $\CGhat$. 
So also in $\Ghat$
\[
g_1 g_2^{-1} =  [u_{01}]_{\Ghat}^{-1}[u_{02}]_{\Ghat}  
\]
is represented by a composition of  an $\I[\alpha_1]$-walk
with an $\I[\alpha_2]$-walk inside the embedded skeleton
$\Hhat[\alpha]$. In particular, as $g_1 g_2^{-1}$ connects $v_1$ to
$v_2$ in $\CE(\Hhat[\alpha],\Ghat, \alpha)$, 
$v_0 := v_1  [u_{01}]_{\Ghat}^{-1}$ witnesses the non-empty overlap 
between $\alpha_1(v_1)$ and $\alpha_2(v_2)$ in the skeleton
$\Hhat[\alpha]$. Moreover,
$g = v_1  [u_{01}]_{\Ghat}^{-1} g_0$ is $\alpha_0$-connected to this overlap 
by an $\alpha_0$-walk. This implies freeness of $\Hhat[\alpha]$ in
$\CGhat$ and that $\CE(\Hhat[\alpha],\Ghat,\alpha)$ embeds 
into $\CGhat$, so far as a weak substructure. 

\medskip\noindent
\emph{Claim~2.}
$\Hhat[\alpha] \subsetw \CGhat$ is bridge-free.%
\footnote{This also stengthens Claim~1 to $\CE(\Hhat[\alpha],\Ghat,\alpha) \subset \CGhat$, cf.\ Observation~\ref{CEexistobs}.}

\medskip
For Claim~2 we need to show that elements of
$\CE(\Hhat[\alpha],\Ghat,\alpha)\subsetw \CGhat$ 
that are members of the same $\beta$-coset of $\Ghat$  
for some $\beta \strictsubset \alpha$ must be
$\beta$-connected inside $\CE(\Hhat[\alpha],\Ghat,\alpha)$
in the manner expressed in criterion~$(\ast)$ of
Definition~\ref{bridgefreedef}, as illustrated in Figure~\ref{bridgefreefig}. 
In the following we use $\K:=\CE(\H[\alpha],\G,\alpha)$ to denote the
small coset amalgam at the level of $\G$, and 
$\Khat := \CE(\Hhat[\alpha],\Ghat,\alpha)\subsetw \CGhat$ for the
small coset amalgam at the level of $\Ghat$.

Consider elements $v_i g_i \in v_i \CGhat[\alpha_i]$ in $\Khat$ 
for $v_i$ from $\Hhat[\alpha] \subsetw \CGhat[\alpha]$, 
represented as members of constituents $v_i \CGhat[\alpha_i]$ over
minimal supporting sets $\alpha_i \strictsubset \alpha$ (existence and
uniqueness of minimal supports for elements is intrinsically guaranteed
by definition of $\CE(\Hhat[\alpha],\Ghat,\alpha)$; here we could 
meanwhile also argue by freeness of $\Hhat[\alpha]$). 

\begin{figure}
  \[
\nt\hspace{-1cm}
\xymatrix{
 &\circ \ar@{->}[rr]|*{[w]_\K}&& \circ & &&  \circ \ar@{->}[rr]|*{[w]_{\Ghat}} && \circ 
 \\
\K && && \ar@{<~>}[r] &&& \qquad
&& \makebox(20,0)[l]{$\Khat \subsetw \CGhat$}
\\
\H[\alpha]
&
v_1' \ar@{->}[uu]|*+{g_1}
\ar@{->}[rr]|*{[u]_\K}
&& v_2' \ar@{->}[uu]|*+{g_2}
&&&
v_1 \ar@{->}[uu]|*+{g_1}
\ar@{->}[rr]|*{u}
&& v_2 \ar@{->}[uu]|*+{g_2}&\Hhat[\alpha]
}
\]
\caption{The argument for bridge-freeness.}
\label{bridgefreelemfig}
\end{figure}

Assume that $v_2g_2$ is linked to $v_1g_1$ by some
$\beta$-walk, i.e.\  that $v_2g_2 = v_1g_1 h$ for some $h \in
\Ghat[\beta]$, say $h = [w]_{\Ghat}$ for $w \in \beta^\ast$. 
Let $v_2 = v_1 [u]_{\Ghat}$ for some $u \in \alpha^\ast[\I]$ tracing
an $\I[\alpha]$-walk from $v_1$ to $v_2$ in $\Hhat[\alpha]$
(cf.\ Figure~\ref{bridgefreelemfig}, right).
It follows that 
\[
g_1 \,[w]_{\Ghat} = [u]_{\Ghat}\, g_2 \; \mbox{ in } \Ghat.
\]

Consider vertices
$v_i'$ in $\H[\alpha] \subsetw \K = \CE(\H[\alpha],\G,\alpha)$
that are linked by a $u$-labelled walk, 
such that in particular $v_2' = v_1' [u]_{\K}$ (cf.\ the lefthand
side of Figure~\ref{bridgefreelemfig}). By compatibility of $\Ghat$
with $\K$, the operation of 
$[w]_{\K}$ maps $v_1'g_1$ to $v_2' g_2$ in $\K$. Therefore  
$v_1'g_1$ and $v_2' g_2$ are $\beta$-connected in $\K$.
As $\K$ has the cluster property (cf.\ Lemma~\ref{firstclusterlem}),
the $\beta$-component of
$v_1'g_1$ and $v_2' g_2$ in $\K$ either intersects the skeleton
$\H[\alpha] \subsetw \K$ or else has a unique minimal supporting
constituent $v_0' \CG[\alpha_0]$. In the latter case, the cluster property
yields some element $v_0' g_0 \in v_0' \CG[\alpha_0]$, in which the
cluster can be rooted; and as elements of this cluster the $v_i'g_i$
must be of the form $v_i'g_i \in v_0'g_0
\G[\beta_i']$ where $\beta_i := \beta \cap \alpha_i'$ for suitable
$\alpha_i' \strictsubset \alpha$ with $\alpha_i \subset
\alpha_i'$ (generator sets $\alpha_i$ are minimal for $v_ig_i$ in
$\Khat$, and hence also for $v_i'g_i$ in $\K$). In the former case,
i.e.\ if the $\beta$-cluster degenerates to a $\beta$-constituent
of the form $v_0' \CG[\beta]$, we may use
$v_0' g_0 = v_0' 1 \in \CG[\alpha_0]$ for $\alpha_0 = \emptyset$, and
similarly obtain that $v_i'g_i \in v_0'g_0
\G[\beta_i']$ where $\beta_i := \beta \cap \alpha_i'$ for suitable
$\alpha_i' \strictsubset \alpha$ with $\alpha_i \subset \alpha_i'$
(in this case, indeed, minimality of $\alpha_i$ implies that
$\alpha_i \subset \beta$ so that $\beta_i = \alpha_i = \alpha_i'$ can be used). 
In either case, therefore, $v_i' g_i = v_0'g_0 [w_i]_{\G}$ in $\K$,
for suitable $w_i \in \beta_i^\ast$. And as the operation of $[w]_{\K}$ maps 
$v_1'g_1$ to $v_2' g_2$ in $\K$,   $v_0'g_0$ is a fixpoint of 
the operation of 
\[
[w_1]_{\K}[w]_{\K}[w_2^{-1}]_{\K}. 
\]

The same argument applies to $v_0'g_0$ as an element 
of the free amalgam of $\K$ with a copy of $\CG[\beta]$ attached at 
$v_0'g_0$,
\[
\K, v_0'g_0 \oplus \CG[\beta].
\] 

Because $\Ghat$ is assumed to be also compatible with that,
also here the operation of $[w]$ must map 
$v_1'g_1$ to $v_2' g_2$. In $\CG[\beta] \subset \K, v_0'g_0 \oplus
\CG[\beta]$,
however,  $w_1 w w_2^{-1} \in \beta^\ast$ can only have
$v_o' g_0$ as a fixpoint if 
$[w_1 w w_2^{-1}]_\G = 1$. So $[w]_\G = [w_1]_\G^{-1} [w_2]_\G$ and
therefore also $[w]_{\Ghat} = [w_1]_{\Ghat}^{-1}[w_2]_{\Ghat}$ 
in $\Ghat$, since 
$\G \preceq_\alpha \Ghat$ and $\beta \strictsubset \alpha$.
As the isomorphism type of the overlap between
$\alpha_1$- and $\alpha_2$-cosets is the same in $\CG$ and $\CGhat$, 
it follows that also in $\Ghat$ 
there is $v_0 \in \Hhat[\alpha] \subsetw \Khat$ such that 
\[
v_1g_1 [w_1]_{\Ghat}^{-1} = v_0 g_0 = v_2g_2[w_2]_{\Ghat}^{-1}
\]
so that $v_1g_1$ and $v_2g_2$ are $\beta$-related to $v_0g_0$ within 
$v_1\CG[\alpha_1']$ and $v_2\CG[\alpha_2']$, respectively, as required
by bridge-freeness.
\eprf

\section{Construction of N-acyclic E-groups over I}
\label{Iacycsec}

Again, we fix a constraint graph $\I$ and only consider
$\E$-groups $\G$ that are compatible with
$\I$ according to Proviso~\ref{constraintpatternprov}.
We collect the desirable qualities of $\E$-groups towards
freeness and cluster properties at levels defined in terms of the sets
$\Gamma_k$, i.e.\ w.r.t. the numbers of generators from $E$ that are involved.
This is geared towards an inductive construction
w.r.t.\ $k$, 
see Lemma~\ref{goodinductionlem} 
and Figure~\ref{indictionfig}. In the inductive process we may 
assume global $N$-acyclicity (for the target $N \geq 3$) for
the essential intermediate $\E$-groups $\G$,
by interleaving steps based on Proposition~\ref{inductiveNacycprop}.

\bD 
\label{goodatlevelkdef}
Let $N \geq3$. For $k \geq 1$,  
$\G$ is \emph{$N$-good for $\I$ at level $\Gamma_k$} if
\bre
\item
  $\G$ is $N$-acyclic,
\item
  the embedded $\I$-skeletons $\H[\alpha]$ are bridge-free
  in $\CG[\alpha]$, for $\alpha \in \Gamma_{k}$, 
\item
 $\G$ is compatible with 
\bre
\item[--]
  amalgamation chains of lengths up to $N-2$ over
 $\Gamma_k(\CG)$,
 \item[--]
 free amalgamation clusters
 over $\Gamma_k(\CG)$, 
\item[--]
free amalgams of such clusters with 
  an extra $\CG[\alpha'] \in \Gamma_k(\CG)$.
\ere
\ere
 \eD

$N$-goodness at level $\Gamma_1$ is just $N$-acyclicity: 
conditions~(ii) and~(iii) are vacuous for $\alpha = \emptyset$.  
Goodness at level~$\Gamma_{|E|+1}$ is the eventual goal. 
The following lemma allows us to increase the level of $N$-goodness by one.
So, by induction, we obtain an $N$-acyclic $\G$ that is (bridge-)free
over its embedded $\I[\alpha]$-skeletons for all $\alpha \subset E$;
compare Figure~\ref{indictionfig}.

\begin{figure}
  \[
\xymatrix{ 
*++{\makebox(30,10)[l]{\nt\quad$\G$: good at level $\Gamma_k$}} 
\ar@{=>}[ddd]|*+{\rotatebox{-90}{$\preceq_{k}$}}
^*{\makebox(190,0)[l]{\nt\;\;\;\small Corollary~\ref{closurecor}}} 
 \ar@{->}[drr]%
^*{\makebox(170,0){\small Observation~\ref{CEexistobs}}} 
&
\\
& &
*++{\makebox(10,10)[l]{$\CE(\H[\alpha],\G,\alpha)$ for $\alpha \in \Gamma_{k+1}$}} 
 \ar@{->}[d]
^*{\mbox{\small Lemma~\ref{firstclusterlem}}}
\\
& & *++{\makebox(10,10)[l]{$\CE(\H[\alpha],\G,\alpha) \oplus \CG[\alpha']$ for 
$\alpha \in \Gamma_{k+1},\alpha' \in \Gamma_k$}} 
 \ar@{->}[dll]
^*{\makebox(120,0){\small Lemma~\ref{bridgelem}}}
\\
*++{\makebox(30,10)[l]{\nt\quad$\Ghat$: good at level $\Gamma_{k+1}$}} 
& \qquad & \qquad& \qquad & \qquad & \qquad& \qquad & \qquad&
}
\]
\caption{Inductive pattern: from level $\Gamma_k$ to level
  $\Gamma_{k+1}$ for Lemma~\ref{goodinductionlem}.}
\label{indictionfig}
\end{figure}

\bL
\label{goodinductionlem}
Let $N \geq 3$, and let $\G$ be $N$-good for $\I$ at level $\Gamma_{k}$.
Let $\Ghat^\ast$ be obtained as 
$\Ghat^\ast := \sym(\K)$ for $\K$ 
with components $\CG$ and, 
  for all embedded $\I$-skeletons $\H[\alpha] \subsetw
  \CG[\alpha]$ where $|\alpha| = k$, i.e.\ for $\alpha \in
  \Gamma_{k+1} \setminus \Gamma_k$: 
\[
\CE(\H[\alpha],\G,\alpha),v
 \oplus \CG[\alpha']
\; \mbox{ for all $v$ and all $\alpha' \strictsubset \alpha$.}
\]
Then $\Ghat^\ast \succeq_k \G$ and any $\Ghat \succeq_{k+1} \Ghat^\ast$
that is $N$-acyclic and compatible with free amalgamation clusters 
over $\Gamma_{k+1}(\Ghat^\ast)$ 
and free amalgams of such clusters with an extra 
$\CG[\alpha'] \in \Gamma_{k+1}(\CGhat^\ast)$
is $N$-good for $\I$ at level $\Gamma_{k+1}$.
\eL

\prf
Small coset amalgams $\CE(\H[\alpha],\G,\alpha)$
for $\alpha \in \Gamma_{k+1}$ 
as well as their free amalgams $\CE(\H[\alpha],\G,\alpha),v
 \oplus \CG[\alpha']$ with $\CG[\alpha']$ for $\alpha' \strictsubset
 \alpha$ are well-defined by Corollary~\ref{closurecor} (itself
 essentially based on Lemma~\ref{firstclusterlem} and Observation~\ref{CEexistobs}).
These ingredients in $\K$ for $|\alpha|=k$ together with  
$\G$ being good for $\I$ at level $\Gamma_k$, preserve compatibility
of $\G$ with all $\beta$-components of 
$\K$ for  $\beta \in \Gamma_k$, by Corollary~\ref{closurecor}. 
Therefore $\CG$ as an ingredient in $\K$ entails that $\Ghat^\ast =
\sym(\K)$ unfolds $\G$ in a manner that preserves $\Gamma_k(\G)$, i.e.\
such that $\G \preceq_k \Ghat^\ast$. It essentially remains to argue that the embedded
$\alpha$-skeletons $\Hhat^\ast[\alpha] \subsetw \CGhat^\ast[\alpha]$
for $|\alpha| = k$ are  bridge-free: 
condition~(ii) for goodness at level~$\Gamma_{k+1}$
boils down to the case for $|\alpha|=k$ as this condition is preserved
by $\preceq_k$ for smaller $\alpha$. But compatibility of $\Ghat^\ast$
with the $\CE(\H[\alpha],\G,\alpha),v \oplus \CG[\alpha']$, which are part
of $\K$ for these $\alpha$ and $\alpha' \strictsubset \alpha$, guarantees 
bridge-freeness by Lemma~\ref{bridgelem}.

It is clear that $\Ghat \succeq_{k+1} \Ghat^\ast$  is good at level $\Gamma_{k+1}$
whenever $N$-acyclicity for condition~(i) and compatibility
conditions~(iii) for goodness at level $\Gamma_{k+1}$ are met.
This can be achieved in a $\preceq_{k+1}$-preserving unfolding
chain, since already $\G \preceq_k \Ghat^\ast$ satisfies these
conditions at level $\Gamma_k$, so that the new compatibility
requirements can be implemented without affecting subgroups 
at level $\Gamma_{k+1}$. This uses the conservative nature of 
the functors (cf.\ Definition~\ref{conservedef}) that generate  
amalgamation chains, free amalgamation clusters and 
amalgams of free amalgamation clusters with one extra $\CG[\alpha]$:
\bre
\item[--]
Corollary~\ref{chainconservecor} for amalgamation chains, 
\item[--]
Corollaries~\ref{clusterconservecor}
and~\ref{clusterplusextraconservecor} for
plain and extended amalgamation clusters.
\ere
\eprf

Similar to Theorem~\ref{symmacycEgroupthm} we sum up the 
construction and emphasise its preservation of symmetries over $E$.
For $N$-acyclicity over $\I$ compare Definition~\ref{INacycdef} and 
its connection with freeness from Lemma~\ref{cosettocosetoverIlem}.

\bT
\label{symmacycIEgroupthm}
For every finite $\E$-group $\G$ that is compatible with the
constraint graph $\I$ and every $N \geq 2$ there is 
a finite $\E$-group $\Ghat \succeq \G$ that is $N$-acyclic over $\I$
and fully symmetric over $\G$ in the sense that every permutation 
$\rho \in \mathrm{Sym}(E)$ 
of the generator set $E$ that is a symmetry of $\I$ and $\G$ extends to a
symmetry of $\Ghat$: $\G^\rho \simeq \G \Rightarrow \Ghat^\rho \simeq \Ghat$.
\eT

In particular, we may obtain finite $\E$-groups 
$\Ghat$ that are $N$-acyclic over $\I$ and 
fully symmetric over $\I$ in the sense of admitting every symmetry of
$\I$ as a symmetry. For this we may start e.g.\ from 
$\G := \sym(\H)$ where $\H$ is the disjoint union of $\I$ and 
the hypercube $\H = 2^E$.

\section{From groups to groupoids}
\label{groupoidsec}

In terms of the combinatorial action of the generators $e \in E$ on 
an $\E$-graph $\H$, and by extension of the monoid structure of $E^\ast$
on $\H$, the involutive nature of $\pi_e \in \mathrm{Sym}(V)$ 
is closely tied to the undirected nature of $e$-edges in
$\E$-graphs. We want to overcome this restriction by 
allowing for directed $e$-edges. At the same time we may want to relax
the strictly prescribed uniformity between vertices. The latter has  
already been achieved in the context of involutive generators 
with constraint graphs $\I$ in \S~\ref{constraintsec},
if we think of $S$-markings, e.g.\ in $\I$-skeletons according to
Definition~\ref{Ialphaskeletondef}, as marking out distinct sorts of vertices. 
So now we want to allow for vertices of different sorts 
with \emph{directed} transitions via $e$-edges between vertices of 
specific sorts. Some applications of 
related notions of acyclicity in graph and hypergraph 
structures inspired by Cayley graphs 
are very naturally cast in 
terms of multi-sorted multi-graph structures and  
related groupoids, cf.~\cite{Otto12JACM,arXiv2018}.
In this section we directly reduce 
the construction of groupoids with the desired 
coset acyclicity properties to the construction for groups with
involutive generators from the previous sections.

\subsection{Constraint patterns for groupoids}
\label{constraintpatternsec}

In the following we consider groupoid structures with 
a specified pattern of sorts (types of elements, objects) and generators
(for the groupoidal operation, morphisms). Groupoids in our sense can also be  
associated with inverse semigroups or monoids of correspondingly restricted
patterns. 
We choose a format for the specification of their sorts 
that is very similar to the format of $\E$-graphs,
and call this specification a \emph{constraint pattern}.
The corresponding structures generalise the constraint graphs 
of Section~\ref{constraintsec} in the desired direction.
Such a template will be a directed multi-graph with edge set $E$ and vertex
set $S$, but unlike $\E$-graphs considered so far, the edges $e \in E$
are directed, with an explicit operation of edge reversal.

\bD[constraint pattern $\I$]
\label{Idef}
\mbox{}\\
A \emph{constraint pattern} is a multi-graph 
$\I = (S,E, \iota_1,\iota_2, \cdot^{-1})$, which we formalise as a
two-sorted structure with a set $S$ of vertices and a set $E$ of edges 
as sorts, linked by surjective maps 
$\iota_i \colon E \rightarrow S$ 
that associate a source and target vertex with every edge 
$e \in E$, and a fixpoint-free and involutive operation of 
edge reversal $e \mapsto e^{-1}$ on $E$ that is compatible 
with the $\iota_i$ in the sense that $\iota_1(e^{-1}) = \iota_2(e)$. 
\eD

In the following we mostly abbreviate the notation for a
constraint pattern $\I$ as above to 
just $\I = (S,E)$, leaving the remaining structural details implicit.

For $s,s' \in S$, we let $E[s,s'] := \{ e \in E \colon \iota_1(e)= s,
\iota_2(e)=s' \}$ be the set of edges linking source $s$ to target $s'$. 

In order to extend the notion of  $\I$-reachability (based on an undirected
constraint graph $\I$ in Section~\ref{constraintsec}) 
to a similar concept of $\I$-reachability
w.r.t.\ a constraint pattern $\I$ we
consider words that label \emph{directed} walks in $\I$. 
It is constructive to compare 
related notions in Section~\ref{constraintsec}.
Relevant subsets $\alpha \subset E$ will in the following always
be closed under edge reversal: $\alpha = \alpha^{-1}$.
A \emph{reduced} word over $E$ now is a word in which no 
$e \in E$ is directly followed or preceded by its inverse $e^{-1}$.

\bD[$\I$-words and $\I$-walks]
\label{Jalphawalksdef}
\mbox{}\\
Natural sets of (reduced) $\alpha$-words that occur as edge label sequences
along directed walks on the constraint pattern $\I$ 
are defined as follows:
\bre
\item[--]
$\alpha^\ast[\I] \subset E^\ast$ consists of those 
$w \in \alpha^\ast$ that label a walk in $\I$;
\item[--]
$\alpha^\ast[\I,s] \subset E^\ast$
consists of those $w \in \alpha^\ast$ that label a walk 
from $s$ in $\I$;
\item[--]
$\alpha^\ast[\I,s,t] \subset E^\ast$
consists of those $w \in \alpha^\ast$ that label a walk 
from $s$ to $t$ in $\I$.
\ere
\eD

In particular we write
$E^\ast[\I]$ for the set of all (reduced) words over $E$ that label walks in $\I$, 
and naturally extend the $\iota$-maps to
all of $E^\ast[\I]$ as follows. Since a walk from $s$ to $t$ in $\I$ is a 
sequence $s= s_0,e_1,s_1,\ldots,e_n,s_n=t$ such that 
$\iota_1(e_{i}) = s_{i-1}$ and $\iota_2(e_i) = s_{i}$ for $1 \leq i
\leq n$, this walk is fully determined by the sequence of edges and
can be identified with the word $w = e_1\ldots e_n \in E^\ast$. So we think
of $E^\ast[\I]$ as the set of all words $w = e_1\ldots e_n$ with 
$\iota_2(e_i) = \iota_1(e_{i+1})$ for $1\leq i < n$, and put
$\iota_1(w) := \iota_1(e_1)$ and $\iota_2(w) := \iota_2(e_n)$. This word
$w$ labels a walk  in $\I$ from the source vertex $\iota_1(w)$ to the target
vertex $\iota_2(w)$. Correspondingly 
\[
E^\ast[\I,s,t] = \{ w \in E^\ast[\I] \colon \iota_1(w) = s,
\iota_2(w) = t\}.
\]

Concatenation between (reduced) words $w_1$ and $w_2$ is defined 
as $w_1w_2 \in E^\ast[\I,\iota_1(w_1),\iota_2(w_2)]$ whenever their 
$\iota$-values match in the sense that $\iota_2(w_1) = \iota_1(w_2)$.
This concatenation operation reflects composition of walks in $\I$.

\medskip
We think of the vertex set $S$ of $\I$ as a set of \emph{sites} or
\emph{sorts} marked by \emph{vertex colours} and of the edge set $E$ as a set of 
\emph{links} or \emph{edge colours} that will govern the r\^oles of
elements and generators in corresponding groupoids, as in the
following definition. 
A groupoid is viewed as a group-like structure with groupoid elements of 
specified sorts. These sorts are pairs of sites and specify the source
and the target site of the groupoid element. The
groupoidal composition operation, which is partial overall, is fully 
defined for pairs of elements that share the same interface site.

\subsection{I-groupoids and their Cayley graphs}  
\label{Igroupoidcayleysec}

\bD[$\I$-groupoid]
\label{groupoiddef}
\mbox{}\\
An \emph{$\I$-groupoid} 
based on the constraint pattern $\I = (S,E)$ is a groupoid
structure of the form 
$\G = (G,(G_{s,t})_{s,t \in S}, \cdot\,, (1_s)_{s \in S}, (g_e)_{e \in E})$ where 
\bre
\item
the family $(G_{s,t})_{s,t \in S}$ partitions the universe
$G$ of groupoid elements;%
\footnote{Some of the sets $G_{s,t}$ may be empty as $\I$ is not
required to be connected.}
\item
$\cdot\;$ is a groupoidal composition operation mapping any pair of elements 
in $G_{s,t} \times G_{t,u}$ to an element of $G_{s,u}$, for all
combinations of $s,t,u \in S$;
\item
$1_s \in G_{s,s}$ is a left and right neutral element w.r.t.\ $\cdot\,$, 
for every $s \in S$;
\item
$G$ is generated by the family of pairwise distinct elements
$g_e \in G_{\iota_1(e),\iota_2(e)}$ for $e \in E$,
where $g_{e^{-1}}$ is the groupoidal inverse of $g_e$ w.r.t.\ $\cdot\;$:
$g_{e^{-1}} = g_e^{-1}$ in the sense that 
$g_{e} \cdot g_{e^{-1}} = 1_s$ for $s = {\iota_1(e)}$
and $g_{e^{-1}} \cdot g_e= 1_{s'}$ for $s' = {\iota_2(e)}$.
\ere
\eD

The set $E^\ast[\I,s,t]$ of those (reduced) words over $E$ that
label walks from $s$ to $t$ in $\I$ now suggests an interpretation of
$w \in E^\ast[\I,s,t]$ as a product of generators that represents a
groupoid element in $G_{s,t}$. With $w = e_1\cdots e_n \in
E^\ast[\I,s,t]$ we associate the groupoid element  
$[w]_\G$ that is the 
groupoidal composition 
\[
\textstyle
[w]_\G 
:= \prod_{i=1}^n g_{e_i} = g_{e_1} \cdots g_{e_n} \in G_{s,t}.
\]

Note that $[w]_\G \in G_{s,t}$ precisely for $s = \iota_1(w)$ and $t=
\iota_2(w)$. 
The $\iota$-maps extend to the elements of an 
$\I$-groupoid $\G$ according to
$\iota_i(g) = s_i$ for $i =1,2$ 
if, and only if, $g \in G_{s_1,s_2}$
if, and only if, 
$g = [w]_\G$ for some $w \in E^\ast[\I,s_1,s_2]$.

There is an obvious notion of homomorphisms between
$\I$-groupoids, which needs to
respect the r\^ole of the distinguished generators.
The existence of a homomorphism $h \colon \Ghat \rightarrow \G$ (uniquely
determined and surjective if it exists) 
is expressed as $\Ghat \succeq \G$. 
The following is analogous to Observation~\ref{freeinvgroupobs}.

\bO
\label{freeIgroupoidobs}
The quotient of $\I$-walks w.r.t.\ cancellation of direct edge 
reversal, as represented by the set of reduced words in 
$E^\ast[\I]$ 
as label sequences, forms an $\I$-groupoid with concatenation. This  
can be regarded as the \emph{free $\I$-groupoid}, which has any other 
$\I$-groupoid as a homomorphic image.
\eO

Analogous to Definition~\ref{Esymmdef} we also consider 
symmetries that are induced by admissible re-labellings of the 
sites and links. The natural candidates stem from symmetries 
of the constraint pattern $\I$. A symmetry of $\I$ is just an 
automorphisms in the usual sense for $\I$ as a multi-sorted structure,
induced by matching permutations of the sets $E$ and $S$ that are
compatible with the $\iota_i$ and with edge reversal.

\bD[symmetry over $\I$]
\label{Isymmgroupoiddef}
\mbox{}\\
An automorphism $\rho$ of $\I$ is a \emph{symmetry} of an $\I$-groupoid $\G$
if the renaming of sorts and generators according to $\rho$,
$s \mapsto \rho(s)$ and $g_e \mapsto g_{\rho(e)}$, 
yields an isomorphic $\E$-groupoid, $\G^\rho \simeq \G$.  
\eD

In an $\I$-groupoid $\G$, the set $\alpha^\ast[\I,s,t]$ of (reduced) words over a
subset $\alpha = \alpha^{-1} \subset E$ carves out 
a \emph{generated subgroupoid}  
$\G[\alpha] \subset \G$, as well as corresponding 
\emph{groupoidal cosets} 
at $g \in G$. These are defined in the obvious manner as
\[
\barr{@{}rcl@{}}
\G[\alpha] &=& \bigcup_{s,t} \G[\alpha,s,t]
\mbox{ where }
\\
\hnt
\G[\alpha,s,t] &=& \{ [w]_\G \in G \colon w \in \alpha^\ast[\I,s,t] \},
\\
\hnt
\mbox{ and } \;
g \G[\alpha] &=& \bigcup_t \{ g \cdot [w]_\G \colon  w \in \alpha^\ast[\I,\iota_2(g), t] \}.
\earr
\]

As the constraint pattern $\I$ will mostly be fixed, we shall often
suppress its explicit mention and write, e.g., just $E^\ast[s,t]$, or 
$\alpha^\ast[s,t]$, just as we already wrote $\G[\alpha]$ or
$\G[\alpha,s,t]$ when $\I$ was implicitly determined by the
$\I$-groupoid $\G$.%
\footnote{Note that $g\CG[\alpha]$ is always connected while $\G[\alpha]$ and
its Cayley graph (defined below) may or may not be.}

\medskip
The notion of a Cayley graph for a groupoid $\G$
encodes the operation of generators on groupoid elements, by right
multiplication, as with Cayley graphs of groups (cf.\ Definition~\ref{Cayleydef}).

\bD[Cayley graph of an $\I$-groupoid]
\label{groupoidCayleydef}
\mbox{}\\
The \emph{Cayley graph} of an $\I$-groupoid $\G = (G,(G_{s,t})_{s,t\in S},\cdot\,,(1_s)_{s
  \in S}, (g_e)_{e \in E})$ is the directed 
edge-coloured graph $\CG := \mathrm{Cayley}(\G)= (G,(R_e)_{e \in E})$ 
with vertex set $G$ and edge sets of colour $e \in E$ according to  
\[
R_e := \{ (g,g \cdot g_e) \colon g \in G_{s,t} \mbox{ for some } s \in S \mbox{
  and } \iota_1(e) = t \}.
\]
\eD

As with Cayley graphs for $\E$-groups, the Cayley graphs of
$\I$-groupoids are more homogeneous than the underlying groupoid.
In particular the neutral elements $1_s$ are not identified as
individual constants in $\CG$. 
What is still recognisable in $\CG$, for an $\I$-groupoid $\G$, is 
membership in the sets 
\[
G[s,\ast]:= \iota_1^{-1}(s) 
\; \mbox{ and } \;
G[\ast,s]:= \iota_2^{-1}(s), 
\]
which are identified by the existence of corresponding 
incoming or outgoing $R_e$-edges for $e$ 
with $\iota_2(e) = s$ or $\iota_1(e) = s$, respectively. 
The algebraic structure of the $\I$-groupoid $\G$
is still fully determined by its Cayley graph $\CG$
through the corresponding action of partial permutations.
In the terminology of~\cite{Lawson} it can be recovered 
as a groupoid embedded in the full 
\emph{symmetric inverse semigroup} $I(G)$ 
over its vertex set $G$. In analogy with
the case of groups and their Cayley graphs, where the group is 
realised as a subgroup of the full symmetric group of global
permutations of the vertex set, the groupoid is realised
as a subgroupoid of the set of all bijections between 
the relevant sets $G[\ast,s]$.

\bO
\label{groupoidbackfromCayleyobs}
The $\I$-groupoid $\G$ is isomorphic to the 
$\I$-groupoid generated by the following bijections 
$\pi_e$ for $e \in E[s,s']$:
\[
\barr{@{}rcl@{}}
\pi_e \colon G[\ast,s] &\longrightarrow& G[\ast,s']
\\
g &\longmapsto& g \cdot g_e,
\earr
\]
where $g \cdot g_e$ is identified as the unique vertex $g'$ of $\CG$
for which $(g,g') \in R_e$. Here $1_s = \mathrm{id}_{G[\ast,s]}$ is 
the identity on $G[\ast,s]$. 
\eO

Note that $\pi_e$ is a partial bijection of the set $G$ but total 
in restriction to the indicated domain and range.

The analogy is carried further in the following
(cf.\ Definitions~\ref{Egraphdef},~\ref{symdef} and~\ref{compatdef} in
connection with $\E$-graphs and $\E$-groups.)
A major difference is the absence of a simple \emph{completion}
operation for $\I$-graphs.%
\footnote{Indeed the proposal of a naive completion operation
accounts for the major flaw in~\cite{lics2013}.}

\bD[$\I$-graph]
\label{Igraphdef}
\mbox{}\\
An \emph{$\I$-graph}, for a constraint pattern $\I = (S,E)$, 
is a vertex- and edge-coloured directed graph 
$\H = (V,(V_s)_{s \in S}, (R_e)_{e \in E})$, whose vertex set $V$ is partitioned into 
non-empty subsets $V_s$ of vertices of colour $s \in S$, 
with edge sets $R_e \subset V_{\iota_1(e)}\times V_{\iota_2(e)}$ 
of colour $e$ for $e \in E$ such that $R_{e^{-1}} = R_e^{-1}$. 
The $\I$-graph $\H = (V,(V_s)_{s \in S}, (R_e)_{e \in E})$ is
\emph{complete} if each $R_e$ is a complete matching between
$V_{\iota_1(e)}$ and $V_{\iota_2(e)}$ (i.e.\ the graph of a bijection
$\pi_e \colon V_{\iota_1(e)} \rightarrow V_{\iota_2(e)}$). 
An automorphism of $\I$ is a
\emph{symmetry} of the $\I$-graph $\H$ if its operation
as a renaming on $\H$ yields an isomorphic $\I$-graph:
$\H^\rho \simeq \H$.
\eD

Clearly the Cayley graph of an $\I$-groupoid is a complete 
$\I$-graph, which also shares every symmetry of the $\I$-groupoid
(cf.\ Definition~\ref{Isymmgroupoiddef} for those symmetries). 
Conversely, any complete $\I$-graph determines an $\I$-groupoid
in the manner indicated for this special case in 
Observation~\ref{groupoidbackfromCayleyobs} above.
In a complete $\I$-graph $\H$ as in Definition~\ref{Igraphdef}, 
the composition of the $\pi_e$ along $w = e_1 \cdots e_n \in
E^\ast[\I,s,t]$, $\pi_w = \prod_{i=1}^n \pi_{e_i} = \pi_{e_n} \circ
\cdots \circ \pi_{e_1}$,
induces a bijection $\pi_w \colon V_s \rightarrow
V_t$, which we denote as $[w]_\H$. The natural composition
operation on matching interface sites induces the structure of an
$\I$-groupoid $\G$ on the set $G = \{ \pi_w \colon w \in E^\ast[\I] \}$.

\bD[$\sym(\H)$]
\label{Isymdef}
\mbox{}\\
From a complete $\I$-graph $\H$ as in Definition~\ref{Igraphdef},
with induced partial bijections $\pi_w$ for $w \in E^\ast[\I]$, we 
obtain the $\I$-groupoid
\[
\sym(\H) := \bigl(G,(G_{s,t}), \,\cdot\,,(1_s), (\pi_e)\bigr)
\]
with $G_{s,t} = \{ \pi_w \colon w \in E^\ast[\I,s,t] \}$, 
composition of partial bijections (which is full composition in
matching sites to match concatenation of labelling sequences)
and identities in corresponding sites as neutral elements.
\eD

In these terms, Observation~\ref{groupoidbackfromCayleyobs} can be restated as
$\sym(\CG) \simeq \G$ if $\G$ is an $\I$-groupoid with Cayley graph $\CG$.

\bD[compatibility]
\label{Igraphcompatdef}
\mbox{}\\
An $\I$-groupoid $\G$ is \emph{compatible} with the complete 
$\I$-graph $\H$ if $\G \succeq \sym(\H)$, i.e.\ if
for all $w \in E^\ast[\I,s,s]$ 
\[
[w]_\G = 1_s\; \Rightarrow \; [w]_\H = \mathrm{id}_{V_s}.
\]
\eD

The following illustrates these concepts and their far-reaching analogy
with the situation for $\E$-groups from Section~\ref{genpatternsec}.

\bO
\label{groupoidcompatobs}
Any $\I$-groupoid $\G$ is compatible with its Cayley graph.
Another $\I$-groupoid $\Ghat$ is compatible with the Cayley graph
$\CG$ of $\G$ if, and only if, $\Ghat \succeq \G$, if, and only if
the map $h \colon \Ghat \rightarrow \G$ which maps 
$[w]_{\Ghat}$ to $[w]_{\G}$ is well-defined (and thus the homomorphism
in question).
\eO

\subsection{Coset acyclicity for groupoids} 
\label{groupoidcosetacycsec} 

Also the following are straightforward analogues of the corresponding notions
for $\E$-groups in Definitions~\ref{cosetcycledef} and~\ref{Nacycdef}.

\bD[coset cycles]
\label{groupoidcosetcycledef}
\mbox{}\\
Let $\G$ be an $\I$-groupoid, $n \geq 2$. 
A \emph{coset cycle of length $n$} in $\G$ is a cyclically
indexed sequence of pointed cosets $(g_i \G[\alpha_i], g_i)_{i \in \Z_n}$
such that, for all $i$,
\bre
\item
(connectivity)
$g_{i+1} \in  g_i \G[\alpha_i]$, i.e.\ $g_i \G[\alpha_i] = g_{i+1} \G[\alpha_i]$; 
\item 
(separation)
$g_i \G[\alpha_{i,i-1}] \cap g_{i+1} \G[\alpha_{i,i+1}] = \emptyset$,
\ere
where $\alpha_{i,j} := \alpha_i \cap \alpha_j$.
\eD

\bD[$N$-acyclicity]
\label{groupoidnacyclicdef}
\mbox{}\\
For $N \geq 2$, an $\I$-groupoid $\G$ 
is \emph{$N$-acyclic} if it admits no 
coset cycles of lengths up to $N$.
\eD

\section{Construction of N-acyclic I-groupoids}
\label{Nacycgroupoidsec}

We associate with a constraint pattern $\I = (S,E)$ for $\I$-groupoids $\G$ 
a set $\Ehat$ of involutive generators and a constraint graph $\Ihat$ 
so that $\I$-groupoids of interest can be identified within suitable 
$\EEhat$-groups $\Ghat$ that are compatible with $\Ihat$.
More specifically, 
we aim for a low-level interpretation of Cayley graphs  of
$\I$-groupoids $\CG$ within the direct product of 
$\Ihat$ with the Cayley graph $\CGhat$ of an $\EEhat$-group $\Ghat$ 
that is compatible with $\Ihat$. Compare Definition~\ref{directproddef} for this
direct product.

\medskip
Firstly, we interpret the directed multi-graph structure of the
constraint pattern
\[
\I = (S,E)
= (S,E,\iota_1,\iota_2,\cdot^{-1})
\] 
in the structure of a constraint graph 
\[
\hat{\I} = (\hat{S}, (R_{\hat{e}})_{\hat{e} \in \hat{E}})
\]
for a set $\Ehat$ of involutive generators. 
Recall that the edge relations $R_e$ of the latter are undirected
while the edges $e \in E$ of the former are directed.
To this end,  associate with every $e \in E$
three new edge labels $\{e\}$, $\{ e,e^{-1} \}$ and $\{ e^{-1} \}$ in $\Ehat$,
as well as $2$ new vertices $s_e$ and $s_{e^{-1}}$ in $\Shat$.
On the basis of 
\[
\barr{rcl}
\Ehat &:= &\bigl\{
\{e\}, \{ e,e^{-1} \}, \{ e^{-1} \}   \colon e \in E \bigr\},
\\
\Shat &:=& S \;\dot{\cup}\; \bigl\{ s_e, s_{e^{-1}} 
 \colon e \in E \bigr\},
\earr
\]
we represent directed $e$-(multi-)edges
as walks of length~$3$ in an $\EEhat$-graph 
$\Ihat$ as follows. 
We replace the directed edge $e \in E[s,s']$
and its inverse $e' := e^{-1} \in E[s',s]$ in $\I$ 
by a succession of $3$ undirected edges with labels 
$\{e\}$, $\{e,e'\}$ and $\{e'\}$ that link $s$ and $s'$ 
via the two new intermediate vertices $s_{e}$ and $s_{e'}$:
\[
\xymatrix{
s \ar@{-}[rr]^{\{e\}}
&& s_{\ssc e} \ar@{-}[rr]^{\{e,e'\}}
&& s_{\ssc e'}
\ar@{-}[rr]^{\{e'\}}
&& s'
}
\]

By the same token, a loop $e\in E[s,s]$ at $s$ and its 
inverse $e' := e^{-1}$ get replaced by a cycle of $3$ undirected edges with labels 
$\{e\}$, $\{e,e'\}$ and $\{e'\}$:
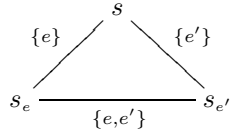
\[
\xymatrix{
& s \ar@{-}[ld]_{\{e\}} \ar@{-}[rd]^{\{e'\}}&
\\
s_{\ssc e} 
\ar@{-}[rr]_{\{e,e'\}}
&& s_{\ssc e'}}
\]

Note that these replacements are inherently symmetric w.r.t.\
edge reversal in the sense that the replacements really concern the edge
pair $\{e,e^{-1}\}$. The direction of $e$ is encoded 
in the directed nature of the walk
\[
s,\{e\},s_{\ssc e}, \{e,e^{-1}\}  ,s_{\ssc e^{-1}},\{e^{-1}\},s',
\] 
whose reversal exactly is the corresponding walk for $e^{-1}$.
The resulting $\EEhat$-graph $\Ihat$ is special also in that each one of
its edge relations $R_{\hat{e}}$ for $\hat{e} \in \Ehat$ consists 
of a single undirected edge. Any automorphism of the constraint
pattern $\I$ turns into a symmetry of the $\EEhat$-graph $\Ihat$, which
is the desired constraint graph.

\medskip
We use this simple schema to associate \emph{$\I$-reachability} w.r.t.\
the constraint pattern $\I$ for $\I$-groupoids and their Cayley graphs
with \emph{$\hat{\I}$-reachability} w.r.t.\ the constraint graph
$\hat{\I}$ for $\hat{\E}$-groups and their Cayley graphs.
Overall, this will allow us to directly extract $\I$-groupoids 
from suitable $\hat{\E}$-groups, in a manner that preserves symmetries
and the desired acyclicity properties.

\medskip
For $\Ehat, \Shat$ and  $\Ihat$ as just constructed from $\I$,
there is a one-to-one correspondence between reduced words in 
\[
\Ehat^\ast[\Ihat,s,t] := 
\{ w \in \Ehat^\ast \colon \mbox{ $w$ labelling a walk from
  $s$  to $t$ in $\Ihat$ }  \} 
\]
and reduced words in $E^\ast[\I,s,t]$ 
that label directed walks from $s$ to $t$ in $\I$.  
In other words, for all $s,t \in S$, 
and modulo passage to reduced words,
the natural replacement map
\[
\barr{@{}rcl@{}}
\hat{\nt\;\nt} \; \colon \; E^\ast[\I,s,t] &\longrightarrow&  \Ehat^\ast[\hat{\I},s,t]
\\
\hnt
w = e_1\cdots e_n &\longmapsto& \hat{w} := 
\{ e_1 \} \{ e_1,e_1^{-1}\}\{ e_1^{-1}\}\cdots \{ e_n \} \{
e_n,e_n^{-1}\}\{ e_n^{-1}\}
\earr
\]
induces a bijection.
For this observation it is essential that reduced words in $\Ehat^\ast[\Ihat]$
can only label walks that link vertices from $S$ if they consist of
concatenations of triplets corresponding to admissible orientations 
of $E$-edges. In connection with the reduced nature of the words involved, 
note on one hand that an immediate concatenation of a triplet for $e \in E$
with the triplet for $e^{-1}$ would not be a reduced $\Ehat$-word. 
On the other hand, the only non-trivial $\{ \{ e \}, \{
e,e^{-1}\},\{e^{-1}\}\}$-component of $\hat{\I}$ consists 
of $\{\iota_1(e),\iota_2(e), s_{\ssc e}, s_{\ssc e^{-1}}\}$. 
The only manner in which a reduced $\Ehat$-word 
can label a walk in $\Ihat$ that exits this $\{ \{ e \},\{
e,e^{-1}\},\{e^{-1}\}\}$-component of $\Ihat$ 
is via $\iota_1(e)$ or $\iota_2(e)$, which are both in $S$.

For notational convenience we also denote as $\hat{\nt\;\nt}$ the
incarnation of the replacement map at the level of reduced words
and at the level of subsets $\alpha \subset E$ that
are closed under edge reversal:
\[
\barr{@{}rcl@{}}
\alpha &\longmapsto& \hat{\alpha} :=
\bigl\{ \{ e \}, \{ e,e^{-1}\}, \{ e^{-1}\} \colon e \in \alpha
\bigr\}.
\earr
\]

\subsection{Groupoids from groups}
\label{extractgroupoidsec}

For a constraint pattern $\I = (S,E)$ and its representation within a
constraint  graph $\Ihat$ for $\EEhat$-graphs according 
to the above translation, consider now an $\EEhat$-group 
$\Ghat$ that is compatible with $\Ihat$. Let $\CGhat$ be the 
Cayley graph of this $\EEhat$-group.
Recall Definition~\ref{directproddef} for the definition of a direct
product, which we now apply to the $\EEhat$-graphs $\Ihat$ and
$\CGhat$:
\[
\Ihat \otimes \CGhat
\]
is an $\EEhat$-graph which reflects $\Ihat$-reachability in the sense 
that $(\hat{s}',\hat{g}')$ is in the $\hat{\alpha}$-connected
component of $(\hat{s},\hat{g})$ if, and only if,  $\hat{g}'$ is in the
$\Ihat[\hat{\alpha},\hat{s}]$-component 
$\CGhat[\Ihat,\alpha,s;\hat{g}]$ 
of $\hat{g}$.

\medskip
We next extract an $\I$-groupoid $\G$
from any $\EEhat$-group $\Ghat$ that is compatible with the
constraint graph $\Ihat$. More specifically, the Cayley graph of the
target $\I$-groupoid 
$\G := \Ghat[\I]$ is 
identified as a definable structure 
within the direct product
$\Ihat \otimes \CGhat$ (technically, a first-order interpretation). The idea is to single out the vertices of
$\Ihat \otimes \CGhat$ with $\Shat$-component in $S \subset \Shat$,
and to replace $\{ e \} \{ e,e^{-1} \}
\{e^{-1}\}$-walks of length $3$ between them by directed
$E$-edges. In essence this is a reversal of the translation that led
from $E$ to $\Ehat$ and from $\E$-graphs to $\EEhat$-graphs.

We define $\G$ in terms of its generators $e \in E$, 
which are interpreted as partial bijections on  
the vertex set of $\Ihat \otimes \CGhat$.
We restrict attention to vertices  
$\{ (s,[\hat{w}]_{\Ghat}) \in \hat{\I} \otimes \Ghat
\colon \hat{w} \in \Ehat^\ast[\Ihat,s,t] \}$ for $s,t \in S \subset \Shat$, and put 
\[
G_{s,t} := \{ (s, [\hat{w}]_{\Ghat}) \colon 
\hat{w} \in \Ehat^\ast[\Ihat,s,t] \} 
= \{ (s,[\hat{w}]_{\Ghat}) \colon w \in E^\ast[\I,s,t] \}.
\]

The second equality appeals to the identification of reduced
words in $\Ehat^\ast[\Ihat,s,t]$ and $E^\ast[\I,s,t]$ for $s,t
\in S \subset \Shat$. The sets $G_{s,t}$ are subsets of the 
vertex set of  $\Ihat \otimes \CGhat$.  
They are disjoint by compatibility of $\Ghat$  with $\Ihat$
and thus partition 
\[
G :=
\dot{\bigcup}_{s,t \in S} G_{s,t}
\]
into subsets (not all necessarily non-empty unless $\I$ is
connected). We write $G_{\ast,t}$ for the union 
$G_{\ast,t} := \bigcup_{s\in S} G_{s,t}$. 
With $e \in E[t,t']$ we associate the 
following partial bijection on the vertex set of $\hat{\I}
\otimes \CGhat$, with domain and image as indicated:
\[
\barr{@{}rcl@{}}
g_e \colon G_{\ast,t} &\longrightarrow& 
G_{\ast,t'}
\\
\hnt
(s,[\hat{w}]_{\hat{\G}}) 
&\longmapsto& (s, [\hat{w}\hat{e}]_{\Ghat}) 
= (s,[\hat{w}]_{\Ghat} \cdot \{ e \} \cdot \{ e,e^{-1} \} \cdot \{ e^{-1} \}),
\earr
\]
where $w \in E^\ast[\I,s,t]$ and $w e \in E^\ast[\I,s,t']$.
Concatenation (and reduction) of corresponding words or walks in $\hat{\I}$ 
induces a well-defined groupoid operation according to
\[
\barr{@{}rcl@{}}
\cdot \; \colon \; G_{s,t} \times G_{t,u} &\longrightarrow& G_{s,u}
\\
\hnt
((s,[\hat{w}_1]_{\Ghat}),(t,[\hat{w}_2]_{\Ghat})) &\longmapsto&
(s, [\hat{w}_1\hat{w}_2]_{\Ghat}),
\earr
\]
where the concatenation relies on the condition that 
$\iota_2(w_1) = t = \iota_1(w_2)$. 
The neutral element in $G_{s,s}$ is $1_s := (s,[\lambda]_{\Ghat}) =
(s, 1_{\Ghat})$.
With these stipulations, 
\[
\G := \Ghat/\I = (G, (G_{s,t})_{s,t \in S}, \cdot, (1_s)_{s
  \in S}, (g_e)_{e \in E})
\]
becomes an $\I$-groupoid with generators 
\[
g_e := [e]_{\G} := 
(\iota_1(e),[\hat{e}]_{\Ghat})
\in G_{\iota_1(e),\iota_2(e)}.
\] 

The induced homomorphism from the free $\I$-groupoid 
(cf.\ Observation~\ref{freeIgroupoidobs})
onto $\G$ maps 
\[
w \in E^\ast[\I,s,t] \longmapsto [w]_\G := (\iota_1(w),[\hat{w}]_{\Ghat}) \in G_{s,t}.
\]

For further analysis we also isolate the induced subgraph 
on those connected components of the $\EEhat$-graph 
$\Ihat \otimes \CGhat$ that embed $\G$: 
\[
\barr{r@{}c@{}l}
\Hhat_0 &\;=\;& (\Ihat \otimes \CGhat) \restr V_0
\,\subset\, \Ihat \otimes \CGhat 
\\
&&
\hnt
\mbox{where }
V_0 := \{ (s,[u]_{\Ghat}) \in \Shat \times \Ghat \colon  
s \in S, u \in \Ehat^\ast[\Ihat,s] \}.
\earr
\]

The set $V_0$ is the vertex set of the union of the connected 
components of the vertices $(s,1)$ in  $\Ihat \otimes \CGhat$
(i.e.\ of the neutral elements $(s,1_s) \in \G$).
Restricting further to vertices of $G \subset V_0$ and 
linking two such vertices by an $e$-edge 
if, and only if, they are 
linked by an $\hat{e} = \{ e \} \{ e,e^{-1}\}\{e^{-1}\}$-labelled walk of 
length $3$ in $\Hhat_0 \subset \Ihat \otimes \CGhat$, we obtain 
an $\I$-graph $\H_0$ that is interpreted in the $\EEhat$-graph
$\Hhat_0$. This $\I$-graph $\H_0$ is (isomorphic to) 
the Cayley graph of the $\I$-groupoid $\G$: 
\[
\barr{r@{}c@{}l}
\H_0 &\;=\;&\bigl( G, (G_{\ast,s})_{s \in S}, (R_e)_{e \in E} \bigr) 
\quad \mbox{ where, for $e = (s,s')$, }
\\
\hnt
R_e &\;=\;& \bigl\{ ((t,[u]_{\Ghat}),(t,[u \{ e \} \{ e,e^{-1}\}\{e^{-1}\} ]_{\Ghat})) \colon
u \in \Ehat^\ast[\Ihat,t,s] \bigr\}.
\earr
\]

\bO
\label{groupoidCayleyobs}
Let $\Ghat$ be an $\EEhat$-group that is compatible with $\Ihat$.
Then the Cayley graph $\CG$ of the $\I$-groupoid $\G = \Ghat/\I$ as just
constructed from $\CGhat$ is (isomorphic to) the 
$\I$-graph $\H_0$ interpreted in $\Hhat_0 \subset \Ihat \otimes \CGhat$. 
\eO

\subsection{Transfer of acyclicity, compatibility and symmetries}
\label{transfersec}

The following is the main technical result of this section. It
reduces the construction of $N$-coset acyclic groupoids to the 
construction of Cayley groups with involutive generators that are
$N$-acyclic over some constraint graph.

\bP
\label{maingroupgroupidprop}
For a constraint pattern $\I$ and its translation into a constraint
graph $\Ihat$ as above, let $\Ghat$ be an $\EEhat$-groupoid 
that is compatible with $\Ihat$. 
Let $\G := \Ghat/\I$ the $\I$-groupoid whose Cayley graph $\CG$ 
is realised as $\H_0$ within $\Ihat \otimes \CGhat$ as discussed above.
\bre
\item
If $\Ghat$ is $N$-acyclic over the constraint graph $\hat{\I}$,  then  
$\G$ is $N$-acyclic.
\item
If $\Ghat$ is compatible with the $\EEhat$-translation of a
complete $\I$-graph $\H$, then $\G$ is compatible with $\H$.
\item
Any symmetry $\rho$ 
of $\I$ induces a permutation $\hat{\rho} \in \mathrm{Sym}(\Ehat)$
that is a symmetry of $\Ihat$; if $\hat{\rho}$ is a symmetry of
$\Ghat$ then $\rho$ is a symmetry of $\G$.
\ere
\eP

The main claim, concerning $N$-acyclicity, 
follows directly from the compatibility
of the corresponding notions of cycles with the interpretation of $\CG
\simeq \H_0$ in $\Hhat_0 \subset \Ihat \otimes \Ghat$. This 
is expressed in the following lemma; the arguments
towards compatibility with a given $\H$ and compatibility of the whole
construction with symmetries are straightforward.

\bL
\label{cyclestocycleslem}
In the situation of Proposition~\ref{maingroupgroupidprop}
there is a natural translation of coset cycles in the $\I$-groupoid $\G
=\Ghat/\I$ [based on the map $\hat{\nt\;\nt}$ for generator sets],
which translates coset cycles in the groupoid $\G$
into $\Ihat$-coset cycles of the same length in the group $\Ghat$.
It follows that $\G$ is $N$-acyclic if 
$\Ghat$ is $N$-acyclic over $\hat{\I}$.
\eL

\prf
Let 
\[
(\ast) \qquad
(g_i\G[\alpha_i],g_i)_{i \in \Z_n}
\]
be a coset cycle in the groupoid $\G$, 
according to Definition~\ref{groupoidcosetcycledef}, 
viewed in $\H_0$.
The connectivity condition for the cycle $(\ast)$ and the 
manner in which $\H_0$ is interpreted in $\Hhat_0 \subset \Ihat \otimes
\CGhat$ implies that there is an 
$\hat{\alpha}_i$-walk from $\hat{g}_i = [\hat{w}_i]_{\Ghat}$
to $\hat{g}_{i+1}  = [\hat{w}_{i+1}]_{\Ghat}$, labelled 
by the $\hat{\nt\;\nt}$-translation of an 
$\alpha_i$-word of generators representing $g_{i}^{-1} g_{i+1} \in \G$. 
The natural $\hat{\nt\;\nt}$-translation of the cycle $(\ast)$ into 
$\Ghat$ is
\[
(\ast\ast) \qquad (\Ghat[\Ihat,\hat{\alpha}_i, s_i; \hat{g}_i],\hat{g}_i)_{i \in
  \Z_n}, 
\]
where the labels $s_i \in S \subset \Shat$ 
are determined by the sorts of the 
$g_i$ according to $s_i = \iota_2(w_i)$.  
This translation in effect replaces the subsets $g_i\G[\alpha_i]$ by their closures 
$\Ghat[\hat{\I},\hat{\alpha}_i, s_i;\hat{g}_i]$ w.r.t.\ 
$\Ihat$-reachability inside their $\hat{\alpha}_i$-coset.
This closure is obtained as the union of all 
$\{\{e\},\{e,e^{-1}\},\{e^{-1}\}\}\}$-connected components
in $\Hhat_0 \subset \Ihat \otimes \CGhat$ 
that contain at least one element of the representation of 
$g_i\G[\alpha_i]$ in $\H_0$.

It is clear that $(\ast\ast)$ has the format of a potential
$\Ihat$-coset cycle of length $n$ over $\Ihat$ in $\Ghat$ 
in the sense of Definition~\ref{Icosetcycledef}. 
It remains to show that the separation
condition in Definition~\ref{Icosetcycledef} for
$(\ast\ast)$ follows from the analogous condition in
Definition~\ref{groupoidcosetcycledef} for $(\ast)$.

Suppose that, in violation of the separation condition
for $(\ast\ast)$, 
\[
(\dagger) \qquad
\hat{g} \in \hat{g}_i\Ghat[\hat{\I},\hat{\alpha}_{i-1,i}, s_i] \cap 
\hat{g}_{i+1}\Ghat[\Ihat,\hat{\alpha}_{i,i+1}, s_{i+1}],
\]
where $\alpha_{i,i\pm 1} = \alpha_i \cap \alpha_{i\pm 1}$.  

We analyse this non-trivial intersection in terms of the 
representation of $\G$ in $\Hhat_0 \subset \Ihat \otimes \CGhat$.
Let $\hat{g} = \hat{g}_i [w_{i,i-1}]_{\Ghat}  = \hat{g}_{i+1}
[w_{i,i+1}]_{\Ghat}$ for suitable $w_{i,j} \in \hat{\alpha}_{i,j}^\ast[\Ihat,s_j]$.
By the separation condition for $(\ast)$, $\hat{g}$ is not 
represented as an element of $\G$ or
$\H_0$ in $\Hhat_0$, so that $\iota_2(w_{i,j}) \in \Shat \setminus S$, 
i.e.\ $\iota_2(w_{i,j}) \in \{ s_{\ssc e},s_{\ssc e^{-1}} \}$ 
for some $e \in E$. 

But in $\{\{e\},\{e,e^{-1}\},\{e^{-1}\}\}$-components 
of elements of $G$ or vertices of $\H_0$ in $\Hhat_0$,
any vertex with $\iota_2$-value outside $S$ is isolated 
from all vertices with $\iota_2$-value in $S$ by 
$\{e\}$- and $\{e^{-1}\}$-edges (just as vertices in $\Shat\setminus S$ 
are isolated from $S$ in $\Ihat$). 
So $(\dagger)$ implies that 
$e,e^{-1} \in \alpha_{i,j}$ for $j = i\pm 1$.
This would imply that there also is an $e$-link between 
the elements of that component that represent elements of $G$.
So elements of $g_i\G[\alpha_{i,i-1}]$ and of
$g_{i+1}\G[\alpha_{i,i+1}]$ occur in the same $\{ e, e^{-1}\}$-component,
which would violate the separation condition for $(\ast)$ since 
$e \in \alpha_{i,i-1}$ and $e \in \alpha_{i,i+1}$.
\eprf

\bT
\label{acycgroupoidthm}
For every finite constraint pattern $\I = (S,E)$, every 
complete $\I$-graph $\H$ and  every $N \geq 2$ 
there is a finite $N$-acyclic $\I$-groupoid 
$\G$ that is compatible with $\H$. Such $\G$ can be chosen to be 
fully symmetric w.r.t.\ the given data, i.e.\ such that every 
symmetry $\rho$ of $\I$ that induces a symmetry of the $\I$-graph 
$\H$ is also a symmetry of the $\I$-groupoid $\G$:
$\H^\rho \simeq \H \Rightarrow \G^\rho \simeq \G$.
\eT

Choosing the Cayley graph of a given $\I$-groupoid $\G_0$ for $\H$,
we obtain a fully symmetric $N$-acyclic $\I$-groupoid $\G \succeq
\G_0$. For $\H := \sym(\I)$ (regarding $\I$ as a complete $\I$-graph
according to Definition~\ref{Igraphdef}) one obtains $N$-acyclic
$\I$-groupoids that are fully symmetric over $\I$.

\section{Conclusion and primary applications}

The generic constructions of the preceding chapters show
the versatility of the fruitful idea to go back and forth
between group-like structures (monoids and groups as well as groupoids)
and graph-like structures (graphs and multi-graphs, undirected as well
as directed, and possibly vertex- or edge-coloured). In one
direction the passage involves the familiar encoding of algebraic
structures in the graph-like representation of generators, as in the
classical notion of Cayley graphs for groups; in the converse
direction, permutation groups are induced by various operations
on graph-like structures. We have here tried to contribute to these
connections with a special emphasis on strong algebraic-combinatorial 
criteria of graded acyclicity in finite structures. The constructions 
presented here extend techniques for the construction
of $N$-acyclic groups with involutive generators 
from~\cite{Otto12JACM} to yield a conceptual improvement and correction
of the proposed constructions for groupoids from~\cite{lics2013}. 
Due to its symmetry preserving generic character, 
the new presentation also supports the use of these groupoids 
in~\cite{arXiv2018} where symmetry considerations are of the essence
towards lifting local symmetries to global symmetries in finite
structures. 
In a different direction, $2$-acyclic finite groupoids have
also been used to resolve an open problem of a purely 
semigroup-theoretic nature in Bitterlich~\cite{Bitterlich}, 
which has meanwhile been treated comprehensively in~\cite{ABO}. 

To conclude the present treatment we briefly look at the most 
salient application for finite groups and groupoids of 
graded coset-acyclicity. This concerns the construction of 
finite coverings of graphs and hypergraphs that unravel short
cycles. 
\bne
\item
Natural, unbranched finite coverings of graphs 
by graphs with interesting acyclicity properties can be obtained as weak
subgraphs of the Cayley graphs of suitable $\E$-groups where $E$ 
is the set of edges of the graph to be covered (individually labelled
as it were). While similar constructions have been used 
in~\cite{Otto04,Otto12JACM} 
and a precursor for special graphs in~\cite{CO17}, we illustrate the key to the new 
generalisation in Proposition~\ref{graphcoverprop} below. 
\item
Natural reduced products with $N$-acyclic $\I$-groupoids 
yield finite branched $N$-acyclic coverings of hypergraphs $\V= (V,S)$
where 
a constraint pattern $\I = (S,E)$ is induced by the intersection graph
of $\V$ that encodes the intersection 
pattern between hyperedges in the given hypergraph.
\item
A new and more direct approach to finite branched  $N$-acyclic coverings of
hypergraphs $\V = (V,S)$ can be based on $\I$-products between a constraint graph 
$\I = (S,E)$ induced by the intersection graph of $\V$ and suitable
$\E$-groups that are not just $N$-acyclic but $N$-acyclic over $\I$; 
cf.~Proposition~\ref{hypcoverprop} below.
\ene

Of these fundamental applications, (2)
has been explored in stages in~\cite{Otto12JACM,lics2013,arXiv2018}.
Application~(1) is new in its form that involves the new notion of $N$-acyclicity
of groups over a constraint graph $\I$, rather than working with weak substructures
of direct products with (unconstrained) $N$-acyclic Cayley groups.
Application~(3) similarly
supersedes~(2). 
Recall from Section~\ref{constraintsec} how control of cyclic 
configurations can be extended 
to configurations governed by reachability patterns w.r.t.\ a given
constraint graph $\I$. While we have seen in Section~\ref{groupoidsec} 
how such groups can yield coset acyclicity in groupoids as used in~(2), 
the underlying groups can also be put to use directly in~(1) and~(3).

\paragraph*{Graph coverings.}
For a finite simple graph $\V = (V,E)$ consider, as a set $E$ of
involutive generators for $\E$-groups, the set of all edges  
$e=\{v,v'\} \in E$, and as a constraint graph $\I$ the $\E$-graph 
$\I = (V, (\{ e \})_{e \in E})$ ($\V$ with individually labelled
edges). For any $\E$-group $\G$ that is compatible with $\I$  
consider the direct product $\Vhat = \I \otimes \CG$
of the constraint graph $\I$ with the Cayley graph $\CG$ of $\G$
according to Definition~\ref{directproddef}. Then the 
natural projection
\[
\pi \colon \barr[t]{@{}rcl@{}} 
\Vhat  &\longrightarrow& \V
\\
\hnt
(v, g) &\longmapsto& v 
\earr 
\]
provides an unbranched covering of $\V$ by 
$\hat{\V} =\I \otimes \CG$. Recall that the connected components of 
$\I \otimes \CG$ are isomorphic to weak subgraphs of $\CG$ 
(cf.\ Observation~\ref{directprodobs}). 

\bP
\label{graphcoverprop}
Let $\V = (V,E)$ be a connected finite simple graph, $E$ associated
with its edge set $E$ as above and $\G$ an $\E$-group that is 
compatible with the $\E$-graph $\I := (V,(\{e\})_{e \in E})$.  
Then each connected component $\H$ of the direct product 
$\I \otimes \CG$ 
\bre
\item
is realised as a weak subgraph of the Cayley graph $\CG$ of $\G$ and
\item
is an 
unbranched finite covering w.r.t.\ the natural projection 
$\pi \colon (v,g) \mapsto v$.
\ere
This covering graph $\H$ inherits the acyclicity properties of $\CG$: 
if $\G$ is $N$-acyclic over $\I$, 
then $\H$ admits no cyclic configurations of length up to~$N$ 
of overlapping $\alpha_i$-connected components with the natural
separation condition for subsets $\alpha_i \subset E$.
\eP

\paragraph*{Hypergraph coverings.}
With a finite hypergraph $\V = (V, S)$ with $S \subset
\mathcal{P}(V)\setminus \{ \emptyset \}$
associate its \emph{intersection graph}
$\I = (S,E)$ where 
\[
E = \{ \{s,s'\} \in S^2 \colon  s \not= s'  , s \cap s' \not=\emptyset\}.
\]

If $\G$ is an $\E$-group that is compatible with $\I$,
then the direct product $\I \otimes \CG$ of the intersection graph $\I$
with the Cayley graph $\CG$ of $\G$ (cf.\ Definition~\ref{directproddef})
gives rise to a finite branched hypergraph covering
of $\V = (V,S)$ by a hypergraph $\Vhat = (\hat{V},\hat{S})$ as follows. 
Consider the disjoint union of $\I \otimes \CG$-tagged copies
of the hyperedges of $\V$,
\[
  \bigcup_{(s,g) \in \I \otimes \CG} \!\!\! s \times \{ (s,g) \} 
  =
  \bigl\{ (v,s,g) \in V \times S \times \G \colon
  v \in s 
  \bigr\}
\]
and its quotient w.r.t.\ the equivalence relation $\approx$ 
induced by identifications 
\[
(v,s,g) \approx (v,s',ge) \; \mbox{ for } \; 
e = \{s,s'\} \in E \mbox{ if } v \in s \cap s'. 
\]

The induced equivalence is such that $(v,s,g) \approx (v',s',g')$
if, and only if, $v'= v  \in s \cap s'$ and $g' = g h$ for some
$h = [w]_\G$ with $w  \in  \alpha_v^\ast[\I,s,s']$,
$\alpha_v := \{ e = \{t,t'\} \in E \colon v\in t \cap t'  \}$,
which implies $g' \in \CG[\I,\alpha_v,s;g]$.
We write $[v,s,g]$ for the $\approx$-equivalence class of $(v,s,g) \in s \times\{
(s,g) \}$. These equivalence classes form the vertices of
$\hat{\V}$; its hyperedges the subsets induced by the $s \in S$,
denoted as
\[
[s,g] := \{ [v,s,g] \colon v \in s \} \; \mbox{ for } \; (s,g) \in \I \otimes \CG.
\]

In the $\approx$-quotient, the $e$-edge between 
$(s,g)$ and $(s',ge)$ in $\I \otimes \CG$ for $e = \{s,s'\}$ becomes 
an intersection of the copies $[s,g]$ and $[s',ge]$ of corresponding 
hyperedges above $s$ and $s'$ in the covering hypergraph.
This covering hypergraph, which we denote as $\hat{\V} := \V \otimes_\I \CG$
is
\[
  \hat{\V} = (\hat{V},\hat{S}) =: \V \otimes_\I \CG
  \; \mbox{ where }
  \left\{ \;
\barr{@{}l@{\;:=\;}l@{}}
\hat{V} & \{ [v,s,g] \colon v \in s, (s,g) \in \I \otimes \CG \}
\\
\hnt
\hat{S} & \{ [s,g] \colon (s,g) \in \I \times \CG \}
\earr
\right.
\]
with covering projection 
\[
\barr{@{}rcl@{}}
\pi \colon \hat{\V} = (\hat{V},\hat{S}) &\longrightarrow& 
\V = (V,S)
\\
\hnt
[v,s,g] &\longmapsto& v.
\earr
\]

\bO
\label{hypcovobs}
Let $\hat{\V} := \V \otimes_\I \CG$ be
defined as above for an $\E$-group $\G$ that is
compatible with the intersection graph $\I = (S,E)$ of $\V = (V,S)$.
Then hyperedges $[s_1,g_1]$ and $[s_2,g_2]$ of $\hat{\V}$ intersect precisely in
vertices $[v,s_1,g_1] = [v,s_2,g_2]$ with $v \in s_1 \cap s_2$ where 
$g_2 \in \CG[\I,\alpha_v,s_1;g_1]$,  $\alpha_v := \{ e=\{s,s'\} \in E \colon v \in s \cap
s' \}$ 
(i.e.\ $g_2 = g_1 [w]_\G$ for some $w$ labelling an $\I$-walk from
$s_1$ to $s_2$ in $\I$, or from
$(s_1,g_1)$ to $(s_2,g_2)$ in $\I \otimes \CG$: 
$w \in \alpha_v^\ast[\I,s_1,s_2]$, cf.\
Definition~\ref{Ialphawalksdef} and Observation~\ref{directprodobs}).
\eO

\bP
\label{hypcoverprop}
Let $(V,S)$ be a finite hypergraph, 
$\I = (S,E)$ its intersection graph, 
$\G$ an $\E$-group that is compatible with $\I$. 
Then the hypergraph $\hat{\V} := \V \otimes_\I \CG$ as defined above
gives rise to a finite branched hypergraph covering
$\pi \colon \hat{\V} \longrightarrow \V$.
\\
The 
covering hypergraph $\hat{\V}$
inherits the acyclicity properties of $\CG$
in the following sense: if $\G$ is $N$-acyclic over $\I$, 
then every induced sub-hypergraph on up to $N$ vertices is acyclic 
in the sense of classical hypergraph theory.  
\eP 

For acyclicity in 
hypergraph terminology (conformality and chordality and
tree-decomposability), 
compare~\cite{BeeriFaginetal,Berge}.

\prf
Consider the hypergraph $\V \otimes_\I \CG$, as defined above, for 
an $\E$-group $\G$ that is compatible with the intersection graph $\I
= (S,E)$ of $\V$ and $N$-acyclic for some $N \geq 3$. 
The covering property follows from Observation~\ref{hypcovobs}, and 
it remains to argue for $N$-acyclicity of $\hat{\V}$. 
We show that the Gaifman graph of 
$\hat{\V}$ cannot have chordless
cycles of lengths $n$ for $3 < n \leq N$
($N$-chordality), nor can it have cliques of size up to $N$ that are
not contained in a single hyperedge ($N$-conformality).

\medskip\noindent
\emph{$N$-chordality.} Suppose $(\hat{v}_i)_{i \in \Z_n}$ is a
chordless cycle of length $n > 3$ in the Gaifman graph 
of $\hat{\V} = (\hat{V},\hat{S})$,
and let hyperedges $[s_i,g_i] \in \hat{S}$ with $\hat{v}_{i} \in
[s_i,g_i] \cap [s_{i+1},g_{i+1}]$ witness the connectivity condition in this cycle.
This implies that $\hat{v}_{i}$ can be represented as   
$\hat{v}_{i} = [v_i, s_i ,g_i] = [v_i,s_{i+1},g_{i+1}]$ 
for some $v_{i} \in s_i \cap s_{i+1}$ and that
$g_{i+1} = g_i [w_i]_\G$ for some $w_i \in
\alpha_i^\ast[\I,s_i,s_{i+1}]$ where $\alpha_i = \{ e = \{s,s'\} \in E \colon 
v_i \in s \cap s'  \}$. 
We claim that 
\[
(\CG[\I,\alpha_i,s_i;g_i],g_i)_{i \in   \Z_n}
\] 
is an $\I$-coset cycle in $\G$, in the sense of 
Definition~\ref{Icosetcycledef}. Then $n > N$ follows
from $N$-acyclicity of $\G$ over $\I$. 
Of the two conditions in Definition~\ref{Icosetcycledef},
connectivity is obvious; it remains to check the 
separation condition:
\[
\CG[\I,\alpha_{i,i-1},s_i;g_i] \cap \CG[\I,\alpha_{i,i+1},s_{i+1};g_{i+1}] = \emptyset,
\]
where $\alpha_{i,j} := \alpha_i \cap \alpha_j$.
But this follows from chordlessness of the given cycle. Suppose $g$ were a
member of this intersection, i.e.\ 
$g =g_i [w]_\G$ for some 
$w \in \alpha_{i,i-1}^\ast[\I,s_i,s]$ and 
$g := g_{i+1}[w']_\G$ for some 
$w' \in \alpha_{i,i+1}^\ast[\I,s_{i+1},s]$ (the same $s$,
due to compatibility of $\G$ with $\I$). 
Then $\hat{v}_{i-1} = [v_{i-1},s_{i-1},g_{i-1}] = [v_{i-1},s,g]$ because 
$w \in \alpha_{i-1}^\ast$ and $\hat{v}_{i-1} \in [s_i,g_i]$, which implies 
$\hat{v}_{i-1} \in [s,g]$. Similarly,
$\hat{v}_{i+1} = [v_{i+1},s_{i+1},g_{i+1}] = [v_{i+1},s,g]$ because
$w' \in \alpha_{i+1}^\ast$, which implies that  
$\hat{v}_{i+1} \in [s,g]$, too. So the given cycle would have a chord
linking $\hat{v}_{i-1}$ to $\hat{v}_{i+1}$.

\medskip\noindent
\emph{$N$-conformality.}
Suppose $m = \{ \hat{v}_i \colon 1 \leq i \in \Z_n \}$ forms a
clique of size~$n \geq 3$ in the Gaifman graph of $\hat{\V} = (\hat{V},\hat{S})$ 
such that every subset
$m_i := m \setminus \{ \hat{v}_{i} \}$ of size $n-1$ is contained
in some hyperedge (a minimal violation of conformality). 
For $i \in \Z_n$, let $[s_i,g_{i}] \in \hat{S}$ be a hyperedge
that contains $m_i = m \setminus \{\hat{v}_{i}\} =
\{  \hat{v}_{j} \colon j \not= i \}$.
So, for $j \not= i$, $\hat{v}_j \in [s_i,g_{i}]$ can be represented as  
$\hat{v}_j = [v_j,s_i,g_i] \in [s_i,g_i]$.
Let $\gamma_i := \alpha_{v_i} = \{ e = \{s,s'\} \in E \colon v_i \in s \cap s'\}$; 
and $\gamma_{i,i'} :=\bigcap_{j \not= i,i'}\gamma_j$, so that $\I$-walks
in $\gamma_{i,i'}$ preserve all $\hat{v}_j$ for $j \not= i,i'$,
in the sense that $[v_j,s_j,g_j] = [v_j,s,g_j h]$ for $h = [w]_\G$,
$w \in \gamma_{i,i'}^\ast[\I,s_j,s]$.
In particular, $g_{i+1} \in  \CG[\I,\alpha_{i},s_i;g_i]$
for $\alpha_i := \gamma_{i,i+1}$,
as $\hat{v}_j \in [s_i,g_{i}] \cap [s_{i+1},g_{i+1}]$ for $j \not= i,i+1$, 
cf.\ Observation~\ref{hypcovobs}.
Consider then 
\[
  (\CG[\I,\alpha_i,s_i;g_i],g_i)_{i \in   \Z_n} \; \mbox{ for } \;
  \alpha_i := \gamma_{i,i+1} = \!\!\!\bigcap_{j \not= i,i+1}\!\!\!\gamma_j
\]
as a candidate for an $\I$-coset cycle. We show that if this were not
an $\I$-coset cycle, then the whole of $m$ would be contained in some
hyperedge $[s,g] \in \hat{S}$. 
The connectivity condition on $\I$-coset cycles 
from Definition~\ref{Icosetcycledef} is just that
$g_{i+1} \in  \CG[\I,\alpha_{i},s_i;g_i]$.
The separation condition now is that 
\[
\CG[\I,\alpha_{i,i-1},s_i; g_i] \cap \CG[\I,s_{i+1},\alpha_{i,i+1},s_{i+1};g_{i+1}] = \emptyset,
\]
where $\alpha_{i,i-1} = \alpha_{i} \cap \alpha_{i-1}
=\bigcap_{j \not= i} \alpha_j
= \{ e = \{s,s'\} \in E \colon m_i \subset s \cap s' \}$ and
$\alpha_{i,i+1} = \alpha_{i} \cap \alpha_{i+1}
=\bigcap_{j \not= i+1} \alpha_j
= \{ e = \{s,s'\} \in E \colon m_{i+1} \subset s \cap s' \}$.

Assume there were some $g$ in the above intersection, i.e.\
$g = g_i h$ for some $h = [w]_\G$ with $w \in \alpha_{i,i-1}^\ast[\I,s_i,s]$ 
and $g =  g_{i+1} h'$ for some $h' =[w']_\G$ with
$w' \in \alpha_{i,i+1}^\ast[\I,s_{i+1},s]$, for the same $s \in S$
(due to compatibility of $\G$ with $\I$).
We claim that this would imply $m \subset [s,g]$, contradicting the assumptions.
For $j \not= i$, $\hat{v}_j \in [s_i,g_i]$ implies 
$\hat{v}_j = [v_j,s_i,g_i]$, which, together with $h = [w]_\G$,
$w \in \alpha_{i,i-1}^\ast[\I,s_i,s]$ entails that $\hat{v}_j \in [s,g]$ for $j \not= i$.
Similarly, for $j \not= i+1$, $\hat{v}_j \in [s_{i+1},g_{i+1}]$ implies 
$\hat{v}_j = [v_j,g_{i+1}]$, which, together with $h' = [w']_\G$,
$w' \in \alpha_{i,i+1}^\ast[\I,s_{i+1},s]$
entails that $\hat{v}_j \in [s,g]$ for $j \not= i+1$. As $\G$ is $N$-acyclic,
$n > N$ follows.  
\eprf

\paragraph*{Acknowledgements.}
This paper has a complicated history. It is based on 
ideas revolving about acyclicity in hypergraphs and Cayley groups
that I first expounded in \cite{Otto12JACM} (the journal version of 
a LICS 2010 paper). That paper deals with the construction of finite 
$N$-acyclic groups  and applies it to the finite model theory of 
the guarded fragment. The promising extension to the groupoid
situation was seemingly achieved with~\cite{lics2013}, 
which also paved the way to more fundamental
applications in hypergraph coverings and extension properties for
partial automorphisms. These applications were successively elaborated 
in a series of arXiv preprints leading to~\cite{arXiv2018} and stand
as key motivations and achievements also of my DFG project
on \emph{Constructions and Analysis in Hypergraphs of Controlled
  Acyclicity}, 2013--18. The proposed construction of $N$-acyclic 
groupoids in~\cite{lics2013}, however, contained a serious flaw that
was also carried in the arXiv preprints. 
This mistake from~\cite{lics2013} was finally found out in 2019 by Julian
Bitterlich \emph{after}
he had shown, as part of his PhD work, that those results would 
positively resolve a longstanding conjecture by Henkell and Rhodes 
in semigroup theory, which had attracted the attention of, 
among others, Karl Auinger. I am deeply indebted to Julian
Bitterlich for identifying and clarifying the gap in my earlier
attempts at the construction of finite $N$-acyclic groupoids.
Initial concerns about the use of my sketchy and, as it later turned
out, truly flawed construction from~\cite{lics2013} in Julian 
Bitterlich's contributions had earlier been raised by Jiri~Kadourek. In a
way, it was his initial challenge of the result that triggered the 
intense re-investigation that led to Julian's discovery of the actual 
flaw. I am also very grateful to Karl Auinger and Julian Bitterlich 
for helpful discussions in these difficult matters at the time 
of Julian's breakthrough and discovery of my mistake. Karl also 
provided extremely valuable critical comments on earlier versions of
this paper. Those helped to eliminate some gaps and mistakes
besides triggering a more specific treatment towards the Henckell--Rhodes
conjecture in~\cite{ABO}, which also recasts key concepts
from the current paper in a more algebraic framework.

\end{document}